\documentclass[11pt]{amsart}

\usepackage{amscd}
\usepackage{amssymb}
\usepackage{fancyheadings}
\usepackage{enumerate}

\theoremstyle{plain} \newtheorem{theorem}{Theorem}[subsection]
\theoremstyle{plain} 
\theoremstyle{plain} \newtheorem{lemma}[theorem]{Lemma}
\theoremstyle{plain} \newtheorem{proposition}[theorem]{Proposition}

\newtheorem{corollary}[theorem]{Corollary}

\newcommand{\nr}{\addtocounter{theorem}{1}  
                   \noindent {\thetheorem .}}
\newcommand{\defi}{\medskip \noindent {\it Definition \nr} }
\newcommand{\defifin}{\medskip}
\newcommand{\eks}{\medskip \noindent {\it Example \nr} }
\newcommand{\eksfin}{\medskip}
\newcommand{\rem}{\medskip \noindent {\it Remark \nr} }
\newcommand{\remark}{\rem}
\newcommand{\remfin}{\medskip}
\newcommand{\remarkfin}{\remfin}

\newcommand{\nota}{\medskip \noindent {\it Notation \nr} }
\newcommand{\notafin}{\medskip}
\newcommand{\llabel}{\addtocounter{theorem}{-1}
\refstepcounter{theorem} \label}

\newcommand{\pv}{{\bf P}(W)}
\newcommand{\pw}{{\bf P}(U)}
\newcommand{\pvw}{{\bf P}(W/U)}
\newcommand{\pvv}{{\bf P}(W/(w))}
\newcommand{\pvwo}{{\bf P}(W/U_0)}
\newcommand{\op}{\mathcal U}

\newcommand{\Bc}{{\bf B}}
\newcommand{\Ac}{{\bf A}}

\newcommand{\coh}{{{\text{{\rm coh}}/\pv}}}
\newcommand{\cohvw}{{{\text{{\rm coh}}/\pvw}}}
\newcommand{\qc}{{{\text{{\rm qc}}/\pv}}}
\newcommand{\qco}{{\text{{\rm qc}}/\op}}
\newcommand{\qcw}{{\text{{\rm qc}}/\pw}}
\newcommand{\qcvv}{{{\text{{\rm qc}}/\pvv}}}
\newcommand{\vb}{{{\text{{\rm vb}}/\pv}}}
\newcommand{\vbvw}{{{\text{{\rm vb}}/\pvw}}}
\newcommand{\gco}{{G{\text{{\rm -coh}}/\pv}}}
\newcommand{\gqc}{{G{\text{{\rm -qc}}/\pv}}}
\newcommand{\gvb}{{G{\text{{\rm -vb}}/\pv}}}

\newcommand{\dbvb}{D^b(\vb)}
\newcommand{\dbco}{D^b(\coh)}
\newcommand{\dbcq}{D^b_{\text{coh}}(\qc)}
\newcommand{\dcq}{D_{b, \text{coh}}(\qc)}
\newcommand{\dcqvv}{D_{b, \text{coh}}(\qcvv)}
\newcommand{\dbcovw}{D^b(\cohvw)}
\newcommand{\dbvbvw}{D^b(\vbvw)}
\newcommand{\dbqc}{D^b(\qc)}
\newcommand{\dbqco}{D^b({\rm qc}/\op)}
\newcommand{\dbcow}{D^b({\rm coh}/\pw)}
\newcommand{\dbqcw}{D^b({\rm qc}/\pw)}
\newcommand{\dbcqw}{D^b_{\text{coh}}({\rm qc}/\pw)}
\newcommand{\dbcoM}{D^b_{-}(\coh)}
\newcommand{\dbqcM}{D^b_{-}(\qc)}

\newcommand{\dbqcoM}{D^b_{-}({\rm qc}/\op)}

\newcommand{\modv}[1]{{#1}\text{-{mod}}}

\newcommand{\vmod}[1]{{#1}\text{-{mod}}}
\newcommand{\vfmod}[1]{{#1}\text{-{fmod}}}
\newcommand{\modstab}[1]{{#1}{\text{-\underline{fmod}}}}
\newcommand{\svgmod}{\sv,G\,\text{-module}}

\newcommand{\til}{\sim}

\newcommand{\La}{{E}}
\newcommand{\sv}{{S(W)}}
\newcommand{\sw}{{S(U)}}
\newcommand{\lv}{{\omega_E}}
\newcommand{\lw}{{\omega_{E(U^*)}}}
\newcommand{\lvd}{{\La(V)}}
\newcommand{\svw}{{S(W/U)}}
\newcommand{\svv}{{S(W/(w))}}
\newcommand{\lwd}{{\La(U^*)}}
\newcommand{\lWd}{{\La(W^*)}}
\newcommand{\lvwd}{{\La((W/U)^*)}}
\newcommand{\lvvd}{{\La((W/(w))^*)}}

\newcommand{\gd}{\circledast}
\newcommand{\ad}{{A^{!}}}

\newcommand{\bgd}{{B^{\gd}}}

\newcommand{\gopv}{{\mathcal O}_{\pv}}

\newcommand{\gopvv}{{\mathcal O}_{\pvv}}
\newcommand{\gopvp}{{\mathcal O}_{\pv, P}}
\newcommand{\gF}{{\mathcal F}}
\newcommand{\gG}{{\mathcal G}}
\newcommand{\gQ}{{\mathcal Q}}
\newcommand{\gE}{{\mathcal E}}
\newcommand{\gI}{{\mathcal I}}
\newcommand{\gJ}{{\mathcal J}}
\newcommand{\gR}{{\mathcal R}}
\newcommand{\gK}{{\mathcal K}}
\newcommand{\gB}{{\mathcal B}}
\newcommand{\gA}{{\mathcal A}}
\newcommand{\gC}{{\mathcal C}}
\newcommand{\gT}{{\mathcal T}}
\newcommand{\gO}{{\mathcal O}}

\newcommand{\fr}{{-F}}
\newcommand{\cfr}{{-cF}}

\newcommand{\acf}{A\!\cfr}

\newcommand{\lvdcf}{\lvd\!\cfr}
\newcommand{\lwdcf}{\lwd\!\cfr}
\newcommand{\svf}{\sv\!\fr}
\newcommand{\svwf}{\svw\!\fr}
\newcommand{\lvwdcf}{\lvwd\!\cfr}
\newcommand{\svgf}{\sv,G \fr}
\newcommand{\lvdgcf}{\lvd,G\!\cfr}

\newcommand{\Kom}{\text{Kom}}

\newcommand{\koml}{\text{Kom}(\lvd)}
\newcommand{\koms}{\text{Kom}(\sv)}
\newcommand{\KomLin}{\text{KomLin}}
\newcommand{\klvdcf}{K(\lvdcf) }
\newcommand{\kolvdcf}{K^{\circ}(\lvdcf)}
\newcommand{\kolwdcf}{K^{\circ}(\lwdcf)}
\newcommand{\kolvdcfMw}{K^{\circ}_{-}(\lvdcf)}
\newcommand{\kolvwdcf}{K^{\circ}(\lvwdcf)}
\newcommand{\kocolvdcf}{K^{\circ}_{\text{coh}}(\lvdcf)}
\newcommand{\komolvdcf}{\Kom^{\circ}(\lvdcf)}
\newcommand{\ksv}{K(\sv)}
\newcommand{\koacf}{K^{\circ}(\acf)}

\newcommand{\res}{\text{res}}
\newcommand{\Hom}{\text{Hom}}

\newcommand{\gext}{\gE xt}
\newcommand{\id}{\text{{id}}}
\newcommand{\im}{\text{im}\,}
\newcommand{\codim} {\text{codim}\,}
\newcommand{\resol}{\text{resol}\,}
\newcommand{\rank}{\text{rank}\,}
\newcommand{\lpd}{\text{lpd}\,}
\newcommand{\coker}{\text{coker}\,}
\newcommand{\Tor}{\text{Tor}}

\newcommand{\gs}{G_{\sv}}
\newcommand{\gsw}{G_{\sw}}
\newcommand{\fl}{F_{\lvd}}
\newcommand{\tilfl}{{\overset{\til}{F}}_{\lvd}}
\newcommand{\tilflvw}{{\overset{\til}{F}}_{\lvwd}}
\newcommand{\tilov}[1]{{\overset{\til}{#1}}}

\newcommand{\tilflmin}{{\overset{\til}{F}}_{\lvd,min}}
\newcommand{\gsmin}{G_{\sv, min}}
\newcommand{\gswmin}{G_{\svw, min}}
\newcommand{\flmin}{F_{\lvd, min}}
\newcommand{\gsvw}{G_{\svw}}
\newcommand{\gsvv}{G_{\svv}}

\newcommand{\flvw}{F_{\lvwd}}

\newcommand{\lder}{{\bf L}}
\newcommand{\lgvw}{\lder \gG_{|\pvw}}
\newcommand{\lfvw}{\lder \gF_{|\pvw}}

\newcommand{\ytrev}{T}
\newcommand{\ytrevw}{T^\prime}
\newcommand{\Id}{I^{\cdot}}
\newcommand{\Pd}{P^{\cdot}}

\newcommand{\Ga}{\Gamma}
\newcommand{\gstar}{\Ga_*}
\newcommand{\rgstar}{{\bf R}\gstar}
\newcommand{\trgstar}{\tau{\bf R}\gstar}
\newcommand{\gaf}{\Gamma_*(\gF)}
\newcommand{\gags}{\Gamma_{G,*}}
\newcommand{\rgags}{{\bf R}\Gamma_{G,*}}
\newcommand{\trgags}{\tau{\bf R}\Gamma_{G,*}}

\newcommand{\delsv}{\Delta_{\sv}}
\newcommand{\delkg}{\Delta_{k[G]}}
\newcommand{\strdelsv}{\overline{\Delta}_\sv}

\newcommand{\rderhom}{R \mathcal{H}om}
\newcommand{\homkn}{\mathcal {H}om}
\newcommand{\supp}{\text{Supp}\,}

\newcommand{\sus}{\subseteq}
\newcommand{\larr}{\longrightarrow}
\newcommand{\pil}{\rightarrow}
\newcommand{\lpil}{\larr}
\newcommand{\prpil}{\dashrightarrow}
\newcommand{\inpil}{\hookrightarrow}
\newcommand{\pils}{\twoheadrightarrow}
\newcommand{\spil}{\pils}
\newcommand{\bihom}[2]{\overset{#1}{\underset{#2}{\rightleftarrows}}}

\newcommand{\mto}[1]{\stackrel{#1}\longrightarrow}

\newcommand{\iso}{\cong}
\newcommand{\te}{\otimes}

\newcommand{\equi}{\backsimeq}
\newcommand{\isopil}{\overset{\cong}{\lpil}}
\newcommand{\ispil}{\isopil}
\newcommand{\vvi}{\langle}
\newcommand{\hvi}{\rangle}

\newcommand{\pn}{{\bf P}^r}
\newcommand{\pfour}{{{\bf P}^4}}
\newcommand{\pfire}{{{\bf P}^4}}
\newcommand{\ptre}{{{\bf P}^3}}
\newcommand{\pto}{{{\bf P}^2}}
\newcommand{\pdim}{v}

\newcommand{\hele}{{\bf Z}}
\newcommand{\ove}[1]{\bar{e}_{#1}}
\newcommand{\cove}[1]{\hat{\bar{e}}_{#1}}

\newcommand{\injEn}{{1.5.2}}
\newcommand{\FGdef}{2.1.4 }

\newcommand{\koseksplTre}{3.1}
\newcommand{\adjTo}{2.1.6}
\newcommand{\qisFem}{5.1.2}
\newcommand{\VWdiagTre}{3.5.4}

\newcommand{\finkdekvSyv}{7.2.3}
\newcommand{\dersvSyv}{7.2.4}
\newcommand{\LRtriangFem}{5.2.4 }
\newcommand{\hovedekvFem}{5.2.5' }
\newcommand{\ekvSyv}{7.1.4 }

\newcommand{\komlinisoAtte}{8.1.4 }
\newcommand{\filtAtte}{8.2.2 }

\newcommand{\tenqisSeks}{6.3.1 }
\newcommand{\bimodSeks}{6.3.2 }

\begin{document}

\title[Coherent sheaves via Koszul duality]
{Describing coherent sheaves on projective spaces via 
        Koszul duality}

\author { Gunnar Fl{\o}ystad}
\address{ Matematisk Institutt\\
          Johs. Brunsgt. 12 \\
          5008 Bergen \\
          Norway}   
        
\email{ gunnar@mi.uib.no }

\maketitle

\section*{Introduction}
  
It is well known that there is a close connection between coherent
sheaves on a projective space $\pv$ 
where $W$ is a a vector space over a field $k$,
and finitely generated
graded modules over the symmetric algebra $S(W)$. The
Bernstein-Gel'fand-Gel'fand (BGG) correspondence \cite{BGG} from 1978
relates coherent sheaves on $\pv$ with graded modules over the 
exterior algebra $\lvd = \oplus \wedge^i(V)$
where $V$ is the dual vector space of $W$. This
correspondence
may be seen as a composition of the first connection and the 
correspondence between (complexes of) graded modules over $S(W)$ and
(complexes of) graded modules over $\lvd$ 
coming from the fact that $S(W)$ and $\lvd$ are {\it dual} Koszul
algebras. This latter correspondence, called {\it Koszul duality}, 
stems again from \cite{BGG}, and is treated subsequently in \cite{BGS} and
\cite{Fl}.

Of course the relationship between coherent sheaves on $\pv$ and
graded modules
over the symmetric algebra $\sv$ has been widely used. In this paper
we shall investigate in detail the BGG-correspondence between 
complexes of coherent sheaves and complexes of graded modules over the
exterior algebra. Our claim is that for algebraic purposes, complexes
of graded modules over the exterior algebra $\lvd$ may be a more
natural tool for
investigating complexes of coherent sheaves on $\pv$ than are complexes
of graded modules over the symmetric algebra $\sv$.
To mention some applications of this we give a strikingly simple algebraic
construction of the Horrocks-Mumford bundle on $\pfire$, we get a very
natural proof of the Castelnuovo-Mumford theorem \cite{Mu} on the 
regularity of coherent sheaves
and we also give a generalization of a theorem of Barth \cite{Ba} 
on stable rank
two sheaves on $\pto$, to a form 
which holds for {\it all} coherent sheaves on a projective space
$\pv$ (and from which Barth's theorem is an immediate corollary).

\medskip

The BGG-correspondence states more precisely that there is an 
equivalence of categories between the bounded derived category of 
coherent sheaves on $\pv$ and the {\it stable} module category of 
finitely generated graded left $\lvd$-modules
\begin{equation}
\dbco \equi \modstab{\lvd}. \label{inBGGur}
\end{equation}

We shall study this correspondence from a slightly different angle.
Let $\gF$ be a coherent sheaf on $\pv$. Then $\oplus_{n \in \hele}
H^0 \gF(n)$ is a graded $\sv$-module. Let 
\[ \lv = \oplus_{i=0}^{\dim_k V} \Hom(\wedge^i(V),k) \]
be the graded dual of $\lvd$. Then $\lv$ is a graded left
$\lvd$-module.
Via Koszul duality the graded $\sv$-module 
$\oplus_{n \in \hele} H^0 \gF(n)$ gives rise to a
complex of graded left $\lvd$-modules
\begin{equation}
   \cdots \lpil \lv(p) \te_k H^0 \gF(p) \mto{d^{p}}             
   \lv(p+1) \te_k H^0 \gF(p+1) \lpil \cdots.  \label{inH0F}
\end{equation} 
This complex is exact in the $p$'th component when $p > $ regularity 
of $\gF$. Now note that $\lv$ is 
isomorphic as a graded left $\lvd$-module to $\lvd$ shifted $\dim_k V$ degrees 
to the left. Thus $\lv$ is a free $\lvd$-module.
Let $a \geq $ regularity $\gF$.
 We may then construct a 
minimal free resolution $P^{\cdot}$ of $\ker d^a$ where each component 
$P^p$ is a direct sum of (shifts of) $\lv$'s. Splicing this together
with the complex (\ref{inH0F}) truncated in components of degrees
$\geq a$, gives us an {\it acyclic} complex which can be 
identified as a complex (let $\pdim = \dim \pv$)
\begin{equation}
\cdots \pil \oplus_{i=0}^{\pdim} \lv(p-i) \te_k H^i \gF(p-i)
\mto{d^p}  \oplus_{i = 0}^{\pdim} \lv(p+1-i) \te_k H^i \gF(p+1-i)
\pil \cdots . \label{inHiF}
\end{equation}
Now let $\kolvdcf$ be the homotopy category of acyclic complexes
$\Id$ whose components $I^r$ are of the form
$ \oplus_{i \in \hele} \lv(-i) \te_k V_i^r $ with $\sum_{i \in \hele} 
\dim_k V_i^r$ finite.
Letting $\coh$ be the category of coherent sheaves on $\pv$, then the above 
construction (\ref{inHiF}) gives us a functor
\[ \coh \mto{T_0} \kolvdcf \]
which extends to a functor
\begin{equation} \dbco \mto{T} \kolvdcf \label{inekv}
\end{equation}
which becomes an {\it equivalence} of categories.
We call $\ytrev(\gG)$ the {\it Tate resolution} of $\gG$.
 This equivalence is the same as the
BGG-correspondence since it is standard that there is an equivalence
of categories (the Tate correspondence)
\[ \kolvdcf \equi \modstab{\lvd}. \]
However the version (\ref{inekv}) has several advantages compared to 
the version (\ref{inBGGur}). For instance there is the explicit way that
$\ytrev(\gF)$, given by (\ref{inHiF}), 
is related to the cohomology of $\gF$. It is this that
enables the above mentioned generalization of Barth's theorem on 
stable rank two sheaves on $\pto$, to a form which holds for all coherent
sheaves (see Remark \ref{3barth}).

Another reason, which also serves to illustrate the naturality of 
the exterior complex in the study of coherent sheaves on $\pv$,
comes from looking at how a coherent sheaf $\gF$ may be represented 
by a complex
of free $\sv$-modules. Let us look at some of these ways.

The most familiar may be the minimal free
resolution of $\oplus_{n \in \hele} H^0 \gF(n)$. The minimal free
resolution has certain ``higher'' versions, Walter complexes
$F_r^\cdot$, one for each integer $r \geq 0$ such that
$r$ is less than a certain integer depending on (the local projective
dimension of) $\gF$. They
are characterized by the facts that i) length $F_r^\cdot \leq
\pdim - 1$, ii) $H^i (F_r^\cdot) = \oplus_{n \in \hele} H^i \gF(n)$ for
$i = 0, \ldots, r$, iii) $H^i (F_r^\cdot) = 0$ otherwise, and 
iv) the sheafification
of ${F_r^{\cdot}}$ is quasi-isomorphic to $\gF$. 

Another way of representing a coherent sheaf by complexes of free
$\sv$-modules
comes from the Beilinson correspondence \cite{Bei}
(originating in the same journal edition as the BGG-correspondence)
which gives an equivalence of categories
\begin{equation} \dbco \equi K[-r-\pdim, -r] \label{inbei}
\end{equation}
where $K[-r-\pdim,-r]$ is the homotopy category of 
bounded complexes of finite rank free $\sv$-modules of the form
$\oplus_{i = r}^{r+\pdim} \sv(-i) \te_k V_i$. 
If $\gG$ corresponds to $F^\cdot$ in $K[-r-\pdim,-r]$, then again the 
sheafification of ${F}^\cdot$ is quasi-isomorphic 
to $\gG$.
This correspondence and the
Beilinson spectral sequence derived in the original proof of 
(\ref{inbei}) 
has been used extensively for example 
in the construction of vector bundles on $\pv$, in the study
of their moduli \cite{OSS}, and in the study of surfaces in $\pfour$
\cite{DSE}.

Given an object $\gG$ in $\dbco$ 
we show that all complexes $F^\cdot$ of free $\sv$-modules such that
the sheafification of ${F^{\cdot}}$ is quasi-isomorphic
to $\gG$ and which are reasonably nice (this includes resolutions, Walter
complexes, Beilinson complexes and rigid complexes), may be 
obtained by simply truncating the Tate resolution $\ytrev(\gG)$ at a 
suitable place and then transforming this via Koszul duality to 
a complex of free $\sv$-modules. Thus although complexes
of free $\sv$-modules representing $\gG$ manifest themselves 
in quite different forms, the corresponding exterior complexes are
{\it basically the same}.

\medskip

We thus propose to study objects in $\dbco$ by studying the 
corresponding exterior complex in $\kolvdcf$. First note
that an object $\Id$ in $\kolvdcf$ is completely determined
by any of its differentials
\begin{equation}
I^a \mto{d_{\Id}^a} I^{a+1} \label{indIa}
\end{equation}
To see this note that $\lv$ is both a projective and injective
$\lvd$-module. Thus the truncation $\sigma^{\leq a-1} {\Id}$ is a 
projective resolution of $\ker d^a_{\Id}$ and the truncation
$\sigma^{\geq a+1} \Id$ is an injective resolution of 
$\im d_{\Id}^a$. Thus $\Id$ is uniquely determined, up to
homotopy, by (\ref{indIa}). 

Conversely, given {\it any} map
\begin{equation}
\oplus_{q \in \hele} \lv(-q) \te_k V_q \mto{d} \oplus_{q \in \hele}
\lv(-q) \te_k W_q \label{indiff}
\end{equation}
between finitely generated graded left $\lvd$-modules. By taking
a projective resolution of $\ker d$ and an injective resolution of 
$\coker d$, with components in the resolutions 
consisting of finite direct sums of (shifts of)
$\lv$, we get an object in $\kolvdcf$ and thus an object in
$\dbco$. This gives, at least in principle, a great amount of
freedom
in constructing objects in $\dbco$. 

\medskip

But how do we determine properties of a coherent sheaf $\gF$ on $\pv$ from
the corresponding Tate resolution $\Id = \ytrev(\gF)$ ? Let us look at some
properties which are often of interest to determine.

Firstly, from the form (\ref{inHiF}) of $\Id$ we see that the 
cohomology of 
$\gF$ may be determined from $\Id$. 

Secondly, the minimal free resolution of $\oplus_{n \in \hele} H^0 \gF(n)$ 
\begin{equation}
 \cdots \pil \oplus_{q \in \hele} \sv(-q) \te_k V_q^p \pil \cdots \pil
\oplus_{q \in \hele} \sv(-q) \te_k V_q^0 \pil
\oplus_{n \in \hele} H^0 \gF(n) \label{inminfri}
\end{equation}
is often of interest because the syzygies $V^p_q$ 
of order $p$ and weight $q$ are attached geometric significance \cite{Gr}. 
Via Koszul duality the exterior complex corresponding to the 
minimal free resolution (\ref{inminfri}) is the
complex
(\ref{inH0F}). The syzygies $V^p_q$ may be computed as the graded 
pieces of the cohomology of the complex (\ref{inH0F}) (after a
suitable re-indexing), see Subsection \ref{4koskohSek}. 
In fact the graded pieces of the cohomology of 
(\ref{inH0F}) are just the
Koszul cohomology groups of $\oplus_{n \in \hele} H^0\gF(n)$ as defined
in \cite{Gr}.


Thirdly, we show, Theorem \ref{5ker}, that
the Hilbert polynomial of $\gF$ may be computed from the
dimensions of the graded pieces of the kernel of $d^a_{\Id}$ 
(see (\ref{indIa})) for any $a$.

As a fourth topic we turn to the question of how local properties of
the coherent
sheaf $\gF$ can be determined from $\Id$. We show, Corollary \ref{6lokalP}, 
that the
rank of $\gF$ at any point $P$ in $\pv$ 
can be found by a quite local computation on $\Id$ 
involving only the terms (for arbitrary $a$)
\begin{equation}
I^{a-1} \mto{d^{a-1}} I^a \mto{d^a} I^{a+1}. \label{indIlok}
\end{equation}
We also show
how to compute the projective dimension of the
localization
$\gF_P$ at any point $P$ in $\pv$
from (\ref{indIlok}).

Conversely, if one starts with an exact sequence (\ref{indIlok}),
complete this to an object $\Id$ in $\kolvdcf$ and we give criteria,
Theorem \ref{6dtilF},
for when the corresponding object $\gG$ in
$\dbco$ is (quasi-isomorphic to) a coherent sheaf.

All this is particularly useful if one wants to determine if a 
complex $\Id$ gives rise to a vector bundle.

\medskip 

Another feature which is easily described via the Tate resolution is 
projection. Let $U \sus W$ be a linear subspace. Then there is a projection
$\pv \prpil \pw$ with $\pvw$ as the center of projection. If $\gF$ is a 
coherent sheaf on $\pv$, with Supp$\,\gF$ disjoint from $\pvw$, then 
$p_* \gF$ is a coherent sheaf on $\pw$. The Tate resolution, $\ytrev(p_* \gF)$,
is then given by $\Hom_{\lWd}(\lwd, \ytrev(\gF))$.

\medskip

A way of constructing a map (\ref{indiff}) would be to take three
representations $A,B$, and $W$ of a group $G$ together with a 
$G$-equivariant map $W \te_k A \pil B$. This gives rise to a 
$G$-equivariant map
\begin{equation}
\lv(-1) \te_k A \pil \lv \te_k B
\end{equation}
If $A,B$ and $W$ belong to a suitable abelian category 
of $G$-representations which is 
semi-simple, this can be
completed to a $G$-equivariant exterior complex $\Id$. The associated
object $\gG$ in $\dbco$ should then be a $G$-equivariant coherent sheaf.
We develop the theory for such conclusions and this gives us an
equivalence of categories between the $G$-equivariant bounded derived
category of coherent sheaves and the $G$-equivariant version of 
$\kolvdcf$
\begin{equation*}
D^b_G(\coh) \equi K^\circ_G(\lvdcf) 
\end{equation*}

\medskip We now proceed to give an overview of how the paper is
organized.

\medskip 

In Section 1 we recall the basic theory of Koszul duality which is 
relevant for this paper. This is mostly the theory developed
in \cite{Fl} rephrased in the case where the dual Koszul algebras are
$\lvd$ and $\sv$. In particular we define the adjoint Koszul functors
\[ \koml \bihom{\fl}{\gs} \koms \]
between the categories of complexes of graded left $\lvd$-modules and 
graded $\sv$-modules respectively.

Let 
\[ \gstar\, : \, \coh \pil \modv{\sv} \]
be the graded global section functor given by
\[ \gaf = \oplus_{n \in \hele} \Gamma(\pv,\gF(n)). \]
Then we get a derived functor
\[ \rgstar \, : \, \dbco \pil \koms. \]
We show in Section 2 that when $\gG$ is in 
$\dbco$, the composition $\gs \circ \rgstar (\gG)$
is an acyclic complex in $\koml$. The proof quickly reduces to 
show that $\gs \circ \rgstar (\gopv)$ is an acyclic complex in
$\koml$.

Section 3 establishes the equivalence of categories
\begin{equation}
\dbco \bihom{T}{Sh} \kolvdcf
\end{equation}
and gives very explicit descriptions of the functors $\ytrev$ and $Sh$.
We also establish the form of $\ytrev(\gF)$ for a coherent sheaf, given
in (\ref{inHiF}).

Section 4 considers bounded complexes of finite rank 
free $\sv$-modules $\Pd$ such
that the sheafification of ${\Pd}$ is quasi-isomorphic to a
given object $\gG$ in $\dbco$. We demonstrate that such $\Pd$'s
are obtained essentially by truncating the Tate resolution
$\ytrev(\gG)$ and then applying the Koszul functor
$\fl$. In particular for certain canonical truncations of 
$\ytrev(\gG)$ we get corresponding canonical complexes $\Pd$,
like minimal free resolutions, Walter complexes, Beilinson complexes
and linear complexes. We also discuss Koszul cohomology and
Castelnuovo-Mumford regularity of coherent sheaves.

Section 5 shows how to compute the Hilbert polynomial of $\gG$
(suitably defined for a complex of coherent sheaves)
from the Tate resolution $\ytrev(\gG)$ or rather from the kernel of 
the differential $d^a_{\ytrev(\gG)}$ for any $a$.

Section 6 is devoted to study how local properties of a coherent
sheaf $\gF$ can be determined from the Tate resolution
$\Id = \ytrev(\gF)$. We show, Corollary \ref{6lokalP}, 
how the rank of $\gF$ at any point $P$ in $\pv$
may be determined from the part
\begin{equation}
I^{a-1} \mto{d^{a-1}} I^a \mto{d^a} I^{a+1} \label{inIlok2}
\end{equation}
of $\Id$ for any $a$.
We also show how the projective dimension of the localization
$\gF_P$ at a point $P$ may be determined from (\ref{inIlok2}).

Conversely, given an exact sequence (\ref{inIlok2}) we give 
sufficient criteria, Theorem \ref{6dtilF},
for the corresponding object $\gG$ in 
$\dbco$ to be (quasi-isomorphic to) a coherent sheaf.

Section 7 studies the projections $\pv \prpil \pw$ for a linear 
subspace $U \sus W$. If $\gF$ is a coherent sheaf on $\pv$ with 
support disjoint form $\pvw$, we show that the Tate resolution
\[ \ytrev(p_*\gF) = \Hom_{\lWd}(\lwd, \ytrev(\gF)).\] 
Actually for the subcategory of $\dbco$ consisting of $\gG$ with 
the support of all $H^i(\gG)$ disjoint from $\pvw$ we show that such
complexes may be projected down to complexes in $\dbcow$ and we show
that the corresponding functor on Tate resolutions is 
$\Hom_{\lWd}(\lwd, -)$.

Section 8 develops some general theory about the correspondence
between $G$-equivariant coherent sheaves on $\pv$ and graded
modules over $\sv$ whose module structure is compatible
with the $G$-action. For the expert on 
representation theory this section is rather obvious and the 
proofs of some propositions may be seen as overdoing it, but
we include it because we don't know of a good reference covering the
cases we need to consider. This section is independent of the rest
of the paper.

Section 9 gives the $G$-equivariant versions of the most important
theorems in this article, Theorem \ref{3ekv} and Theorem \ref{3ekspl}.

Section 10 consists of some examples. We give a construction of 
the Horrocks-Mumford bundle by constructing a part of its exterior 
complex. This seems a very natural way of constructing it and
requires almost no cleverness at all.

We also consider $GL(W)$-equivariant vector bundles on $\pv$ and the 
corresponding $GL(W)$-equivariant exterior complexes.

\medskip \noindent {\it Acknowledgments.} Most of this paper and the 
preceding paper \cite{Fl} were written during my sabbatical at MIT the 
academic year 99/00 and I thank MIT for their hospitality. 

The origins of this paper stems from investigations of a conjecture
in \cite{Fl2} using Macaulay 2 and we would like to state our appreciation
of this program.

\medskip \noindent {\it Note.} 
Much of the material of this paper was independently developed by 
D. Eisenbud and F.-O. Schreyer in a preprint published at the same time
as this paper appeared as a preprint. It then seemed to us to be the most
beneficial for the mathematical community that we cooperate to write a joint
more extended version. Since the two papers were quite distinct in 
approach a direct merger of the two papers did not seem desirable. While
in the present  paper we extensively use the language of derived and 
triangulated
categories, in the preprint be Eisenbud and Schreyer they tried to avoid
this language. We then wrote a join paper \cite{EFS} based on 
the original preprint by Eisenbud and Schreyer. Therefore all the basic
ideas and results in Sections 3,4,5, and 10 should be considered joint
work with Eisenbud and Schreyer, although the specific form and proof
will usually be different from that of \cite{EFS}.

Also the paper \cite{EPY} contains results, Theorem 4.1, which are 
equivalent to the results in Subsection 6.2.

The notation in the present paper has been somewhat changed from the 
original notation, so that it is more aligned with the notation of
\cite{EFS}.

\section{Preliminaries.}

In this section we shall recall the facts from \cite{Fl} which we need. In 
that paper we studied quadratic dual Koszul algebras $A$ and $\ad$. Here we 
shall recall and state results from \cite{Fl} specialized to the case where 
$A$ and $\ad$ are the exterior algebra $\lvd$ and the symmetric 
algebra $\sv$, where $V$ is a finite-dimensional vector space over a field 
$k$ and $W = V^*$ is the dual vector space.

\subsection {Notation.} \label{1not}

As just said, let $k$ be a field and $V$ a finite-dimensional vector space 
over $k$ with $W = V^*$. Let $\pdim + 1 = \dim_k V$ so $\pdim$ is the 
dimension of the projective space $\pv = {\bf Proj}(S(W))$.
When tensoring over $k$ we shall normally drop $k$ as a subscript of
$\te$.

We can form the tensor algebra 
\[ T_k(V) = k \oplus V \oplus (V \te V) \oplus \cdots \oplus V^{\te n} \oplus 
\cdots \]
and also the tensor algebra $T_k(W)$.

In \cite{Fl} we studied quadratic dual algebras $A$ and $\ad$ which are 
quotients of $T_k(V)$ and $T_k (W)$ respectively. In the following we shall 
only be interested in the case where $A$ is the exterior algebra $\lvd 
= \oplus_{i=0}^{\pdim + 1} \wedge^i V$ and 
$\ad$ is the symmetric algebra $S(W) = \oplus_{i \geq 0} 
\text{Sym}^i(W)$. We consider $V$ and $W$ to have degree $1$ so $\lvd$ and
$S(W)$ are positively graded algebras. (This contrasts with the convention
in \cite{EFS}, where $V$ is considered to have degree $-1$.)
\medskip




\medskip

If $\Bc$ is an abelian category we let $\Kom(\Bc)$ be the category of 
{\it complexes of objects} in $\Bc$. We shall also let $K(\Bc)$ be the 
{\it homotopy 
category} of 
complexes of objects in $\Bc$, i.e. $K(\Bc)$ has the same objects as 
$\Kom(\Bc)$ but homotopic morphism are identified.
   
We let $D(\Bc)$ be the {\it derived category} associated to $\Bc$ and $D^b(\Bc) 
\sus D(\Bc)$ the full subcategory consisting of bounded complexes. If $\Ac 
\sus \Bc$ is a thick abelian subcategory, i.e. it is closed under extensions, 
then we let $D^b_{\Ac}(\Bc)$ be the full subcategory of $D^b(\Bc)$ consisting 
of complexes whose cohomology is in $\Ac$.
We also let $D_{b, \Ac}(\Bc)$ be the full subcategory of $D(\Bc)$ consisting 
of complexes $X$ such that $H^i(X)$ is in $\Ac$ for all $i$ and $H^i(X)$ is 
nonzero for only a finite number of $i$.

\medskip

If $M$ is a complex in $\Kom(\Bc)$ and $r$ an integer, 
we let $M[r]$ denote the complex shifted
$r$ places to the left, i.e. $M[r]^p = M^{r+p}$ and
$d_{M[r]}^p = (-1)^r d_M^{p+r}$.



\medskip

Let $B = \oplus_{i \geq 0} B_i$ be a positively graded associative algebra 
with $k = B_0$ a central field. 
We let $\vmod{B}$ be the category of graded left 
$B$-modules with homomorphisms of degree $0$. If $\Bc$ is the category $\vmod{B}$ 
we write $\Kom(B), K(B)$, and $D(B)$ instead of $\Kom(\vmod{B}), K(\vmod{B})$, 
and $D(\vmod{B})$.

\medskip

If $M$ is in $\vmod{B}$ or $\Kom(B)$ and $r$ an integer, we denote by $M(r)$
the module or complex of modules with a shift of $r$ in the $B$-module grading,
i.e. $M(r)^p_q = M^p_{r+q}$ (when $M$ is a complex).

\medskip

A graded left $B$-module of the form 
\begin{equation} \oplus_{q \in \hele} B(-q) \te V_{q} \label{1f} 
\end{equation}
where $V_{q}$ is a vector space over $k$ is called a {\it free} $B$-module. 
Its {\it rank}
is $\sum_{q \in \hele} \dim_k V_{q}$, and we note that (\ref{1f}) is a 
projective module.

A graded left $B$-module of the form
\begin{equation} \Pi_{q \in \hele} \Hom_k(B(q),V_q) \label{1cf}
\end{equation}
is a {\it cofree} $B$-module. Its {\it corank} is 
$\sum_{q \in \hele} \dim_k V_{q}$, and we note that (\ref{1cf}) is an injective
module (see 
\cite[Lemma \injEn]{Fl}).
Let
 $\bgd = \oplus_{i \leq 0}\Hom_k(B_{-i}, k)$. This is the {\it graded 
dual} of $B$ and it is a $B$-bimodule.
If each $B_{-i}$ is finite dimensional then
\begin{equation*} \Pi_{q \in \hele} \Hom_k(B(q),V_q) \iso 
\Pi_{q \in \hele} \bgd(-q) \te V_{q}. 
\end{equation*}

We let $\lvdcf$ be the category of cofree left $\lvd$-modules, and we let 
$\svf$ be the category of free $\sv$-modules.

\subsection{ The graded dual of $\lvd$.} \label{1lv}

 There is a perfect pairing 
\[ \wedge^p(V) \te \wedge^p(V^*) \lpil k \]
given by 
\[ u_1 \wedge \cdots \wedge u_p \te \alpha_1 \wedge \cdots \wedge \alpha_p \mapsto 
\sum_{\sigma} \text{sgn}(\sigma) \alpha_{\sigma(1)}(u_1) \cdots 
\alpha_{\sigma(p)}(u_p) = \text{det}(\alpha_j(u_i)). \]
where $\sigma$ runs over all permutations of $\{ 1, \ldots, p \}$. This 
pairing is denoted by $<\,  , \,  >$.

From this we get an isomorphism 
\[ \wedge^p(V^*) \mto{i} \wedge^p(V)^*. \]

The left $\lvd$-module structure on $\lvd^\gd$ gives maps
\[ \wedge^p(V) \te \wedge^{p+q}(V^*) \pil \wedge^q(V^*). \]
The map $u \te \alpha \mapsto u \alpha$ is determined by $<w, u\alpha> = 
<w \wedge u, \alpha> $ for $w \in \wedge^q(V)$. More explicitly it is given 
as follows. Let $u = u_1 \wedge \cdots \wedge u_p$ and $\alpha = \alpha_1 
\wedge \cdots \wedge \alpha_{p+q}$. Then 
\begin{equation} \label{1au}
u \alpha = \sum_{\sigma} \text{sgn}(\sigma) \alpha_{\sigma(q+1)}(u_{1}) 
\cdots \alpha_{\sigma(q+p)}(u_{p}) \alpha_{\sigma(1)} \wedge  \cdots 
\alpha_{\sigma(q)}, 
\end{equation}
the sum over all permutations $\sigma$ of $\{ 1, \ldots, p+q \}$ preserving 
the order of $\{ 1, \ldots, q\}$.

The graded module $\lvd^\gd$ is the canonical module 
$\lv$ of the finite dimensional $k$-algebra $\lvd$.

\begin{lemma} The natural map $\lvd \te \wedge^{\pdim + 1}(V^*) \pil
\lv$ is an isomorphism of left $\lvd$-modules.
\end{lemma}

\begin{proof} Let $\alpha$ be a generator of the one-dimensional space 
$\wedge^{\pdim + 1}(V^*)$. 
It will be sufficient to show that $u \alpha \neq 0$ for any nonzero $u$ 
in $\wedge^p(V)$. But since 
the pairing 
\[ \wedge^{\pdim + 1-p}(V) \te \wedge^p(V)
\pil \wedge^{\pdim + 1}(V) \]
is perfect we may find $w$ in $\wedge^{v+1-p}(V)$ such that
$w \wedge u \neq 0$. But then
\[ < w, u \alpha > = < w \wedge u, \alpha > \neq 0. \]
\end{proof}

Hence $\lvd$ is a Gorenstein ring. 
As a $k$-algebra it is also called a {\it Frobenius algebra} 
which for positively graded algebras $B$ are algebras with 
$B^{\gd} \iso B(r)$ for some integer $r$, as left $B$-modules.
We see that $\lvd$ is a Frobenius algebra.
A module over a Frobenius algebra is projective 
if and only if it is injective (see
\cite[Theorem 4.2.4]{We}). Also since $B_d$ must be nonzero for only a finite
number of degrees $d$, 
the concepts of free and cofree modules coincide
for a positively graded Frobenius algebra. 
It might therefore seem unnecessary to use the term cofree when speaking of
$\lvd$-modules.
However we shall continue to use
the concepts cofree and corank when speaking of $\lvd$-modules, when 
these are the natural concepts occurring in Koszul duality.

\subsection{ The Koszul functors. } \label{3sekKF}

By \cite{Fl} there are functors
\begin{eqnarray*}
\fl \,  : & \, \koml & \lpil \koms \\
\gs \,  : & \, \koms & \lpil \koml
\end{eqnarray*}
which we shall define. Let us first say that there are two ways of
defining them. One is conceptual and compact but not very explicit.
This is Definition \FGdef  in \cite{Fl}. 
We shall however give a 
more explicit definition \cite[Subsec. \koseksplTre]{Fl}, 
which is also the traditional way of defining these
functors.

Let $M$ be in $\koml$. The graded module in component $p$ is 
denoted $M^p$ and its graded piece of degree $q$ is denoted $M^p_q$. The 
complex $M$ comes with a differential $d_M$. We define $\fl(M)$ to be the 
{\it total direct sum complex}  of the double complex (unbounded in all 
directions, also down and to the left)
\[ \begin{CD} 
\sv \te M^2_0 @>>> \cdots   @.  \\
@AAA @AAA @. \\
 \sv \te  M^1_0 @>>> \sv(1) \te M^1_1 @>>> \cdots \\
 @AA{(d_v)_0^0}A @AAA @AAA \\
\sv \te M^0_0 @>{(d_h)_0^0}>>  \sv(1) \te M^0_1 @>>> \sv(2) \te M^0_2 
\end{CD} \]
where the vertical differential 
\[ (d_v)_q^p = \id_{\sv(q)} \te (d_M)_q^p \]
and the horizontal differential is given by 
\[  (d_h)_q^p(s \te m) = \sum_{\alpha \in A} sv_\alpha \te \check{v}_\alpha m 
\]
where $\{ v_\alpha \}_{\alpha \in A}$ is a basis for $W$ and $\{ 
\check{v}_\alpha \}_{\alpha \in A} $ is a dual basis for $V$.

For $N$ in $\koms$ we define $\gs(N)$ to be the {\it total complex} 
of the double complex (unbounded in all directions, also down and to
the left)
\[ \begin{CD}
\lv \te N^2_0 @>>> \cdots  @. \\
@AAA @AAA @. \\
 \lv \te  N^1_0 @>>> \lv(1) \te N^1_1 @>>>  \cdots \\
 @AA{(d_v)_0^0}A @AAA @AAA \\
\lv \te N^0_0 @>{(d_h)_0^0}>>  \lv(1) \te N^0_1 @>>> \lv(2) \te N^0_2. 
\end{CD} \]
Naturally one should here take the total {\it direct product} complex with 
products 
in the category of graded left $\lvd$-modules. However, since $\lv$ 
is non-zero only in a finite number of degrees, 
in the case above this is the same as the total {\it
direct sum} complex.

The horizontal differential is given by 
\[ (d_h)_q^p = \id_{\lv(q)} \te (d_N)_q^p \]
and the vertical differential is given by
\[ (d_v)_q^p (l \te n) = \sum_{\alpha \in A} l \check{v}_\alpha \te v_\alpha 
n. \]

\eks \llabel{1ekskoskom}
Let $M$ be a graded left $\lvd$- module. If we consider it as a complex 
with $M$ in the component of degree zero, then $\fl(M)$ is the complex
\[ \cdots \pil \sv(-1) \te M_{-1} \pil \sv \te M_0 
\pil \sv(1) \te M_1 \pil \cdots .\]
In particular $\fl(k) = \sv$ and $\fl(\lv)$ is just the Koszul complex and 
thus a free resolution of $k$.

If $N$ is a graded $\sv$-module, then considering it as a complex with $N$ 
in the component of degree zero, $\gs(N)$ is the complex
\[ \cdots \pil \lv(-1) \te N_{-1} \pil \lv \te N_0 
\pil \lv(1) \te N_1 \pil \cdots .\]
In particular $\gs(k) = \lv$ and $\gs(\sv)$ is just a rearrangement of the 
Koszul complex making it a cofree resolution of $k$.

Combining the above we see that $\fl \circ \gs (k) $ is quasi-isomorphic to 
$k$ and 
$\gs \circ \fl (k)$ is quasi-isomorphic to $k$.
\eksfin

The functors $\fl$ and $\gs$ are {\it exact}, i.e. they take short 
exact sequences of complexes to short exact sequences of complexes.
(It is {\it not} true however that $\fl$ takes acyclic complexes to 
acyclic complexes.)


\vskip 2mm
If $M$ is in $\Kom(\lvd)$ and $N$ is in $\Kom(\sv)$
we have 
the following identities 
\begin{eqnarray}  \fl(M(a)[b]) & = & \fl(M)(-a)[a+b] \\
\gs(N(a)[b]) & = & \gs(N)(-a)[a+b]. \label{1tvister} 
\end{eqnarray} 
\medskip

By \cite[Cor. \adjTo]{Fl} the functor $\fl$ is left adjoint to the functor 
$\gs$, i.e. 
for $M$ in $\koml$ and $N$ in $\koms$
there is a natural isomorphism
\begin{equation} \label{1FGadj}
\Hom_{\koms}(\fl(M), N)  \iso \Hom_{\koml}(M , \gs(N)). 
\end{equation}
This adjunction gives natural morphisms
\begin{gather} \label{1quasiso}
\fl \circ \gs (N) \lpil N  \\
M \lpil  \gs \circ \fl (M) \notag
\end{gather}
which are quasi-isomorphisms by \cite[Prop. \qisFem]{Fl}.

\subsection{ Subspaces of $W$.} \label{1subsp}

If $A \pil B$ is a homomorphism of positively graded algebras, we get 
functors between module categories. The functor
\[ \res^B_A \, : \, \, \vmod{B} \lpil \vmod{A} \]
is the restriction functor. It has a left adjoint functor
\[ B \te_A - \, : \, \, \vmod{A} \lpil \vmod{B} \]
and a right adjoint functor
\[ \Hom_A(B, -) \, : \,\, \vmod{A} \lpil \vmod{B}. \]
These functors extend to functors between the categories $\Kom(A)$ and
$\Kom(B)$.

Let $U \sus W$ be a vector subspace. Then we have morphisms of algebras
\begin{gather*}
\sv  \lpil  \svw \\
\lvwd  \lpil  \lWd 
\end{gather*}
There is then a diagram of functors
\begin{equation} \label{1resVW}
\begin{CD} 
\Kom(\lvd) @>{\fl}>>  \Kom(\svf) \\
@V{\res^{\lvd}_{\lvwd}}VV  @VV{\svw \te_{\sv}-}V \\
\Kom(\lvwd) @>{\flvw}>> \Kom(\svwf).
\end{CD} \end{equation}
which by \cite[Prop. \VWdiagTre]{Fl} gives a natural isomorphism of functors
\[ (\svw \te_{\sv} -) \circ \fl \cong \flvw \circ \res^{\lvd}_{\lvwd}. \]
There is also a diagram of functors
\begin{equation} \label{1resVWto}  \begin{CD} 
\Kom(\lvdcf) @<{\gs}<< \Kom(\sv) \\
@A{\Hom_{\lvwd}(\lvd, -)}AA   @AA{\res^{\svw}_{\sv}}A \\
\Kom(\lvwdcf) @<{\gsvw}<< \Kom(\svw) 
\end{CD} \end{equation}
which by \cite[Prop. \VWdiagTre]{Fl} gives a natural isomorphism of functors
\[ \Hom_{\lvwd}(\lvd, -) \circ \gsvw \cong \gs \circ \res^{\svw}_{\sv}. \]

\subsection{ Equivalence of categories. }
\label{1ekvSek}
Let $\vfmod{\sv}$ be the category of finitely 
generated $\sv$-modules
and let $\vfmod{\lvd}$ be the category of finitely generated
$\lvd$-modules.  The traditional equivalence of categories in 
Koszul duality, \cite[Thm. 2.12.6]{Bei} 
says that the functors $\fl$ and $\gs$ descend to give an 
(by abuse of notation we do not change the name of the functors)
{\it equivalence} of categories
\[ D^b(\vfmod{\lvd}) \bihom{\fl}{\tau \gs} 
D^b(\vfmod{\sv}) \]
where $\tau \gs$ is the functor $\gs$ followed 
by a suitable truncation of 
the complex. This result will however not be sufficient for our purposes,
since we both shall consider $\sv$-modules that are not finitely generated,
and consider complexes of $\lvd$-modules that are unbounded.
We therefore have to consider categories containing such objects.

The functors 
\[  \Kom(\lvd) \bihom{\fl}{\gs} \Kom(\sv) \]
takes homotopic morphisms to homotopic morphisms. Thus they "descend" to 
functors
\begin{eqnarray} K(\lvd) \bihom{\fl}{\sv} K(\sv) \label{1l-s} \end{eqnarray}
which are also adjoint. (By abuse of notation we don't change the name of the 
functors.)
Note that the functors $\fl$ and $\gs$ in (\ref{1l-s})
also restrict to give adjoint 
functors
\begin{eqnarray}
K(\lvdcf) & \bihom{\fl}{\gs} & K(\sv)  \label{1cf-s}\\
K(\lvd) & \bihom{\fl}{\gs} & K(\svf) \label{1l-f}\\
K(\lvdcf) & \bihom{\fl}{\gs} & K(\svf) \label{1cf-f}
\end{eqnarray}
It is a remarkable fact, as we shall state shortly, that the functors in 
(\ref{1cf-f}) give an equivalence of categories. Whether this is true for dual 
Koszul algebras in general we do not know but it does hold if one of them is 
finite dimensional \cite[Thm. \finkdekvSyv ']{Fl}.

Now let $N^R(\lvd)$ be the null system (see \cite[Def. 1.6.6]{KaSh} 
for more on this) of 
the triangulated category $K(\lvd)$ consisting of all objects $M$ in 
$K(\lvd)$ such that i. $M$ is acyclic and ii. $\fl(M)$ is acyclic. We get 
a triangulated category $D^R(\lvd) = K(\lvd)/N^R(\lvd)$ 
(\cite[Def. \LRtriangFem]{Fl}).

Similarly $N^L(\sv)$ is the null system of the triangulated category $K(\sv)$ 
consisting of all objects $N$ in $\Kom(\sv)$ such that i. $N$ is acyclic and 
ii. $\gs(N)$ is acyclic.
We get a triangulated category $D^L(\sv) = K(\sv)/N^L(\sv)$
 (see \cite[Def. \LRtriangFem]{Fl}).

By \cite[Remark \dersvSyv]{Fl}, $D^L(S(W))$ is isomorphic to the derived
category $D(S(W))$. It is not true however that $D^R(\lvd)$ is 
isomorphic to the derived category $D(\lvd)$.


The following is Theorem \hovedekvFem and Theorem \ekvSyv  
from \cite{Fl} in the case
where the dual Koszul algebras are $\lvd$ and $\sv$.

\begin{theorem} \label{1ekvkat}
The functors (\ref{1l-s})-(\ref{1cf-f}) all descend to give 
adjoint equivalences of categories (by abuse of notation we don't change the 
name of the functors)
\begin{eqnarray*} D^R(\lvd) & \bihom{\fl}{\gs} & D^L(\sv) \\
K(\lvdcf) & \bihom{\fl}{\gs}& D^L(\sv) \\
D^R(\lvd) & \bihom{\fl}{\gs} & K(\svf) \\
K(\lvdcf) & \bihom{\fl}{\gs} & K(\svf).
\end{eqnarray*} 
Furthermore if 
\begin{eqnarray*}
K(\lvdcf) & \mto{i_{\lvd}} & D^R(\lvd) \\
K(\svf) & \mto{i_{\sv}} & D^L(\sv) 
\end{eqnarray*}
are the inclusion functors, then there are natural isomorphisms of functors 
$i_{\lvd} \pil \gs \circ \fl$ and $\fl \circ \gs \pil i_{\sv}$ so $i_{\lvd}$ 
and $i_{\sv}$ both give equivalences of categories.
\end{theorem}

\subsection{ Filtrations.} \label{1filtsek}
If $M$ is in $\Kom(\lvd)$ (resp. $N$ is in $\Kom(\sv)$) then $\fl(M)$ (resp. 
$\gs(N)$) may be a rather "large" complex. We would like to find a "small" 
version of this complex. In Proposition \ref{1filt} below we give 
sufficient criteria 
for when to do this and also state what this "small" complex looks like.

Let $\KomLin(\lvd)$ be the full subcategory of $\Kom(\lvdcf)$ consisting of 
complexes of the form
\[  \cdots \lv(-1) \te L_{-1} \pil \lv \te L_0 \pil \lv(1) \te L_1 \pil \cdots 
. \]
There is a natural functor
\begin{equation} \label{1gskomlin}
 \vmod{\sv} \mto{\gs} \KomLin(\lvd). \end{equation}
This functor gives an isomorphism of categories 
\cite[Cor. \komlinisoAtte]{Fl}. 

Correspondingly we can define $\KomLin(\sv)$ and there is a natural functor
\[ \vmod{\lvd} \mto{\fl} \KomLin(\sv) \]
which gives an isomorphism of categories.

\medskip
A complex $P$ in $\Kom(\svf)$ is {\it minimal} if the differentials in $k 
\te_{\sv}P$ are zero. A complex $I$ in $\Kom(\lvdcf)$ is {\it minimal} if the 
differentials in $\Hom_{\lvd}(k, I)$ are zero. 

If $P$ in $\Kom(\svf)$ is a minimal complex with 
\[P^p = \oplus_{q \in \hele} \sv(-q) \te V^p_q \]
we may define a filtration
\[ \cdots \sus P\vvi{r-1}\hvi  \sus P\vvi{r}\hvi  \sus P\vvi{r+1}\hvi  \sus \cdots \]
where 
\[(P\vvi{r}\hvi )^p = \oplus_{p+q \leq r} \sv(-q) \te V^p_q. \]
We then get quotient complexes $Q\vvi{r}\hvi  = P\vvi{r}\hvi /P\vvi{r-1}\hvi $ which are complexes
\[ \cdots \sv(p-1-r) \te V^{p-1}_{r+1-p} \pil \sv(p-r) \te V^p_{r-p} \pil 
\sv(p+1-r) \te V^{p+1}_{r-1-p} \pil \cdots . \]
If $I$ in $\Kom(\lvdcf)$ is a minimal complex with 
\[ I^p = \oplus_{q \in \hele} \lv(-q) \te V^p_q \]
we may define a cofiltration
\[ \cdots \spil I\vvi{r-1}\hvi  \spil I\vvi{r}\hvi  \spil I\vvi{r+1}\hvi  \spil \cdots \]
where
\[ (I\vvi{r}\hvi )^p = \oplus_{p+q \geq r} \lv(-q) \te V^p_q. \]
We then get kernel complexes
$K\vvi{r}\hvi  = \ker (I\vvi{r}\hvi  \pil I\vvi{r+1}\hvi )$ which are complexes
\[ \cdots \pil \lv(p-1-r) \te V^{p-1}_{r+1-p} \pil \lv(p-r) \te V^p_{r-p} 
\pil \lv(p+1-r) \te V^{p+1}_{r-1-p} \pil \cdots . \]

The following is Theorem \filtAtte from \cite{Fl}
in the case where the dual Koszul algebras are $\lvd$ and $\sv$.

\begin{proposition} \label{1filt}
a. Let $N$ in $\Kom(\sv)$ be a bounded above complex. Then there is a homotopy
equivalence $I \pil \gs(N)$ where $I$ is a minimal complex, unique up to 
isomorphism in $\Kom(\lvdcf)$. The complex $I$ has bounded above cofiltration
and the kernels in the cofiltration of $I$ are given by 
$K\vvi{r}\hvi  = \gs(H^r(N)) [-r]$.

Conversely given a minimal complex $I$ with a bounded above cofiltration,
then $\fl (I)$ is a bounded above complex and the natural map 
$I \pil \gs \circ \fl (I)$ is a homotopy equivalence.

b. Let $M$ in $\Kom(\lvd)$ be a bounded below complex. Then there is a 
homotopy equivalence $\fl(M) \pil P$ where $P$ is a minimal complex, unique 
up to isomorphism in $\Kom(\svf)$. The complex $P$ has bounded below filtration
and the cokernels in the filtration of $P$ are given by
$ Q\vvi{r}\hvi  = \fl(H^r(M) )[-r]$.

Conversely given a minimal complex $P$ with a bounded below filtration, then
$\gs(P)$ is a bounded below complex and the natural map 
$\fl \circ \gs(P) \pil P$
is a homotopy equivalence.
\end{proposition}

Letting $\Kom^-(\sv)$ be the full subcategory of $\Kom(\sv)$ consisting of
bounded above complexes, we thus get a functor
\[ \gsmin \, : \, \Kom^-(\sv) \lpil \Kom(\lvdcf) \]
given by $N \mapsto I$.

Also, if we let $\Kom^+(\lvd)$ be the full subcategory of $\Kom(\lvd)$
consisting of bounded below complexes, we get a functor
\[  \flmin \, : \, \Kom^+(\lvd) \lpil \Kom(\svf) \]
given by $M \mapsto P$.

\subsection {Subspaces of $W$ and cohomology. } \label{1sekWFG}

Let $P$ and $I$ be complexes in $K(\svf)$ and $K(\lvdcf)$ respectively
such that $I \iso \gs(P)$ and thus $P \iso \fl(I)$. We would like to find 
out more about how these two complexes are related.

Let $U \sus W$  be a subspace. Then $\lwd$ is a $\lvd$-bimodule. 
If $N$ is in $\Kom(\lvd)$, we then
get the morphism complex (see Subsection \ref{1not}) $\Hom_{\lvd}(\lwd, N)$
which will be a complex of left $\lvd$-modules (and thus left $\lwd$-modules).
The following is Theorem \tenqisSeks and Corollary \bimodSeks in
\cite{Fl}.

\begin{proposition} \label{1PIht}
a. Let $P$ be in $\Kom(\svf)$.
There is a "twisted" quasi-isomorphism of complexes
\begin{equation} \label{1WFG}
 \Hom_{\lvd}(\lwd, \gs(P)) \mto{\alpha(P)} \svw \te_{\sv} P. 
\end{equation}
By twisted we mean that 
$ \Hom^p_{\lvd}(\lwd, \gs(P))_q$ maps to $(\svw \te_{\sv} P)^{p+q}_{-q}$.
From the first complex the cohomology comes equipped with a left 
$\lwd$-module structure. From the second complex the cohomology comes
with an $\svw$-module structure. These two actions of $\lwd$ and
$S(W/U)$ commute.

b. Let $I$ be in $\Kom(\lvdcf)$. There is a "twisted" quasi-isomorphism
of complexes
\[ \Hom_{\lvd}(\lwd, I) \mto{\beta(I)} \svw \te_{\sv} \fl(I). \]
The cohomology has a left $\lwd$-module structure and an $\svw$-module 
structure and these two actions commute.

c. If $I = \gs (P)$ then $\beta(I)$ composed with the canonical
map (see (\ref{1quasiso})) 
\[ \svw \te_{\sv} \fl \circ \gs(P) \pil \svw \te_{\sv} P \]
gives the map $\alpha(P)$.

d. If $P = \fl (I)$ then $\alpha(P)$ composed with the canonical map
(see (\ref{1quasiso})) 
\[ \Hom_{\lvd}(\lwd, I) \pil 
\Hom_{\lvd}(\lwd, \gs \circ \fl (I))\]
gives the map $\beta(I)$.
\end{proposition}

In {\it a.} denote the cohomology of the first complex as
\[ {}^I H = \oplus_{p \in \hele} H^p\Hom_{\lvd}(\lwd, \gs(P)) \]
and the cohomology of the second complex as 
\[ {}^{II} H = \oplus_{p \in \hele} H^p(\svw \te_{\sv} P).  \]
Then these modules are related by
\[ {}^I H^p_q = {}^{II} H^{p+q}_{-q} \text{ and } {}^I H^{p+q}_{-q} = 
{}^{II} H^p_q. \]
The module ${}^I H$ is a left $\lwd$-module with $l$ in $\La^d(W^*)$ acting
with bidegree $\binom{0}{d}$ while ${}^I H$ is an $\svw$-module 
with $s$ in $S^d(W/U)$ acting with bidegree $\binom{d}{-d}$.

Similarly ${}^{II} H$ is a left $\lwd$ module with $l$ in $\La^d(W^*)$ acting
with bidegree $\binom{d}{-d}$ while ${}^{II} H$ is an $\svw$ module 
with $s$ in $S^d(W/U)$ acting with bidegree $\binom{0}{d}$.

\eks \llabel{1kFG}
Let $U = W$.  Then from {\it a.} we get that 
\[ H^p(\gs(P))_q = H^{p+q}(k \te_{\sv} P)_{-q}. \]
If $M$ is a complex in $\Kom(\lvd)$ then if we let $P = \fl (M)$
and compose with the quasi-isomorphism $M \pil \gs \circ \fl (M)$ we
get that 
\[ H^p(M)_q = H^{p+q}(k \te_{\sv} \fl(M))_{-q}. \]

Let $U = 0$. Then from {\it b.} we get that 
\[ H^{p+q}\Hom_\lvd(k, I)_{-q} = H^p(\fl(I))_{q}. \]
If $N$ is a complex in $\Kom(\sv)$ then if we let $I = \gs (N)$
and compose with the quasi-isomorphism $ \fl \circ \gs (N) \pil N$ we
get that 
\[ H^{p+q}\Hom_{\sv}(k, \gs(N))_{-q} = H^p(N)_q. \]
\eksfin

\subsection {Group actions} \label{1sekGr}

Let $G$ be a linear algebraic group over the field $k$. We call 
a (possibly infinite dimensional) rational representation of $G$ with 
left $G$ action a {\it $G$-module}. The coordinate ring $k[G]$, is
a Hopf algebra and $W$ is a $G$-module if and only if $W$ is a left 
$k[G]$-comodule. The symmetric algebra $\sv$ becomes a $G$-module
and the algebra map $\sv \te \sv \pil \sv$ is a morphism of $G$-modules.
A graded module $M$ over $\sv$ is an {\it $\sv,G$-module} if $M$ is 
a $G$-module and the module map $\sv \te M \pil M$ is a $G$-module
map. Write $\vmod{\sv,G}$ for the category of $\sv,G$-modules. 
Similarly $\lvd$ is a $G$-module with the algebra map
$\lvd \te \lvd \pil \lvd$ a $G$-module map, and we get a category
$\vmod{\lvd, G}$.

\vskip 2mm

We let $\svgf$ be the full subcategory of $\vmod{\sv,G}$ whose objects are
\[ \oplus_{q \in \hele} \sv(-q) \te W_q \]
where the $W_q$ are $G$-modules. These are the {\it free} 
$\sv,G$-modules. Similarly $\lvdgcf$ is the full subcategory of
$\vmod{\lvd,G}$ whose objects are
\[ \oplus_{q \in \hele} \lv(-q) \te W_q \]
where the $W_q$ are $G$-modules. These are the {\it cofree} $\lvd,G$-modules.

In order for free $\sv,G$-modules to be projective in $\vmod{\sv,G}$
and cofree $\lvd,G$-modules to be injective in $\vmod{\lvd,G}$ we shall
henceforth assume that {\it the category of $G$-modules is semi-simple.}
I.e. short exact sequences of $G$-modules are split.  This holds
for instance if char $k = 0$ and $G$ is a finite or semi-simple group.

\vskip 2mm

For compact notation analogous to conventions in earlier paragraphs,
we denote the categories $K(\vmod{\sv,G})$ and $K(\svgf)$ as
$K_G(\sv)$ and $K_G(\svf)$. There is a forgetful functor
$K_G(\sv) \pil K(\sv)$, but is it faithful ?
Similarly we have categories $K_G(\lvd)$ and $K_G(\lvd \cfr)$. 

We get adjoint functors
\[ K_G(\lvd) \bihom{\fl}{\gs} K_G(\sv). \]
Let $N_G^L(\sv)$ be the null system in $K_G(\sv)$ whose objects are the
$N$ such that {\it i.} $N$ is acyclic and {\it ii.} $G_{\sv}(N) $ is acyclic. 

Similarly we have a null system $N_G^R(\lvd)$ in $K_G(\lvd)$ and we get
quotient triangulated categories
\[ D^R_G(\lvd) = K_G(\lvd)/N_G^R(\lvd), \quad 
D^L_G(\sv) = K_G(\sv)/N_G^L(\sv). \]
The category $D^L_G(\sv)$ is equal to the derived category of 
$\sv,G$-modules.
By \cite[Sec.10]{Fl} we have the analog of Theorem \ref{1ekvkat}.

\begin{theorem}
Assume that the category of $G$-modules is semi-simple. Then there are
adjoint equivalences of categories.
\begin{eqnarray*} D^R_G(\lvd) & \rightleftarrows & D^L_G(\sv) \\
K_G(\lvdcf) & \rightleftarrows & D^L_G(\sv) \\
D^R_G(\lvd) & \rightleftarrows & K_G(\svf) \\
K_G(\lvdcf) & \rightleftarrows & K_G(\svf).
\end{eqnarray*} 
\end{theorem}

Also we have the analog of Proposition \ref{1filt}
in the $G$-equivariant setting
provided the category of $G$-modules is semi-simple. We thus get a functor
\[ G_{\sv,min} : K^{-}_G(\sv) \lpil K_G(\lvd \cfr) \]
where $K^-_G(\sv)$ is the full subcategory of $K_G(\sv)$ consisting of 
bounded above complexes.

\section { Derived categories of sheaves on a projective space. } 
\label{2dersek}

This section is mostly to do preliminary work for the theory 
we develop in the next sections.  Let $\qc$ be the category of {\it 
quasi-coherent
sheaves on $\pv$}. It has full subcategories
\[ \vb \sus \coh \sus \qc \]
consisting of {\it locally free sheaves} (algebraic vector bundles) of
finte rank and 
{\it coherent sheaves} respectively. We then get the {\it derived categories} 
(see 
Subsection \ref{1not}) 
\[ \dbvb, \,\dbco, \, D^b_{\coh}(\qc), \,
D_{b,\coh}(\qc). \]
We shall usually for short write ``coh'' instead of 
``$\coh$'' in the subindex.
The first thing we show is that
these categories are all equivalent. 

There is also a well-known adjunction of functors
\[ \vmod{\sv} \bihom{\til}{\gstar(\pv,-)} \qc \]
where $\til$ is sheafification and 
\[ \gstar(\pv, \gF) = \oplus_{n \in \hele}
\Ga(\pv, \gF(n))\] 
is the graded global section functor. We then
get a derived functor 
\[ \rgstar(\pv,-) \, : \, \dbco \pil D(\sv) \]
which may be composed with the Koszul functor 
\[ D(\sv) \mto {\gs} 
\klvdcf . \]
We demonstrate the basic fact that
$\gs \circ \rgstar(\pv,\gG)$ is {\it acyclic} for all $\gG$ in $\dbco$.

\subsection {Equivalences of derived categories. }
\label{2derekv}
\begin{proposition} \label{2Pderekv} The natural maps
\[ \dbvb \mto{i_1} \dbco \mto{ i_2} \dbcq \mto {i_3} \dcq \]
all induce equivalences of categories.
\end{proposition}

\begin{proof} That $i_3$ induces an equivalence of categories is clear.

We now show first that $i_1$ and $i_2$ are fully faithful. We use criterion 
1.6.10 of  \cite{KaSh}. 
For the first inclusion we then need to show that if 
$\gG$ is a bounded complex of coherent sheaves, then there is a bounded
complex of vector bundles $\gE$ and a quasi-isomorphism $\gE \pil \gG$.
But $\gstar(\gG)$ is a bounded complex of finitely generated $\sv$-modules,
and this category has enough projectives. Thus we may find a quasi-isomorphism
$P \pil \gstar(\gG)$ where $P$ is a complex of free finitely generated
$\sv$-modules. Now sheafifying we get a quasi-isomorphism 
$ \tilov{P} \pil \gG$. Then let $\gE$ be the shafification $\tilov{P}$. 

The functor $j_1 \, : \, \dbco \pil \dbvb$ given by $j_1(\gG) = \gE$ gives
a quasi-inverse to the functor $i_1$, so the categories $\dbco$ and 
$\dbvb$ are equivalent.

\medskip

For the second map $i_2$ we shall show that given a bounded complex
$\gQ$ of quasi-coherent sheaves with coherent cohomology, there is a 
subcomplex $\gG \sus \gQ$ of coherent sheaves such that this inclusion
is a quasi-isomorphism. 
Suppose $\gQ$ is a complex $\cdots \pil \gQ_1 \pil \gQ_0 \pil 0$.
First find a coherent subsheaf $\gG_0 \sus \gQ_0$ such that $\gG_0 \pil H_0 
(\gQ)$ is 
surjective. Suppose by induction we have constructed complexes
\[ \begin{CD}
 @. \gG_k @>{d_k^{\gG}}>> \gG_{k-1} @>>> \cdots @>>> \gG_0 \\
@. @VVV @VVV @. @VVV \\
\gQ_{k+1} @>{d_{k+1}^{\gQ}}>> \gQ_k @>{d_k^{\gQ}}>> \gQ_{k-1}
@>>> \cdots @>>> \gQ_0 
\end{CD} \]
such that i. $\gG_i$ is a coherent subsheaf of $\gQ_i$ for each $i$,
ii. $H_i(\gG) \pil H_i(\gQ)$ is an isomorphism for $i < k$ and
iii. $\ker d_k^{\gG} \pil H_k(\gQ)$ is surjective.
Let $\gG_{k+1}^{\prime} = (d_{k+1}^{\gQ})^{-1}(\ker d_k^{\gG}) \sus
\gQ_{k+1}$ so that we get a pull-back diagram
\begin{equation} \label{2pullback}
\begin{CD}
\gG_{k+1}^{\prime} @>{d_{k+1}^{\gG^{\prime}}}>> \ker d_k^{\gG} \\
@VVV @VVV  \\
\gQ_{k+1} @>{d_{k+1}^{\gQ}}>> \ker d_k^{\gQ}. 
\end{CD} 
 \end{equation}
Note that $\ker d_{k+1}^{\gG^{\prime}} = \ker d_{k+1}^{\gQ}$. 
Consider the map 
\[ \alpha_k \, : \, \coker d_{k+1}^{\gG^{\prime}} \lpil \coker 
d_{k+1}^{\gQ}. \]
Since (\ref{2pullback}) is a pull-back diagram, it is clear that $\alpha_k$
is injective. By assumption {\it iii.} above it is also surjective and hence
an isomorphism. Now choose a coherent subsheaf 
$\gG_{k+1}^{\prime\prime} \sus \gG_{k+1}^{\prime}$ such that the composition
$\gG_{k+1}^{\prime\prime} \pil \gG_{k+1}^{\prime} \pil \im 
d_{k+1}^{\gG^{\prime}}$
is surjective, and choose a coherent subsheaf $\gG_{k+1}^{(3)}
\sus \ker d_{k+1}^{\gG^{\prime}}$ such that the composition
\[ \gG_{k+1}^{(3)} \pil  \ker d_{k+1}^{\gG^{\prime}} \ispil 
\ker d_{k+1}^{\gQ} \pil H^{k+1}(\gQ) \] 
is surjective.
Now let 
\[ \gG_{k+1} = \gG_{k+1}^{(3)} + \gG_{k+1}^{\prime\prime} \sus 
\gG_{k+1}^{\prime}. \]
Thus we may proceed inductively and construct a quasi-isomorphism 
$\gG \sus \gQ$ where $\gG$ is a complex of coherent sheaves.
The functor $j_2 \, : \, \dbcq \pil \dbco$ given by 
$j_2(\gQ) = \gG$ gives a    quasi-inverse to $i_2$.
\end{proof}

\subsection{ The graded global section functor. }

There are functors
\[ \vmod{\sv} \bihom{\til} {\gstar(\pv,-)} \qc \]
where $\til$ is the sheafification and 
\[ \gstar (\pv,\gF)  = \oplus_{n \in \hele} \Ga (\pv,\gF(n)) \]
where $\Ga(\pv,\gF(n))$ are the global sections of $\gF(n)$. 
If it is clear that we are considering coherent sheaves on
$\pv$ we write just $\Ga_*$ for $\Ga_*(\pv,-)$.
The functor
$\gstar$ is right adjoint to $\til$ and so we have an isomorphism
\begin{equation} \Hom_{\qc} (\tilov{M}, \gQ) \cong 
\Hom_{\vmod{\sv}}(M, \gstar(\gQ)). \label{2adj}
\end{equation}
Furthermore the natural map coming from the adjunction 
\[ \til \circ \gstar(\gQ) \pil \gQ \]
is an isomorphism.

We record the following for later use.

\begin{lemma} \label{2ggsadj}
There are adjunctions with $\til$ left adjoint

a. $\koms \bihom{\til}{\gstar} \Kom(\qc) $.

b. $\ksv \bihom{\til}{\gstar} K(\qc) $.
\end{lemma}

\begin{proof} Let $M$ be in $\koms$ and $\gQ$ be in $\Kom(\qc)$. Then
we clearly get an isomorphism of morphism complexes
\[ \Hom_{\qc}(\tilov{M}, \gQ) \cong \Hom_{\sv}(M, \gstar(\gQ)). \]
Taking cycles in degree zero of these complexes we get {\it a.} Taking 
homology
in degree zero of these complexes we get {\it b.}
\end{proof}

Now the category $\qc$ has enough injectives. We may therefore define the 
right derived functor (see \cite[10.5]{We}) 
\[ \rgstar \, : \, \dcq \lpil D(\sv) \]
which is a functor of triangulated categories.
If $\gQ$ is in $\dcq$ and $\gQ \pil \gI$ is a bounded below injective 
resolution, then by definition $\rgstar(\gQ) = \gstar(\gI)$. We denote
the cohomology group $H^i( \rgstar(\gQ))$ as $H^i_* \gQ$.

\begin{lemma} \label{2begfin}
Let $\gQ$ be in $\dcq$. Then $\rgstar(\gQ)$ has bounded 
cohomology
and for every integer $p$, each graded piece $(H^p_* \gQ)_q$ is finite 
dimensional.
\end{lemma}

\begin{proof} This is clearly true if $\gQ$ is a coherent sheaf (viewed
as a complex with this sheaf in the component of degree zero).
Since the coherent sheaves generate the triangulated category
$\dcq$ and $\rgstar$ is a functor of triangulated categories, we
get the lemma.
\end{proof}

Since $\rgstar(\gQ)$ has bounded cohomology we may define a functor
\begin{equation} \label{1rgstar}
\trgstar : \dcq \lpil D^{-}(\sv) 
\end{equation}
where $\trgstar(\gQ)$ is the complex $\rgstar(\gQ)$ suitably truncated
above such that $\trgstar(\gQ)$ is quasi-isomorphic to $\rgstar(\gQ)$.
This will be used in the statement of Theorem \ref{3ekv}.

\subsection{ Acyclicity of the corresponding exterior complex. }

We may compose the functor
\[ \rgstar(\pv,-) \, : \, \dcq \lpil D(\sv) \]
with the Koszul functor 
\[ \gs \, : \, D(\sv) \lpil \klvdcf. \]
We now show the following fundamental observation.

\begin{proposition} \label{2acext}
Let $\gQ$ be in $\dcq$. Then the complex 
$\gs \circ \rgstar(\pv,\gQ)$ is acyclic.
\end{proposition}

\begin{proof}

Note that $\dbco$ is generated, as a triangulated category, by
finite direct sums of $\gopv(-i)$ for integers $i$. Hence $\dcq$ is also 
generated by these objects. It is therefore enough to show
that $\gs \circ \rgstar(\pv,\gopv)$ is acyclic.

First assume that $W = (w)$ is one-dimensional. Then $\pv$ is a point
and $H^0_* \gopv \pil \rgstar(\pv,\gopv)$ is a quasi-isomorphism. Furthermore
$H^0_* \gopv = \oplus_{n \in \hele} k w^n$ as a module over $\sv$.  But then
$\gs \circ H^0_* (\gopv)$ is the complex (where $\lv = kw \oplus k$)
\[ \cdots \pil \lv(p-1) \te w^{p-1} \pil  \lv(p) \te w^p \pil 
\lv(p+1) \te w^{p+1} \pil \cdots \]
which is acyclic.

\medskip

Now let $W$ be arbitrary (finite dimensional) and let $w$ be an element of 
$W$.
Thus there is a short exact sequence
\[ \gopv(-1) \mto{\cdot w} \gopv \pil \gopvv \]
giving a triangle
\begin{eqnarray} \label{2hypertre}
  & &  \rgstar(\pv,\gopv(-1))  \\
& \pil &   \rgstar(\pv,\gopv) 
\pil \rgstar(\pv,\gopvv)  \notag \\
  & \pil & \rgstar(\pv,\gopv(-1))[1]  \notag \end{eqnarray} 
in $D(\sv)$. 
Now we also have a right derived graded global section functor
\[ \rgstar(\pvv,-) \, : \, \dcqvv \lpil D(\svv). \]
By induction we may assume that 
\[ \gsvv \circ \rgstar (\pvv,\gopvv) \]
is acyclic. 

\vskip 2mm
\noindent {\bf Claim.}  $\gs \circ \rgstar(\pv, \gopvv)$ is acyclic.
\vskip 2mm

\begin{proof}[Proof of claim]
By definition $\rgstar(\pv,\gopv)$ is $\gstar(\pv,\gI)$ where $\gI$ is an 
injective resolution
of $\gopv$ in $\qc$ and $\rgstar (\pvv,\gopvv)$ is  $\gstar(\pvv,\gJ)$ where
$\gJ$ is an injective resolution of $\gopvv$ in 
$\qcvv$. But each $\gJ^i$ is then flasque, \cite[III.2]{Ha}. 
Thus $\gJ$ is a flasque resolution of $\gopvv$ in $\qc$.
There will then be a
quasi-isomorphism of $\sv$-modules
\[ \res^{\svv}_{\sv} \circ \rgstar(\pvv,\gopvv) \lpil \rgstar(\pv,\gopvv) \]
(see \cite[III.2.5]{Ha} and \cite[III.1.2]{Ha}). 
By Theorem \ref{1ekvkat} there is then a quasi-isomorphism 
\begin{eqnarray*} &&\gs \circ \res^{\svv}_{\sv} \circ 
\rgstar(\pvv,\gopvv) \\
& \lpil & \gs \circ \rgstar(\pv,\gopvv) . \label{2grres} \end{eqnarray*}
By (\ref{1resVWto}) in Subsection \ref{1subsp} the former is isomorphic to
\begin{equation} \Hom_{\lvvd}(\lvd, \gsvv \circ \rgstar(\pvv,\gopvv)). 
\label{2homgr} 
\end{equation}
Since $\gsvv \circ \rgstar (\pvv,\gopvv)$ is acyclic by induction and
$\lvd$ is projective as a $\lvvd$-module,  (\ref{2homgr})
will be acyclic. Thus the last part of (\ref{2grres}) is also acyclic,
proving the claim.
\end{proof}

Now from the triangle (\ref{2hypertre}) we get a triangle
(recall (\ref{1tvister}))
\begin{eqnarray}
&&\gs \circ \rgstar(\pv, \gopvv)[-1]  \label{2Ghyptriang} \\
\pil && \gs \circ \rgstar(\pv, \gopv)(1)[-1] \pil 
\gs \circ \rgstar(\pv, \gopv) \notag \\
\pil && \gs \circ \rgstar(\pv, \gopvv) \notag
\end{eqnarray}

Let $I$ be $\gs \circ \rgstar(\pv, \gopv)$. From the triangle
above (\ref{2Ghyptriang}) we get
$H^{p-1}(I)_{q+1} = H^p(I)_q$. The following claim will then 
prove the proposition.


\vskip 2mm
\noindent {\bf Claim.} $H^p(I) = 0$ for $p \gg 0$.
\vskip 2mm

\begin{proof}[Proof of claim] 
 In the following let $\rgstar(-) = \rgstar(\pv,-)$.
Let $\tau^{\leq i}$ and $\tau^{>i}$ be the truncation functors
on $D(\sv)$ (see \cite[1.2.7]{We}). 
For a complex $N$ in $\koms$, let $w^{\leq i} N$
be the complex with $(w^{\leq i} N)^p = \oplus_{q \leq i} N^p_q$
and $w^{>i} N$ be the kernel of $N \pil w^{\leq i} N$. 
Then $w^{\leq i}$ and $w^{>i}$ are also functors on $D(\sv)$.

For a complex $N$ in $D(S)$ these truncation functors give triangles
\begin{eqnarray*}
\tau^{\leq i} M \pil &M \pil &\tau^{>i} M \pil \tau^{\leq i} M[1] \\ 
w^{>i} M \pil & M \pil & w^{\leq i} M \pil w^{>i} M [1].
\end{eqnarray*}
We thus get a triangle
\begin{equation} H^0_* \gopv \pil \rgstar(\gopv) \pil \tau^{>0} 
\rgstar(\gopv)
\pil H^0_* \gopv[1]. \label{2st0} \end{equation}
Since $H^p_* \gopv = 0$ for $p > n$, we have that 
$\tau^{>0} \rgstar(\gopv)$ is quasi-isomorphic to 
$ \tau^{\leq n} \tau^{>0} \rgstar(\gopv)$.
There is also a triangle
\begin{eqnarray*} &&w^{>q} \tau^{\leq n} \tau^{>0} \rgstar(\gopv) \\
&\pil & 
\tau^{\leq n} \tau^{>0}
\rgstar(\gopv) \pil w^{\leq q} \tau^{\leq n} \tau^{>0} \rgstar(\gopv) \\
& \pil & w^{> q} \tau^{\leq n} \tau^{>0} \rgstar(\gopv) [1] . 
\end{eqnarray*}
 Since $H^p \gopv(q)$ vanishes for $p > 0$ and  
$q \gg 0$ we get
$w^{> q} \tau^{\leq n} \tau^{>0} \rgstar(\gopv) $ isomorphic to $0$ 
in $D(\sv)$, and
the middle terms above are then isomorphic in $D(\sv)$ for $q \gg 0$.
 
 Now $\gs $ is a functor of triangulated categories. Applying $\gs$ to 
 (\ref{2st0}) we get a triangle
 \begin{eqnarray} 
&& \gs (\sv) \pil \gs \circ \rgstar(\gopv)  \pil  \gs \circ 
\tau^{>0} \rgstar(\gopv)  \label{2gst0} \\
& \pil &  \gs(\sv) [1]. \notag  \end{eqnarray}
 But now
 \[ \gs \circ \tau^{>0} \rgstar(\gopv) \iso 
 \gs \circ w^{\leq q} \tau^{\leq n} \tau^{>0} \rgstar(\gopv). \]
 By the definition of $\gs$ in Subsection \ref{3sekKF}
the latter complex is zero in large component 
degrees.
 Also, by Example \ref{1kFG}, $\gs(\sv)$ is exact in component degrees greater
 than zero.
 
Since we get from the triangle (\ref{2gst0}) a long exact sequence on 
 cohomology
 \begin{eqnarray*} \cdots & \pil H^p(\gs (\sv)) 
& \pil  H^p(\gs \circ \rgstar(\gopv))  \\
& \pil  H^p(\gs \circ \tau^{>0} \rgstar(\gopv))   
& \pil   H^{p+1} (\gs (\sv)) \pil \cdots 
 \end{eqnarray*}
 we see that $\gs \circ \rgstar(\gopv)$ is acyclic in large component
 degrees which proves the claim and also the proposition.
\end{proof} 
\end{proof}
 
 \rem In \cite[Section 9]{Fl} we showed that the category $D^L(\sv)$ 
which is the derived category  $D(\sv)$ 
 has two $t$-structures, which we called the {\it inner} and {\it outer} 
$t$-structures. 
 We get in this way to cohomological functors
 \begin{gather*}
 H^0_{in} \,  : \, D^L(\sv)  \lpil \vmod{\sv} \\
 H^0_{ou} \,  : \, D^L(\sv)  \lpil \vmod{\lvd}. 
 \end{gather*}
 The inner $t$-structure is just the standard $t$-structure and
 the functor $H^0_{in}$ is just the standard cohomological functor
 $H^0_{in} (N) = H^0(N)$. 

Also the category $D^R(\lvd)$ has two $t$-structures,
 the {\it inner} and {\it outer} $t$-structure. Via the equivalence of 
triangulated
categories
\[ D^L(\sv) \bihom{\gs}{\fl} D^R(\lvd) \]
these two $t$-structures are interchanged. The outer (non-standard) 
$t$-structure on
$D^L(\sv)$ corresponds to the inner  (standard) $t$-structure on $D^R(\lvd)$, 
and the inner (standard) $t$-structure on  
$D^L(\sv)$ corresponds to the outer (non-standard) $t$-structure on 
$D^R(\lvd)$.
Thus if $N$ is $\gs (M)$ then $H^i_{in}(N) \iso H^i_{ou}(M)$.
Proposition \ref{2acext} above says that if $\gQ$ is in $\dcq$ then 
the outer cohomology groups $H^p_{ou} (\rgstar(\gQ))$ vanish for all
integers $p$.

\section{ Complexes of coherent sheaves described as acyclic complexes
over the exterior algebra. }

  This section contains the main results of this article, Theorem \ref{3ekv} 
  and Theorem \ref{3ekspl}.
  
  \defi The category $\komolvdcf$ is the full subcategory of $\Kom(\lvdcf)$
  consisting of {\it acyclic} complexes whose components have 
  {\it finite corank}.
  The category $\kolvdcf$ is the full subcategory of 
\linebreak $K(\lvdcf)$ consisting
  of the same objects as $\komolvdcf$. 
 Note that it is a triangulated category. The objects of these categories
are called {\it Tate resolutions}.
 \defifin
 
 The result we shall prove, Theorem \ref{3ekv}, is that there is an equivalence
 of categories
 \[ \dbco \iso \kolvdcf. \]
 
Moreover we shall give a very explicit description of this correspondence,
Theorem \ref{3ekspl}. This result is closely related to the result
of Bernstein, Gel'fand, and Gel'fand from 1978 \cite{BGG}, 
which shows that $\dbco$ is 
equivalent to the stable module category of finitely generated left
$\lvd$-modules. We discuss this more at the end of this section. Also 
Beilinson \cite{Bei} gave in 1978 a description of $\dbco$ in terms of 
complexes of free $\sv$-modules, which we discuss in Subection \ref{4BeiSek}.

   The interesting aspect of our description compared to the original of 
   Bernstein, Gel'fand, and Gel'fand is the explicit way it is related
   to the cohomology of coherent sheaves. Namely if $\gF$ is 
a coherent sheaf on $\pv$, we show that the corresponding object in
$\kolvdcf$ is a minimal complex $I$ with $p$'th component
\[ I^p = \oplus_{ i = 0}^{\pdim} \lv(p-i) \te H^i \gF(p-i). \]
This seems to be a quite useful way in computing cohomology of coherent sheaves.
Our description also has several other advantages which will be explored
in this and the subsequent sections.

\subsection{ Compositions of $\til$ and $\fl$. }

Recall that if $X$ is a complex then the stupid truncations
$\sigma^{\geq n} X $ and $\sigma^{<n} X$ are the complexes
\begin{eqnarray*} 
\cdots \pil & 0 \pil & X^n \pil X^{n+1} \pil \cdots \\
\cdots \pil X^{n-2} \pil & X^{n-1} \pil & 0 \pil \cdots 
\end{eqnarray*}
respectively.

\begin{lemma} \label{3sigtil}
Let $I$ be in $\komolvdcf$ and let $\sigma^{\geq n} I$
be the stupid truncation.

a. The natural map $\fl((\ker d^n_I) [-n]) \pil \fl(\sigma^{\geq n} I)$ is
a homotopy equivalence.

b. The map $\tilfl (\sigma^{\geq n} I) \pil \tilfl (I)$ is
a quasi-isomorphism.
\end{lemma}

\begin{proof}
Since the only cohomology of $\sigma^{\geq n} I$ is $\ker d^n_I$, part {\it a.}
follows from Proposition \ref{1filt}. Now we prove {\it b.} 
There is an exact sequence 
\[ 0 \pil \sigma^{\geq n} I \pil I \pil \sigma^{<n} I \pil 0 \]
giving an exact sequence
\begin{equation}  0 \pil \fl(\sigma^{\geq n} I) \mto{\alpha} \fl(I) 
\pil \fl(\sigma^{<n}I) \pil 0. \label{3fsigI} \end{equation}
Now $H^p(\fl(\sigma^{<n} I))_q = H^{p+q} (\Hom_{\lvd}(k, \sigma^{<n} I))_{-q}$
by Example \ref{1kFG}. Fix any $p$. Then for $q \gg 0$ we see that 
$H^p(\fl(\sigma^{<n}I))_q = 0$. But this means that 
$\tilfl (\sigma^{<n} I)$ becomes acyclic. Hence sheafifying the map $\alpha$ in 
(\ref{3fsigI}) gives a quasi-isomorphism.
\end{proof}

\begin{proposition} \label{3Inullht}
Let $I$ be in $\komolvdcf$. If $\tilfl (I)$ is
acyclic then $I$ is nullhomotopic.
\end{proposition}

\begin{proof} Choose an integer $n$. 
By the assumption above and the previous Lemma \ref{3sigtil},
$\tilfl ((\ker d^n_I)[-n])$ is acyclic. Now we know the following fact.
If $\gA \pil \gB \pil \gC $ is a complex of coherent sheaves, exact in the 
middle, then $\gstar(\gA) \pil \gstar(\gB) \pil \gstar(\gC) $ 
is exact in the middle in sufficiently high degrees.

Thus since $\tilfl ((\ker d_I^n)[-n])$ 
is a bounded complex of coherent sheaves,
there is a $q_0$ such that $\fl((\ker d^n_I)[-n])_q$ 
is acyclic for $q \geq q_0$.
Let 
\[ J = \gsmin \circ \fl((\ker d^n_I)[-n]).\]
By Proposition \ref{1filt}, $J$ has a cofiltration
\[ \cdots J\vvi{r-1}\hvi \spil J\vvi{r}\hvi \spil J\vvi{r+1}\hvi \spil  \cdots \]
where 
\[ \ker (J\vvi{r}\hvi \spil J\vvi{r+1}\hvi) = \gs (H^r (\fl ( (\ker d^n_I)[-n])))[-r]. \] 
We then see that $J^p = 0$ for $p \gg 0$. We now have
homotopy equivalences
\[ \sigma^{ \geq n } I \iso \gs \circ \fl (\sigma^{\geq n} I) \iso 
\gs \circ \fl ((\ker d^n_I)[-n]) \iso J. \]
The cone $C(\alpha)$ of the composition $\sigma^{\geq n} I \mto{\alpha} J$ 
is then nullhomotopic. Since $C(\alpha)^p = (\sigma^{\geq n} I)^p$ for
$p \gg 0$ we get that $\im d^p_I \pil I^{p+1}$ is a split injection for
$p \gg 0$ and so $\im d^p_I$ is injective (and thus also projective) .
Recall the truncation functors $\tau^{\leq p}$ and $\tau^{ > p}$ from
\cite[1.2.7]{We}. They give an exact sequence
\begin{equation} 0 \pil \tau^{\leq p} I \pil I \pil \tau^{>p} I \pil 0
\label{3tr} \end{equation}
where $\tau^{>p} I$ is the complex
\[ 0 \pil \im d^p_I \pil I^{p+1} \pil I^{p+2} \pil \cdots \]
which is an acyclic bounded below complex consisting of injectives and is thus 
nullhomotopic.
Similarly $\tau^{ \leq p} I$ is an acyclic bounded above complex
consisting of projectives and thus nullhomotopic. Since the sequence
(\ref{3tr}) is componentwise split exact, $I$ becomes nullhomotopic
(see \cite[Sec. I.4]{Iv}.
\end{proof}

\subsection { The equivalence of categories. }

Recall from Proposition \ref{2Pderekv} 
that there is an equivalence of categories 
\[ \dbco \iso \dcq . \]
Also recall the functor $\trgstar$ in (\ref{1rgstar}).
Together with the explicit description in Theorem \ref{3ekspl} 
the following is our main result.

\begin{theorem} \label{3ekv}
There is a functor
\[ \gsmin  \circ \trgstar \, : \, \dcq \lpil \kolvdcf  \]
and a functor 
\[ \til \circ \fl \, : \, \kolvdcf \lpil \dcq. \]
These functors give an adjoint equivalence of triangulated categories,
with $ \til \circ \fl$ left adjoint. 
(Thus we get an equivalence of categories,
by \cite[IV.4.1]{MacCW}. )
\end{theorem}

\begin{proof} We first show that these functors are well-defined.
Firstly, if $\gQ$ is in $\dcq$ then $\gsmin \circ \trgstar(\gQ)$ is a complex
consisting of cofree modules of finite corank, by Lemma \ref{2begfin}.
It is also homotopy equivalent to 
$\gs \circ \rgstar(\gQ) $ by Proposition \ref{1filt}, 
and the latter is acyclic by 
Proposition \ref{2acext}. Hence $\gsmin \circ \trgstar(\gQ)$ is in $\kolvdcf$. 

Let $I$ be in $\kolvdcf$. By Lemma \ref{3sigtil} we have a quasi-isomorphism 
\[ \tilfl ((\ker d^n_I)[-n]) \pil \tilfl(I). \]
Since $\fl( (\ker d^n_I)[-n])$ is a bounded complex of finitely generated
$\sv$-modules, we see that $\tilfl(I)$ is in $\dcq$.

\medskip

Now we show that the functors are adjoint. Let $\gQ$ be in $\dcq$ and
let $\gI$ denote a bounded below injective resolution of $\gQ$ of
quasi-coherent sheaves. Let $I$ be a complex in $\kolvdcf$. Then
\begin{eqnarray*} 
\Hom_{\dcq}(\tilfl (I), \gQ) & \iso & 
\Hom_{K(\qc)}(\tilfl(I), \gI) \\
 & \iso & \Hom_{K(\sv)} (\fl(I) , \gstar(\gI)) \\
 & \iso & \Hom_{K(\lvd)} ( I, \gs \circ \gstar(\gI)) \\
 & \iso & \Hom_{\kolvdcf}(I, \gsmin \circ \trgstar(\gQ)).
 \end{eqnarray*}
 The first isomorphism is by \cite[Lemma I.4.5]{HaRD}, 
the second isomorphism by 
 Lemma \ref{2ggsadj}, and 
 the third by the adjunction (\ref{1l-s}) in Subsection \ref{1ekvSek}. 
  Next we need to show that the natural morphisms
 \begin{gather}
   (\til \circ \fl) \circ  (\gsmin \circ \trgstar) (\gQ) \lpil \gQ  
   \label{3natv} \\
   I \lpil (\gsmin \circ \trgstar) \circ (\til \circ \fl) (I) \label{3nath}
  \end{gather}
  coming from the adjunction, are isomorphisms.
Start with a $\gQ$ in $\dcq$ with bounded below injective resolution $\gI$.
Let $I$ be $\gsmin \circ \trgstar(\gQ)$. Then there is a homotopy equivalence
\[ \fl(I) \lpil \fl \circ \gs \circ \gstar(\gI) \] and a quasi-isomorphism
\[ \fl \circ \gs \circ \gstar (\gI) \pil \gstar(\gI). \]
Hence there is a quasi-isomorphism $\fl (I) \pil \gstar(\gI)$. If we 
sheafify this, we get a quasi-isomorphism 
$\tilfl (I) \pil \gI$. This shows that (\ref{3natv}) is an 
isomorphism.

Now start with an $I$ in $\kolvdcf$ and let $\gQ = \tilfl(I)$.
By the adjunction the identity $\tilfl(I) \mto{=} \gQ$
corresponds to a map $I \mto{\alpha} \gsmin \circ \trgstar(\gQ)$.
Let $C(\alpha)$ be the cone of $\alpha$. 
We must show that $\alpha$ is a homotopy equivalence. The map
$\alpha$ gives a triangle in 
$\kolvdcf$

\[ I \mto{\alpha} \gsmin \circ \trgstar(\gQ) \pil C(\alpha) \pil I[1]. \]
This gives a triangle in $K(\svf)$
\begin{eqnarray*} 
& & \fl(I) \mto{\fl(\alpha)} \fl \circ \gsmin \circ \trgstar (\gQ)
 \lpil  \fl(C(\alpha)) \\
& \lpil & \fl (I) [1] \end{eqnarray*}
and thus a triangle in $\dcq$
\begin{eqnarray*}  
& & \tilfl (I) \mto{\tilfl (\alpha)} \tilfl \circ \gsmin 
\circ \trgstar(\gQ) \lpil  \tilfl(C(\alpha)) \\
& \lpil & \tilfl (I) [1]. 
\end{eqnarray*}
Since the natural map 
\[\fl \circ \gsmin \circ \trgstar(\gQ) \pil 
\trgstar(\gQ) \] 
is a quasi-isomorphism and $\til \circ \trgstar(\gQ) \iso \gQ$, we get that
$\tilfl(\alpha)$ is a quasi-isomorphism.

Then $\tilfl(C(\alpha))$ is acyclic and by 
Proposition \ref{3Inullht} we get that 
$C(\alpha)$ is nullhomotopic. Hence $\alpha$ is an isomorphism in $\kolvdcf$
and thus (\ref{3nath}) is an isomorphism.
\end{proof}

\nota \llabel{3ytrev}
It will be convenient to have a more compact notation for the functor
$\gsmin \circ \trgstar$. For an object $\gQ$ in $\dcq$ we let its
{\it Tate resolution} be $\ytrev (\gQ) = \gsmin \circ \trgstar(\gQ)$.
\notafin

\begin{corollary} There is an equivalence of triangulated categories
\[ \dbco \bihom{ \ytrev  }{ \til \circ \fl \circ \ker d^0 }\kolvdcf . \]
\end{corollary}

\begin{proof} This is because the natural map
\[ \tilfl (\ker d^0_I ) \pil \tilfl(I) \]
is a quasi-isomorphism by Lemma \ref{3sigtil}.
\end{proof}

\subsection {Cohomology of coherent sheaves and the exterior complex. }
Now we shall be considerably more explicit about how the complex 
$\ytrev (\gF)$ looks. Recall from Section \ref{1filtsek} that there is an 
isomorphism of categories
\[  \vmod{\sv} \mto{\gs} \KomLin (\lvdcf). \]
If $M$ is an $\sv$-module then $\gs(M)$ is the complex
\[ \ldots \pil \lv(-1) \te M_{-1} \pil \lv \te M_0 \pil \lv(1) \te M_1 
\pil \cdots \]
and the maps 
\[ \lv(p) \te M_p \pil \lv(p+1) \te M_{p+1} \]
in degree $-p-1$, which is $V \te M_p \pil M_{p+1}$, 
is just the multiplication 
map for the $\sv$-module $M$.

\begin{theorem} \label{3ekspl}
 Let $\gF$ be coherent sheaf. Then $\ytrev (\gF) $ is a 
minimal complex with
\[ \ytrev (\gF)^p \iso \oplus_{r=0}^{\pdim}  
\lv(p-r) \te H^r\gF(p-r). \]
\end{theorem}

\begin{proof}
By Proposition \ref{1filt} there is 
a bounded above cofiltration 
\[ \cdots \spil \ytrev (\gF)\vvi{r-1}\hvi  \spil \ytrev(\gF)\vvi{r}\hvi  \spil 
\ytrev(\gF)\vvi{r+1}\hvi  
\spil \cdots \]
with the kernel of $\ytrev (\gF)\vvi{r}\hvi \spil \ytrev
(\gF)\vvi{r+1}\hvi $ 
equal to $\gs (H^r_* \gF) [-r]$. 
Since the components in the kernel complexes are all
injectives the statement follows.
\end{proof}

\remark  \llabel{3remkoh}
Since $H^p\gF(q) = 0$ for $p > 0$ and $q \gg 0$ we see that 
$\ytrev(\gF)^p = \lv(p) \te H^0\gF(p)$ for $p \gg 0$ when $\gF$ is a 
coherent sheaf. 

This will now enable us to describe the objects in $\kolvdcf$ which 
correspond to coherent sheaves.
\remarkfin

\defi Let $\kocolvdcf$ be the full subcategory of $\kolvdcf$ consisting of 
objects $I$ such that $I^p = \lv(p) \te  V^p_{-p}$ for $p \gg 0$.
\defifin

\begin{corollary} There is an equivalence of categories
\[ \coh \bihom{ \gsmin \circ \trgstar }{ H^0 \circ \til \circ \fl } 
\kocolvdcf . 
\]
\end{corollary}

\begin{proof} If $\gF$ is a coherent sheaf then $\gsmin \circ
\trgstar(\gF)$ is 
in $\kocolvdcf$ by the above Remark \ref{3remkoh}. 
If $I$ is in $\kocolvdcf$ then by 
Example \ref{1kFG}
\[ H^p (\fl(I))_q = H^{p+q} \Hom_{\lvd}(k, I)_{-q}. \]
If $p \neq 0$ then the latter is zero for $q \gg 0$ and hence 
$\tilfl(I)$ has zero cohomology except in the component of degree $0$. 
But then $H^0(\tilfl (I))$ and $\tilfl(I)$ are isomorphic in $\dcq$.
\end{proof}

\remark \label{3dtilG}
Let $I$ be in $\kolvdcf$. Consider a differential $I^p \mto{d^p_I} 
I^{p+1}$. Since $\lv$ is a cofree left $\lvd$-module, $\sigma^{\geq p+2}I$ 
is a cofree (injective) resolution of $\coker d^p_I$ and is thus uniquely 
determined up to homotopy. Since $\lv$ is a free left $\lvd$-module, 
$\sigma^{ \leq p-1} I$ is a free (projective) resolution of $\ker d^p_I$, and 
is also uniquely determined up to homotopy. Thus $I$ is 
{\it uniquely determined}, 
up to homotopy, by {\it any} of its differentials.

Conversely, if we have given a map
\[ \oplus_{q \in \hele} \lv(q) \te V^0_{-q} \mto{d} \oplus_{q \in 
\hele} \lv(q) \te V^1_{-q} \]
with $\sum_{q \in \hele} \dim_k V^0_{-q}$ and $\sum_{q \in \hele} \dim_k 
V^1_{-q}$ finite, then by taking a cofree resolution of $\coker d$ and a 
free resolution of $\ker d$ with components of finite corank and rank, 
we get an object $I$ in $\kolvdcf$ and thus an 
object in $\dbco$.

\medskip
This gives a great amount of freedom in constructing objects of 
\linebreak $\dbco$.

\medskip
If one  by being a bit clever or maybe lucky in choosing $d$ so that one 
knows that one actually has constructed a coherent sheaf (see Subection 
\ref{6EtilFSek} for 
more on this), then by Theorem \ref{3ekspl} one can compute all the cohomology 
modules of $\gF$ by taking a minimal cofree resolution of $\coker d$ and a 
minimal free resolution of $\ker d$.

On Macaulay 2, version 0.86, there exists such procedures for the
exterior algebra.  

\remarkfin

\eks If $m$ is greater or equal to the regularity of $\gF$ (see Subsection 
\ref{4mreg} 
for more on this), then $H^r\gF (m-r) = 0$ for $r > 0$. Thus by Theorem 
\ref{3ekspl} the differential $d^m_I$ is given by
\[ d^m_I \, : \, \lv(m) \te H^0\gF(m) \pil \lv(m+1) \te H^0\gF(m+1) . \]
This differential is determined by the map in degree $-m-1$ which is just the 
multiplication map \[ W \te H^0\gF(m) \pil H^0\gF(m+1). \]
By taking a minimal free resolution of $d^m_I$ one can compute all the 
cohomology groups  $H^r \gF(n)$ for integers $r \geq 0$ and $n$.
\eksfin

\rem \llabel{3barth} 
In the paper \cite{Ba} by Barth the main result says that if
$\gF$ is a stable rank two sheaf on $\pv = \pto$ with first Chern
class $c_1(\gF) = 0$, then $\gF$ is completely determined by the map
\[ W \te H^1\gF(-2) \lpil H^1\gF(-1). \]
But such a sheaf has $H^0\gF$ and $H^2\gF(-2)$ equal to zero.
Also $H^0\gF(-1)$ and $H^2\gF(-3)$ are zero. Hence the exterior
complex in the components of degree $-1$ and $0$ is \[ \lv(-2) \te H^1\gF(-2) \lpil \lv(-1) \te H^1\gF(-1). \]
and so the result of Barth follows from this. 
\remfin

\subsection{ Duals of Tate resoltions}
Note that $\Hom_{\lvd} (\lv, \wedge^{v+1} W)$ is canonically 
isomorphic to $\lvd \te \wedge^{v+1} W$ an hence to $\lv$. 
Thus the dual of a Tate resolution 
$\Hom_{\lvd}(I, \wedge^{v+1}W)$ is also naturally a Tate resolution.

In fact if $I = T(\gG)$, it is the Tate resolution of 
$\rderhom (\gG, \omega_{\pv})[-v-1]$ as the following shows.
Thus we retrieve the usual statement of Grothendieck duality.

\begin{proposition} The diagram

\[ \begin{CD}
\kolvdcf @>{\til \circ \ker^0}>>  \dbvb \\
@V{\Hom_{\lvd}(-, \wedge^{v+1} W)}VV  @VV{\homkn(-, \omega_{\pv})[-v-1]}V \\
\kolvdcf @>{\til \circ \ker^0}>> \dbvb. 
\end{CD} \]
induces a natural isomorphism of functors betweent the two compositions.
Hence for $\gG$ in $\dbco$ the dual $\Hom_{\lvd}(T(\gG), \wedge^{v+1} W)$
is the Tate resolution of 
$ \rderhom(\gG, \omega_{\pv})[-v-1].$
\end{proposition}

\begin{proof} This is straightforward.
\end{proof}

\subsection{ The description of Bernstein, Gel'fand, and Gel'fand. }

The description we have given of $\dbco$ is closely related to the 
description given in \cite{BGG}. They show that $\dbco$ is equivalent to the 
stable module category of $\lvd$. We recall what this is and describe how 
their result it related to ours.

Let $A$ be a positively graded (associative, but not necessarily commutative) 
Artin algebra. For objects $M$ and $N$ in $\vmod {A}$, let $P(M,N) \sus 
\Hom_{\vmod {A}}(M,N)$ be the subset of all morphisms $M 
\pil N$ which factors as
$ M \pil P \pil N$ 
where $P$ is a projective $A$-module.

The {\it stable module category} $\modstab{A}$ is the quotient category of 
the category of finitely generated graded 
$A$-modules, $\vfmod{A}$, having the same objects as $\vfmod {A}$ but with
\[ \Hom_{\modstab{A}}(M,N) = \Hom_{\vfmod{A}}(M,N)/ P(M,N) . \]

When $A$ is a Frobenius algebra (see end of Subsection \ref{1lv}),  
$\modstab{A}$ is 
(equivalent to) a triangulated category. The translation functor $T$ is (up to 
isomorphism) defined as follows.
For $M$ in $\modstab{A}$ embed $M$ in an injective (= projective) object $I$. 
Then the translation $T(M)$ is the cokernel of $M \inpil I$.

Let $\koacf$ be the category consisting of all acyclic complexes of
cofree (= free) left $A$ modules of finite corank (= finite rank), 
with morphisms up to homotopy. 
The following proposition is well known and we state it without proof.

\begin{proposition}
There are functors
\begin{eqnarray*}
\ker d^0 \, : \, \koacf &\lpil &\modstab{A} \\
 I & \mapsto & \ker d^0_I 
\end{eqnarray*}
and 
\begin{eqnarray*}
\resol \, : \, \modstab{A} & \lpil & \koacf \\
M & \mapsto & I
\end{eqnarray*}
where the $I$ is constructed  by letting
$ M \pil I^0 \pil I^1 \pil \cdots$ be an injective resolution of $M$ and $ 
\cdots \pil I^{-2} \pil I^{-1} \pil M$ be a projective resolution.

These functors give an equivalence of triangulated categories.
\end{proposition}

\begin{corollary}[\cite{BGG}]
There is an equivalence of triangulated categories
\[ \dbco \iso \modstab{ \lvd}. \]
\end{corollary}

\section { Symmetric complexes representing complexes of coherent sheaves. }

We use the term symmetric complex to denote a complex of free
$\sv$-modules.

\medskip

If $\gF$ is a coherent sheaf on $\pv$, then a common way to represent it, is 
by means of a minimal free resolution of the $\sv$-module $H^0_* \gF$. 
However this is but one of many ways to represent $\gF$ by complexes of free 
$\sv$-modules. For instance, if $\gF = \gI_{C / \ptre}$ which is 
the ideal sheaf of a 
space curve $C \sus \ptre$, then one has the minimal free resolution of 
$H^0_* \gI_{C/\ptre}$, a complex of free $\sv$-modules $A^{2} \pil A^1 \pil 
A^0 $ with $H^0(A) = H^0_* \gI_{C / \ptre}$ but one also has the minimal 
monad, a complex of free $\sv$-modules $B^{-1} \pil B^0 \pil B^1$ with $H^0(B) 
=  H^0_* \gI_{C / \ptre} $ and $H^1(B) = H^1_* \gI_{C/\ptre}$ (and $H^{-1}(B) 
= 0$). The monad had been used much in the study of space curves. These 
monads are generalized by Walter complexes, shortly to be discussed. 

Also Beilinson \cite{Bei} showed that any object $\gG$ in $\dbco$ is 
isomorphic to the sheafification of a symmetric complex $P$ such 
that each $P^p$ is a direct sum of modules $\sv(-q)$ where 
$0 \leq q \leq \pdim$. 
 
In this section we show that such representations of an object $\gG$ in 
\linebreak $\dbco$ by means of symmetric complexes may be derived by means of 
{\it truncating} 
the complex $\ytrev(\gF)$  at different places and then transforming this 
truncated complex to a symmetric complex using the functor $\flmin$.

\subsection {Symmetric complexes as transforms of truncations of the exterior 
complex. }

\begin{proposition} \label{4trunkext}
a. Let $I$ be in $\komolvdcf$ and let $J$ in 
\linebreak $\Kom(\lvdcf)$  be such that
i. $J$ is bounded below, ii. its components have finite corank, and iii. 
$\sigma^{\geq p} J = 
\sigma^{\geq p} I$ for some $p$. Then $\flmin (J)$ is a bounded complex of 
free $\sv$-modules of finite rank (displayed here with components of degree 
$0$ and $1$) 
\[ \cdots \pil \oplus_{q \in \hele} \sv(-q) \te H^q(J)_{-q} \pil \oplus_{q 
\in \hele} \sv(-q+1) \te H^{q}(J)_{q-1} \pil \cdots  \]
and $\tilflmin(J)$ is isomorphic to  $\tilfl(I)$ in $\dbco$.

b. Conversely given a bounded minimal complex $P$ of free $\sv$-modules of 
finite rank. Then there exists $J$ and $I$ fulfilling the criteria above 
together with homotopy equivalences
$P \iso \flmin(J)$ and $\overset{\til}{P} \iso
\tilflmin(I)$. 
\end{proposition}

\begin{proof} 
Condition {\it i.} implies by Proposition \ref{1filt} that 
$P = \flmin (J)$ has a filtration 
$ \cdots \sus P\vvi{r-1}\hvi \sus P\vvi{r}\hvi \sus \cdots $ 
with the cokernels $P\vvi{r}\hvi/P\vvi{r-1}\hvi$ in the filtration of the form
$\fl(H^r(J))[-r]$ which is a complex (displayed with components of degree $0$ 
and $1$)
\[ \cdots \pil \sv(-q) \te H^q(J)_{-q} \pil \sv (-q+1) \te H^q(J)_{-q+1} 
\pil \cdots. \]
 Conditions {\it i},{\it ii.}, and {\it iii.} now imply all but the 
last statement in  {\it a.} 

There is an exact sequence
\[ 0 \pil \sigma^{\geq p} J \pil J \pil \sigma^{<p} J \pil 0 \]
giving an exact sequence
\[ 0 \pil \fl(\sigma^{\geq p} J) \pil \fl(J) \pil \fl(\sigma^{<p} J) \pil 0. 
\]
Now by Example \ref{1ekskoskom}, $\fl(\lv)$ is quasi-isomorphic to $k$
and hence its sheafification $\tilfl ( \lv)$ is zero..
Since $\sigma^{<p} (J)$ is a bounded complex of cofree $\lvd$ modules of 
finite corank, and
$\fl$ is exact, it follows easily that $\tilfl (\sigma^{<p} J) \iso 0$.
Thus 
\[ \tilfl(J) \iso \tilfl(\sigma^{\geq p} J) \iso \tilfl 
(\sigma^{\geq p} I) \iso \tilfl(I), \]
the latter according to Lemma \ref{3sigtil}. Hence the last statement in {\it a.} follows.

{\it b.} Let $J = \gsmin (P)$. Then clearly {\it i.} and {\it ii.} in 
{\it a.} holds. Also
$H^p(J)_q$ is isomorphic to 
$H^{p+q}(k \te_{\sv} P)_{-q}$ by Example \ref{1kFG}. 
Hence the conditions on $P$ gives 
that $H^p(J) = 0$ for $p \geq p_0$ say. 
Taking a projective resolution of $\ker d_J^{p_0}$ we get a complex $I$ in 
$\komolvdcf$ with $\sigma^{\geq p_0} I = \sigma^{\geq p_0} J$. Also there are 
natural maps
\[ \fl(J) = \fl \circ \gsmin (P) \pil \fl \circ \gs (P) \pil P \]
which are homotopy equivalences. By the uniqueness part 
Proposition \ref{1filt} we must have a homotopy equivalence
\[ \flmin (J) \simeq P . \]
\end{proof}

\rem This result suggests that the exterior complex $\ytrev (\gF)$  of a 
coherent sheaf $\gF$,
is in some ways a more basic invariant than the symmetric complexes 
representing $\gF$ since the latter are obtained by transforming different 
"truncations" of the exterior complex $\ytrev (\gF)$ using the functor 
$\flmin$. Note however from the proof that any {\it bounded} complex $J$ 
fulfilling the conditions {\it i.} and {\it ii.} in a. will have 
$\tilflmin(J) \iso 
0$. Thus in general a complex $J$ fulfilling the conditions in a. will 
contain a lot of "junk" in the lower left end, and thus
similarly the transformed complex $\flmin(J)$ contains a lot of "junk". But 
if we perform
certain canonical truncations of an exterior complex $\ytrev (\gG)$ where 
$\gG$ is in 
$\dbcq$ and transform this to symmetric complexes, we get several well-known 
canonical symmetric complexes associated to $\gG$, as we shall now see.
\remfin

\subsection{ Linear complexes. }
We use the term rigid complex to denote bounded complexes of type
(here displayed with components of degree $0$ and $1$)
\begin{equation} \cdots \pil \gopv (-p) \te V^0 \pil \gopv(-p+1) \te
V^1 \pil \cdots \label{4rigid} \end{equation}
where $p$ is some integer and the $\dim_k V^i$ are finite.

The following says what symmetric complexes we may get by performing
the stupid truncation $\sigma^{\geq p} I$ on an exterior complex.

\begin{proposition} Let $\gG$ in $\dbco$ correspond to $I$ in 
$\kolvdcf$. 
Given an integer $p$. Then $\tilflmin(\sigma^{\geq p} I)$ is a rigid
complex of type (\ref{4rigid})
isomorphic in $\dbco$ to $\gG$.
\end{proposition}

\begin{proof} Let $I = \ytrev(\gG)$ and consider it as a complex in 
$\Kom(\lvdcf)$. 
By Lemma \ref{3sigtil} there is a quasi-isomorphism 
\[ \tilflmin(\sigma^{\geq p} I) = (\tilfl((\ker d^p_I)[-p]) \pil \tilfl(I) . \]
The middle complex is (displayed with components of degree $0$ and 
$1$)
\[ \cdots \pil \sv(-p) \te (\ker d^p_I)_{-p} \pil \sv(-p+1) \te (\ker 
d^p_I)_{-p+1} \pil \cdots, \]
and so we get the proposition by sheafifying this complex and letting 
$V^i$ be 
$(\ker d^p_I)_{-p+i}$.
\end{proof}

\subsection{ Beilinson complexes. } \label{4BeiSek}
In \cite{Bei}, Beilinson gave an alternative description of the 
category $\dbco$, representing the objects by certain complexes
of free $\sv$-modules. This is Theorem \ref{4BeiThm} below. During
his original proof of this theorem, he also established the well-known
Beilinson spectral sequence. To get the complexes of Beilinson we
perform the following truncations on an exterior complex $I$.

\defi Let $I$ be in $K(\lvdcf)$ with 
\[ I^p = \oplus_{q \in \hele} \lv(q) \te V^p_{-q}. \]
For an integer $r$ let $b^{\geq r} I$ be the subcomplex of $I^p$ given by 
\[ (b^{\geq r} I)^p = \oplus_{q \geq r} \lv(q) \te V^p_{-q}. \]
We let $b^{< r} I$ be the cokernel of $b^{\geq r} I \pil I$.
Note that these functors, in contrast to $\sigma^{ \geq r}$ and $\sigma^{< r}$,
are functors on the 
homotopy category $K(\lvdcf)$.
\defifin

If $I = \ytrev (\gG)$, which is a minimal complex, then clearly $(b^{\geq r} 
I)^p = 0$ for $p \ll 0$ and $(b^{\geq r} I)^p = I^p$ for $p \gg 0$. By 
Proposition \ref{4trunkext} a. $\flmin (b^{\geq r} I)$ exists and its sheafification is 
isomorphic to $\gG$. 

We have the exact sequence
\[ 0 \pil b^{\geq r} I \pil I \pil b^{< r} I \pil 0. \]
Since $H^p(I) = 0$ for all integers $p$, we see that $H^p(b^{\geq r} I) _{-q}$
is non-zero only if {\it i.} $\lv(r+i)_{-q} $ is non-zero for some $i \geq 0$
which implies $r \leq q$, and {\it ii.} $\lv(r-i)_{-q} $ is non-zero for some
$i > 0$ which implies $q \leq r  + \pdim$. 
%
Thus we must have $r  \leq q \leq r + \pdim$. 

So we see by Proposition \ref{4trunkext}
that $\flmin (b^{\geq r} I) $ is a complex (displayed with 
components of degree $0$ and $1$)
\[ \cdots \pil \oplus_{q = r}^{ r + \pdim} \sv(-q) \te 
H^q(b^{\geq r} I) _{-q} \pil \oplus_{ q = r} ^{r + \pdim} \sv(-q) 
\te H^{q+1}(b^{\geq r} I)_{-q} \pil \cdots .\]

\defi $K[-r - \pdim, -r]$ is the full subcategory of $K(\svf)$ consisting of 
bounded complexes $P$ with 
\[ P^p = \oplus_{q = r}^{r + \pdim} \sv(-q) \te V^p_q \]
and $\dim_k V^p_q$ finite.
\defifin


\begin{theorem}[Beilinson] 
For each integer $r$ there is an equivalence of categories
\[ \dbco \bihom{\flmin \circ b^{\geq r} \circ \ytrev }{\til} 
K[-r - \pdim, -r] \] 
(where $\sim$ is sheafification).
\label{4BeiThm}
\end{theorem}

\begin{proof}
We get a functor
\[  \flmin \circ b^{\geq r} \, : \, \kolvdcf \pil K[-r-\pdim, -r]. \]
We shall prove that this functor is fully faithful and that every object in 
$K[-r-\pdim,-r]$ is isomorphic to $\flmin(b^{\geq r} I)$ for some $I$ in 
$\kolvdcf$. This will establish that $\flmin \circ b^{\geq r} $ is an 
equivalence of categories, by \cite[II.2.7]{MaGe}.

Let $I$ and $J$ be in $\komolvdcf$. The exact sequence
\[ 0 \pil b^{\geq r} J \pil J \pil b^{< r} J \pil 0 \]
gives a triangle in $\kolvdcf$
\[ b^{< r} J[-1] \pil b^{\geq r} J \pil J \pil b^{< r} J \]
and an exact sequence
\begin{eqnarray*} & \Hom_{K(\lvdcf)} (b^{\geq r} I, b^{<r} J[-1]) & \pil 
\Hom_{K(\lvdcf)}(b^{\geq r}I, b^{\geq r} J) \\
  \pil & \Hom_{K(\lvdcf)}(b^{\geq r} 
I, J) & \pil \Hom_{K(\lvdcf)}(b^{\geq r} I, b^{<r } J) \end{eqnarray*}
which gives an isomorphism
\begin{equation}  \Hom _{K(\lvdcf)} (b^{\geq r} I, b^{\geq r } J) \ispil 
\Hom_{K(\lvdcf)} (b^{\geq r} I, J)). \label{4IJ} \end{equation}
We also get from the triangle in $\kolvdcf$
\[ b^{< r} I [-1] \pil  b^{\geq r} I \pil I \pil b^{< r} I  \]
an exact sequence
\begin{eqnarray} \label{4sekIJ}
& \Hom_{K(\lvdcf)}(b^{< r} I, J)  \pil & \Hom_{K(\lvdcf)}(I,J) \\ 
 \pil & \Hom_{K(\lvdcf)}(b^{\geq r} I, J)  
\pil & \Hom_{K(\lvdcf)}(b^{< r} I [-1], J). \notag
\end{eqnarray}
Now since $I$ is homotopy equivalent to a minimal complex, $b^{< r} I$ will be 
homotopy equivalent to a bounded above complex. Since $J$ is acyclic, 
(\ref{4sekIJ}) gives an isomorphism 
\[ \Hom_{K(\lvdcf)}(I,J) \ispil \Hom_{K(\lvdcf)}(b^{\geq r} I, J). \]
Together with (\ref{4IJ}) this gives that 
\[ b^{\geq r} \, : \, \kolvdcf  \lpil K(\lvdcf) \]
is a fully faithful functor. Also
\[ \Hom_{K(\lvdcf)}(b^{\geq r } I, b^{\geq r} J) = 
\Hom_{K(\svf)}(\fl(b^{\geq r} I), \fl(b^{\geq r} J)) \]
since $\fl$ gives an equivalence of categories. Hence $\flmin \circ b^{\geq 
r}$ is fully faithful.

\medskip

Now given $P$ in $K[-r-\pdim,-r]$. 
We want to show that $P = \flmin (b^{\geq r} I)$ for some $I$ in
$\kolvdcf$.

Let $J = \gsmin (P)$. Then $J$ is a 
bounded below minimal complex. Also 
\[ \Hom^p_{K(\lvdcf)}(k,  J)_q = 
H^p \Hom_{K(\lvdcf)}(k,J)_q = H^{p+q}(P)_{-q} . \]
So $\Hom^p_{K(\lvdcf)}(k, J)_q$ is nonzero only if $ -q \geq r$. 
This gives $b^{\geq r} J = J$. 
We now claim that $J$ can be completed to a minimal complex $I$ in
$\kolvdcf$ with $b^{\geq r} I = J$. This follows by looking at
$H^p(J)_q = H^{p+q}(k \te_{\sv} P)_{-q}$. We see that 
\begin{enumerate}[i.]
\item $H^p(J)$ is non-zero for only a finite number of $p$'s.
\item $H^p(J)$ is a finite length module.
\item $H^p(J)_q$ non-zero implies $q \geq - r - \pdim$.
\end{enumerate}
Hence $J$ can be completed to a minimal complex $I$ in
$\kolvdcf$ with $b^{\geq r } I = J$.
But then $P \iso \flmin (b^{\geq r} I)$. This shows that $\flmin 
\circ b^{\geq r} $ is an equivalence of categories.


By Proposition \ref{4trunkext} we now see that if $I$ is in $\kolvdcf$ then
\[ \tilfl (I) \iso \tilflmin(b^{\geq r} I ). \]
By the main Theorem \ref{3ekv} we see that $\til$ is a quasi-inverse
to $\flmin \circ b^{\geq r} \circ \ytrev $. 
\end{proof}

\subsection {Walter complexes. } 
The minimal free resolution of a coherent sheaf has certain
``higher'' versions introduced by C. Walter in the paper \cite{Wa}.
We describe them below in Theorem \ref{4Wa}. They derive from the 
following truncations of a minimal exterior complex.

\defi Let $\gG$ be in $\dbco$ so $\ytrev (\gG)$ is a minimal complex. By 
Proposition \ref{1filt} there is a cofiltration of $\ytrev(\gG)$.
\[ \cdots \spil \ytrev (\gG)\vvi{r}\hvi \spil \ytrev (\gG)\vvi{r+1}\hvi \spil \cdots . \]
Let 
\[ w_{\leq r} \ytrev (\gG)  = 
\ker ( \ytrev (\gG) \pil \ytrev (\gG)\vvi{r}\hvi ). \]
\defifin

If $K\vvi{r}\hvi$ is the kernel of $\ytrev (\gG)\vvi{r}\hvi 
\pil \ytrev (\gG)\vvi{r+1}\hvi $ then 
$K\vvi{r}\hvi$ is the complex 
$\gs (H^r(\rgstar(\gG))) [-r]$ which is a complex 
(displayed with components of degree $0$ and $1$)
\[ \cdots \pil \lv(-r) \te H^r(\rgstar(\gG))_{-r} \pil \lv (-r+1) \te 
H^r(\rgstar(\gG))_{-r+1} \pil \cdots . \]
Thus 
\begin{equation} \label{4wrE}
 (w_{\leq r} \ytrev (\gG))^p = \oplus_{i \leq r} \lv(p-i) \te H^i (\rgstar 
(\gG))_{p-i} .
\end{equation}

Of course we would now like to form the complex 
\[ \flmin (w_{\leq r} 
\ytrev (\gG)) \]
such that its sheafification is isomorphic to $\gG$ in $\dbco$.
The problem is that $w_{\leq r} \ytrev (\gG)$ is not necessarily bounded 
below. Nor needs
$(w_{\leq r} \ytrev (\gG))^p$ be equal to $\ytrev (\gG)^p$ for $p \gg 0$. 
Thus  $\fl (w_{\leq r} \ytrev (\gG))$ 
needs not be homotopic to a minimal complex, and even if it were,  its 
sheafification needs not be isomorphic to $\gG$ in $\dbco$.

However in certain cases the desired properties hold. If $\gG$ is a coherent
sheaf $\gF$ on $\pv$, recall that the {\it local projective dimension} of 
$\gF$, written $\lpd(\gF)$, is the maximum of all projective dimensions of 
$\gF_P$, the localization of $\gF$ at a point $P$ of $\pv$, as a module over 
the local ring $\gopvp$.

It is a well known fact that $\lpd(\gF) \leq k$ iff $H^p \gF(q)$ vanishes
when $q \ll 0$ for $0 \leq p \leq \pdim - 1 - k$. (This is true if $\gF$ is a 
vector bundle, by Serre duality; then take a locally free resolution of 
$\gF$. ) Of course $H^p \gF(q)$ also vanishes for $q \gg 0$ 
when $p > 0$. 
 

\begin{theorem} \label{4Wa}
Let $\gF$ be a coherent sheaf. Assume that $0 \leq r \leq 
\pdim - 1 - \lpd(\gF)$. Then $w_{\leq r} \ytrev (\gF)$ is a bounded below 
complex so the symmetric complex $P = \flmin (w_{\leq r} \ytrev( \gF))$
exists. It has the following properties.
\begin{enumerate}[i.]
\item The sheafification of ${P}$ is quasi-isomorphic to $\gF$.
\item $P^p = 0$ if $p$ is not in the interval $[ r+1-\pdim, r]$, so $P$ 
has length
at most $\pdim-1$. 
\item $H^p(P) \iso H^p_* \gF$ for $p \leq r$.
\end{enumerate}
\end{theorem}

\remark Note that the critical thing which makes such a complex unique
is that the length of $P$ is at most $\pdim - 1$.
If $r=0$ then $P$ is just the minimal free resolution of $H^0_* \gF$.
\remarkfin

\begin{proof}[Proof of Theorem \ref{4Wa}]  Let $I = \ytrev (\gF)$.
By (\ref{4wrE}) and what is stated just before this proposition,
it is clear that $(w_{\leq r} I)^p = I^p $ for $p \gg 0$ if $r \geq 0$ and 
that $(w_{\leq r} I)^p = 0$ for $p \ll 0$ if $0 \leq r \leq \pdim - 1 - 
\lpd(\gF)$. Thus $\flmin (I)$ exists and its sheafification is 
quasi-isomorphic to $\gF$ by
Proposition \ref{4trunkext}.

Now we have the short exact sequence
\begin{equation} 0 \pil w_{\leq r} I \pil I \pil I / w_{\leq r} I \pil 0. 
\label{4wI} \end{equation}
We see by (\ref{4wrE}), that if $H^p(w_{\leq r} I)_q$ is non-zero then 
$ p -i +q \leq 0 $ for some $i \leq r$. This implies $p+q 
\leq r$. Furthermore we see by 
(\ref{4wrE}) that $H^{p-1} (I / w_{\leq r}I)_q = H^p (w_{\leq r} I)_q$ is 
non-zero only if
$p-1-i +q \geq - \pdim - 1$ for some $i > r$ which implies 
$p+q \geq r+1 - \pdim$.
Thus we get that $H^{p+q} (k \te_{\sv} P) _{-q} = H^p(w_{\leq r} I )_q$ is 
zero for $p+q$ not in the interval $[r+1 - \pdim, r]$. Since $P$ is minimal 
this gives $P^{p+q} = 0$ for $p+q$ not in the interval $[r+1-\pdim, r]$. 

Let $P = \flmin(w^{\leq r} I)$. Then 
$w^{\leq r} I = \gsmin(P)$ and by Proposition \ref{1filt} a. we get that the
kernels $K\vvi{p}\hvi$ in the filtration of $w^{\leq r} I$ are isomorphic to
$\gs ( H^p(P))[-p]$. 


But also  $I = \gsmin \circ \trgstar(\gF)$. 
Thus in the cofiltration of $w_{\leq r} I$ the kernels
$K\vvi{p}\hvi = \ker (I\vvi{p}\hvi \spil I\vvi{p+1}\hvi)$ for $p \leq r$ are equal to $\gs (H^p_* 
\gF)[-p]$.
Hence we get $H^p(P) \iso H^p(\rgstar(\gF))$ for $p \leq r$
since the functor (\ref{1gskomlin}) in Subsection \ref{1filtsek}
is an isomorphism of categories.

\end{proof}

\remark
Suppose $C \sus \ptre$ is a space curve (without isolated or embedded 
points). Then let $\gF$ be  $ \gI_{C/ \ptre}$, the ideal sheaf of the curve in 
$\ptre$. Then $\lpd(\gI_{C/ \ptre}) = 1 $.
Thus $\flmin ( w_{\leq r} \ytrev(\gI_{C/ \ptre}))$ exists for $r = 0$ or $1$ 
and it is a complex of length $2$. If $r = 0$ we get the minimal resolution 
$P^{-2} \pil P^{-1} \pil P^0$ of 
$H^0_* \gI_{C/ \ptre}$. If $r = 1$ we get the minimal monad $P^{-1} \pil P^0 
\pil P^1$ with   
$H^0(P) = H^0_* \gI_{C/ \ptre}$ and $H^1(P) = H^1_* \gI_{C/ \ptre}$. The 
minimal monad has been much studied in the field of space curves.
\remarkfin

\subsection { Koszul cohomology.} \label{4koskohSek}
Consider the complex $w_{\leq 0} \ytrev(\gF)$
(displayed with components of degree $p$ and $p+1$)
\[ \cdots  \pil \lv(p) \te H^0 \gF(p) \pil 
\lv(p+1) \te H^0 \gF(p+1) \pil \cdots . \]
In \cite{Gr} the cohomology $H^p(w_{\leq 0} \ytrev(\gF))_q$ of this 
complex is denoted \linebreak 
$\gK_{-p-q,p}(H^0_* \gF,W)$ and called Koszul
cohomology groups. 

On the other hand consider the minimal free resolution $P$ of $H^0_*
\gF$
\begin{equation*}
 \cdots \pil \oplus_{q \in \hele} S(W)(-q) \te V_q^p \pil \cdots \pil
\oplus_{q \in \hele} S(W)(-q) \te V_q^0 \pil 
 H^0_* \gF.
\end{equation*}
The syzygies $V^p_q$ of order $p$ and weight $q$ are often of interest
because they are attached geometric significance. The following is
\cite[Thm. 1.b.4]{Gr} and shows that the syzygies, suitably reindexed,
are isomorphic to the Koszul cohomology groups.

\begin{proposition} \label{4koskoh}
Let $\gF$ be a coherent sheaf with $\lpd(\gF) \leq
\pdim - 1$ and let $P$ be a minimal free resolution of $H^0_* \gF$. 
Then the syzygy of order $p$ and weight $q$,  
$H^p(k \te_{\sv} P)_q$, is isomorphic as a vector space to $H^{p+q} 
(w_{\leq 0} \ytrev (\gF))_{-q}$. 
(Or, in other notation, $V^p_q$ is isomorphic to $\gK_{-p, p+q}
(H^0_* \gF, W)$.)
\end{proposition}

\begin{proof} This follows since $P$ is  
$\flmin  (w_{\leq 0} \ytrev (\gF))$ and 
thus by Example \ref{1kFG}, 
$H^p (k \te_{\sv} P)_q$ is isomorphic to $H^{p+q} (w_{\leq 0} \ytrev 
(\gF))_{-q}$.
\end{proof}

\subsection{Castelnuovo-Mumford regularity. } \label{4mreg}

A result where the form of the exterior complex (see Theorem \ref{3ekspl}) 
of a coherent sheaf  $\gF$ is to a certain extent present, 
is in the concept of 
$m$-regularity of a coherent sheaf as defined originally in 
Mumford's book \cite{Mu}. Recall that 
the component $\ytrev (\gF)^p$ is equal to 
$\oplus_{i = 0}^{\pdim} \lv (p-i) \te H^i  \gF (p-i)$.
The following is the original theorem in \cite[Lec. 14]{Mu}.

\begin{proposition} Given an integer $m$. 
Suppose $H^i \gF(m-i) = 0$ for $i > 0$. Then 
\begin{enumerate}[i.]
\item $H^0_* \gF$ is generated in degree $m$.
\item $H^i \gF(p-i) = 0$ for all $p \geq m$ and $i > 0$.
\end{enumerate}
\end{proposition}

\begin{proof}
The component $\ytrev(\gF)^m$ is $\lv(m) \te H^0 \gF (m)$. Suppose 
that $H^i \gF(p-i) $ is non-zero for some $i > 0$ and $p >m$ and let $p$ be 
minimal such. 
Since $\ytrev(\gF)$ is a minimal complex, 
$\lv(p-i)_{i-p} \te H^i \gF(p-i)$ is in the 
kernel of the differential $d^p_I$. But again since $\ytrev(\gF)$ 
is minimal this is 
not in the image of $d^{p-1}_I$ whose domain is
$\lv(p-1) \te H^0 \gF(p-1)$. Since $\ytrev(\gF)$ is 
acyclic this must mean that $H^i \gF(p-i) = 0$ for all 
$i > 0$ and $p \geq m$ and this proves {\it ii.}

Thus $\ytrev(\gF)^p = \lv (p) \te H^0 \gF(p)$ for $p \geq m$ so the complex
\[ \cdots \pil \lv(p) \te H^0 \gF(p) \pil \lv(p+1) \te H^0 \gF(p+1) \pil 
\cdots \]
is acyclic for $p > m$. Thus 
\[ \lv(p-1)_{-p} \te H^0 \gF(p-1) \pil \lv(p)_{-p} \te H^0 \gF(p)  \] 
is surjective. But this is just the map
\[ W \te  H^0 \gF(p-1) \pil H^0 \gF (p)  \]
and this proves {\it i.}   
\end{proof}

\section {Hilbert polynomials. }

Suppose we have given the exterior complex $\ytrev (\gF)$ of a coherent sheaf  
$\gF$. In this section we show how one may compute the Hilbert polynomial of 
$\gF$ from $\ytrev (\gF)$. As we shall see, this may be done quite locally on 
$\ytrev (\gF)$, in fact from the kernel of any of the differentials of $\ytrev 
(\gF)$.

\medskip
\subsection {Hilbert polynomials.}
Let $\gF$ be a coherent sheaf. Then the Hilbert polynomial 
$\chi \gF$ of $\gF$ is 
\[ \chi \gF (n) = \sum_{i \geq 0} (-1)^i \dim_k H^i \gF(n) . \]
For a bounded complex $\gG$ of coherent sheaves we define the Hilbert 
polynomial to be 
\[ \chi \gG(n) = \sum_{p \in \hele} (-1)^p \chi \gG^p(n) \]
where $\gG^p$ is the coherent sheaf in component $p$ of $\gG$.

If 
\[ 0 \pil \gF^\prime \pil \gF \pil \gF^{\prime \prime} \pil 0 \]
is a short exact sequence of coherent sheaves, then $\chi \gF = \chi 
\gF^\prime + \chi \gF^{\prime \prime} $. 
Thus for every $\gG$ in $\dbco$ we get a polynomial function
$\chi \gG : \hele \pil \hele$
which is additive on triangles, i.e. if 
\[ \gG^\prime \pil \gG \pil \gG^{\prime \prime} \pil \gG^\prime [1] \]
is a triangle, then 
\[ \chi \gG = \chi \gG^\prime + \chi \gG^{\prime \prime}. \]

\begin{lemma}
Let $\gG$ be in $\dbco$. 
Then 
\[ \chi \gG(n) = \sum_{p \in \hele} (-1)^p \dim_k H^p (\rgstar(\gG))_n . \]
\end{lemma}

\begin{proof} If $\gG$ is a coherent sheaf, this is just the definition of 
$\chi \gG$. 
Let 
\[ \gG^\prime \pil \gG \pil \gG^{\prime \prime} \pil \gG^{\prime}[1] \]
be a triangle in $\dbco$.
The functor $\rgstar$ takes triangles to triangles. Taking cohomology we 
then get a long exact sequence of cohomology. Thus if the lemma holds for 
$\gG^\prime$ and $\gG^{\prime \prime}$ it holds for $\gG$. Since the coherent 
sheaves generate the triangulated category $\dbco$, we get the lemma.
\end{proof}

The following is the most useful result for computing the Hilbert 
polynomial of a coherent sheaf from its exterior complex. It shows
how this Hilbert polynomial can be computed from any of the
differentials of the exterior complex.

\begin{theorem} \label{5ker}
Let $I$ in $\kolvdcf$ correspond to $\gG$ in $\dbco$. 
Fix an integer $p$. Then
\[ \chi \gG(n) = \sum_{q \in \hele} (-1)^{q+p} \dim_k (\ker d_{I}^p)_q \cdot 
\chi \gopv(q+n). \]
\end{theorem}

\begin{proof} By Lemma \ref{3sigtil}, 
$ \tilfl((\ker d_{I}^p) [-p]) $ is isomorphic to $\gG$ in $\dbco$. 
Since $\fl ((\ker d^p_{I}) [-p])$ is the complex (displayed with components 
of degree $0$ and $1$)
\[ \cdots \pil \sv(-p) \te \ker (d^p_{I})_{-p} \pil \sv(-p+1) \te \ker 
(d^p_{I})_{-p+1} \pil \cdots \]
we get that 
\[ \chi \gG(n) = \sum_{q \in \hele} (-1)^{q+p} \dim_k \ker (d_{I}^p)_{q} \cdot 
\chi \gopv(q+n). \]
\end{proof}

The following shows how to compute the Hilbert polynomial of a coherent
sheaf from the coranks of the components of the exterior complex.

\begin{theorem} Let $I$ in $\kolvdcf$ correspond to $\gG$ in $\dbco$. Let 
$I^p = \oplus_{q \in \hele} \lv(-q) \te V^p_{q}$. Then
\[ \chi  \gG(n) = \sum_{p \in \hele} (-1)^{p+n} \dim_k V^p_{n} \]
(if this sum if finite).
\end{theorem}

\begin{proof}
If $\gF$ is a coherent sheaf and $I = E(\gF)$ then the equality
clearly holds. Since $\dbco$ is generated as a triangulated category
by coherent sheaves, and the right hand side in the equality above is 
invariant for homotopic complexes, the theorem follows.

\end{proof}

\section { How to determine properties of a coherent sheaf from its Tate 
resolution. } \label{6egenskap}

In this section we shall further investigate how properties of an object 
$\gG$ in $\dbco$ is reflected in the exterior complex $\ytrev (\gG)$. For 
instance, if $\gG = \gF$ is a coherent sheaf ,  then it is often important to 
determine the ranks of the fibers $\gF_{k(P)}$ for $P$ a point in $\pv$ and 
also the dimensions of the loci where $\gF$ has constant rank.
Furthermore it is also often important to determine the projective 
dimensions of the localizations $\gF_P$ for points $P$ in $\pv$. 

We shall show how these invariants may be 
computed from the exterior complex 
$\ytrev ( \gF)$. 
These derivations come from investigations of the following.
We consider a linear subspace $\pvw \sus \pv$. Let 
$\gG$ be in $\dbco$. Then we may obtain a suitable restriction $\lgvw $ in 
$\dbcovw$ (in a derived sense). We shall describe how the 
exterior complexes of
$\lgvw$ and $\gG$ are related. In particular we shall use how these
are related when $\pvw$ is a point $P$.

\subsection { Restriction to linear subspaces. }

Recall from Proposition \ref{2Pderekv} 
that there is an equivalence of categories
\begin{equation}  \dbvb \bihom{ i_{\pv}}{j_{\pv}}  \dbco. \label{6vbco} 
\end{equation}
Let $U \sus W$ be a subspace, so $\pvw \sus \pv$ is a linear subspace. Then 
we have a restriction functor
\[ \dbvb \pil \dbvbvw \]
by letting 
\[ (\cdots \pil \gE^i \pil \gE^{i+1} \pil \cdots) \,   \mapsto  \,  
(\cdots \pil \gE^i_{|\pvw} \pil \gE^{i+1}_{|\pvw} \pil \cdots) . \]
Via the equivalence (\ref{6vbco}) and the corresponding equivalence for 
$\pvw$, we then obtain a restriction functor
\begin{eqnarray*}  \dbco & \pil & \dbcovw \\
\gG & \mapsto & \lgvw. \\
\end{eqnarray*}
Note that if $\gG$ is a coherent sheaf $\gF$, then $\lfvw$ may not be 
isomorphic in $\dbcovw$ to a coherent sheaf. In fact we have the following.

\begin{proposition} \label{6FtilS}
Let $\gF$ be a coherent sheaf and let $S \sus \pv$ be the locus where $\gF$ 
degenerates in rank. Then $\lfvw$ is isomorphic in $\dbcovw$
to a coherent sheaf (necessarily $\gF_{|\pvw}$) if and 
only if $\pvw$ intersects $S$ properly. 
\end{proposition}

\begin{proof} 
Suppose first that $U = (u)$ so $\pvw$ is a hyperplane in $\pv$. 

Let $\gE$ be a locally free resolution of $\gF$, and let $\gK^i = \ker ( 
\gE^i \pil \gE^{i+1})$. Then there are exact sequences 
\[ 0 \pil \gK^i \pil \gE^i \pil \gK^{i-1} \pil 0. \]
Since $\gK^{i-1}$ is a locally free sheaf on 
$\pv \backslash  S$ and no component of $S$ (in reduced form)
is contained in $\pvw$, 
the element $u$ is not contained in any 
associated prime of $\gK^{i-1}$. In the diagram
\[ \begin{CD}
 \gK^i(-1) @>>> \gE^i(-1) @>>> \gK^{i-1}(-1) \\
@VV{\cdot u}V @VV{\cdot u}V @VV{\cdot u}V \\
\gK^i @>>> \gE^i @>>> \gK^{i-1} 
\end{CD} \] 
the right vertical map is therefore injective. Thus the cokernels
\[ 0 \pil \gK^i / u\gK^i \pil \gE^i / u \gE^i 
\pil \gK^{i-1}/ u \gK^{i-1} \pil 0 \]
is an exact sequence. This gives that $\gE/ u\gE$ is a locally free resolution 
of $\gF / u \gF$ and hence $\lfvw$ is isomorphic to a coherent sheaf. 

By cutting down with hyperplanes, this proves the if part of the statement
in generality.

\medskip

To prove the converse we may assume by first cutting down with a linear
subspace intersecting $S$ properly, that a component of $S$ (in reduced form)
is contained in $\pvw$. Let $u$ be in $U$. Then 
\[ \gK^0 /u\gK^0 \pil \gE^0 /u \gE^0 
\pil \gF/u \gF \pil 0 \]
is not left exact. Applying this again to the exact sequence
\[ 0 \pil \gK_1 \pil \gE^0 /u \gE^0  \pil \gF/u \gF \pil 0 \]
where $\gK_1$ is the kernel, and proceeding, we see that the sequence
\[ \gK^0 /U\gK^0 \pil \gE^0 / U\gE^0 
\pil \gF/ U\gF \pil 0 \]
is not left exact. Hence $\gE_{| \pvw}$ 
is not a resolution of $\gF_{| \pvw}$.
\end{proof}

\nota
Recall from Notation \ref{3ytrev} 
that for $\gG$ in $\dbco$ we get the exterior complex
$\ytrev(\gG)$ which is the composition of $\gsmin$ and $\trgstar(\pv,\gG)$. 
If $\gG^\prime$ is an object
in $\dbcovw$ we shall in the following denote the composition of
$G_{\svw, min}$ and $\trgstar (\pvw, \gG^\prime)$ as $\ytrevw (\gG^\prime)$.
\notafin

The following says how restrictions of complexes of vector bundles
from $\pv$ to $\pvw$ translate when looking at the corresponding
exterior complexes in $K(\lvd)$ and $K(\lvwd)$.

\begin{proposition} The diagram
\[ \begin{CD}
\dbvb @>{\ytrev   }>> K(\lvd) \\
@VV{{}_{|\pvw}}V @VV{\res^{\lvd}_{\lvwd}}V \\
\dbvbvw @>{\ytrevw }>> K(\lvwd) 
\end{CD} \]
gives a natural isomorphism of functors
\[ \res^{\lvd}_{\lvwd} \circ \ytrev  \ispil \ytrevw  \circ _{|\pvw} . \]
\end{proposition}

\begin{proof}
We may assume $\pvw \sus \pv$ is a hyperplane inclusion so $U = (u)$ is a 
one-dimensional subspace of $W$. Let $\gE$ be in $\dbvb$ and let $J$ be in 
$K(\lvdcf)$ be such that $\tilfl(J)$ has bounded cohomology. (Note that 
this is always true if $J = \ytrev (\gG)$ for some
$\gG$ in $\dbco$.) 

Then to give a morphism $J \pil \gs \circ \rgstar(\pv,\gE)$ corresponds to a 
morphism 
$\fl (J) \pil \rgstar (\pv,\gE)$ which gives a morphism $\tilfl (J) \pil 
\gE$ in $\dcq$.
Such a map can be represented by a diagram of morphisms of complexes
\[  \begin{array}{ccccc}
 &  & \gQ & &  \\
 &  \overset {\phi} \swarrow & & \searrow & \\
\tilfl(J) & & & & \gE 
\end{array} \]
where $\phi$ is a quasi-isomorphism. By the proof of 
Proposition \ref{2Pderekv} we may 
assume that $\gQ$ is a bounded complex of vector bundles. But then we clearly 
get a diagram
\[ \begin{array}{ccccc}
& & \gQ_{|\pvw} & & \\
& \overset{\phi_{|\pvw}} \swarrow & & \searrow & \\
\tilfl (J) _{|\pvw} & & & & \gE_{|\pvw} 
\end{array} \]
where $\phi_{|\pvw}$ is a quasi-isomorphism.

This gives us a morphism 
\[ \tilfl (J) / u \tilfl(J) \pil \gE_{|\pvw} \] 
in $\dbcovw$ and thus a morphism in $K(\svw)$
\begin{equation} \fl(J)/ u \fl(J) \pil \rgstar (\pvw,\gE_{|\pvw}). 
\label {6fjvw} \end{equation}

By diagram (\ref{1resVW}) in Subsection \ref{1subsp} 
there is an isomorphism of functors
\[ (\svw \te_{\sv} -) \circ \fl \iso \flvw  \circ \res^{\lvd}_{\lvwd} . \]
Thus we get an isomorphism 
\[ \fl(J) / u \fl(J) \iso \flvw \circ 
\res^{\lvd}_{\lvwd} (J).\]

From this isomorphism and  (\ref{6fjvw}) we get, using the adjointness of 
$\fl$ and $\gs$,
a morphism in $K(\lvwd)$
\begin{equation}  \res^{\lvd}_{\lvwd} (J) \mto{\alpha} \gsvw \circ \rgstar 
(\pvw, \gE_{|\pvw}).  \label{6resg} \end{equation}
Taking $J = \ytrev (\gE)$ this gives our natural transformation of functors.

It remains to prove that the natural transformation is an isomorphism. 
There is a homotopy equivalence $\phi$
\[ \gsvw \circ \rgstar (\pvw, \gE_{|\pvw}) \isopil 
      \gswmin \circ \trgstar(\pvw,\gE_{|\pvw}). \]
By the constructions above $\tilflvw (\alpha)$ is an isomorphism in
$\dcq$ and hence $\tilflvw (\phi \circ \alpha)$ is an isomorphism in 
$\dcq$. By Proposition \ref{3Inullht}, $\phi \circ \alpha$ is a 
homotopy equivalence, and hence $\alpha$ is a homotopy equivalence.
\end{proof}



The following now shows how the restriction functor $\lder(-)_{|\pvw}$
translates to a functor from $\kolvdcf$ to $\kolvwdcf$.

\begin{corollary} \label{6LEres}
The diagram
\[ \begin{CD}
\dbco @>{\ytrev}>> \kolvdcf \\
@V{\lder (-)_{|\pvw}}VV @VV{\res^{\lvd}_{\lvwd}}V \\
\dbcovw @>{\ytrevw}>>  \kolvwdcf 
\end{CD} \]
gives a natural isomorphism of functors
\[ \res^{\lvd}_{\lvwd} \circ E \isopil E^\prime \circ \lder (-)_{|\pvw}. \]
\end{corollary}

\begin{proof} Immediate. \end{proof}

The complex $\ytrev(\gG)$ is acyclic. We also get a complex 
\[ \Hom_{\lvd}(\lwd,  \ytrev(\gG))\] 
which is in general not acyclic when $U \neq W$.
The  following tells us what information the cohomology groups of this complex
gives.

\begin{proposition} \label{6kohytre}
Let $\gG$ be in $\dbco$.  Then there is a natural 
isomorphism of $\svw$-modules 
\begin{eqnarray*} & \oplus_{p,q \in \hele} & H^p(\rgstar(\pvw,\lgvw))_q \\
\iso & \oplus_{p,q \in \hele} & H^{p+q} \Hom_{\lvd} (\lwd, \ytrev (\gG))_{-q} 
\end{eqnarray*}
where the $\svw$-module structure on the latter bigraded group are given as 
in Subsection \ref{1sekWFG}.
\end{proposition} 


\begin{proof}
By Proposition \ref{1PIht} there is a "twisted" quasi-isomorphism
\[ \Hom_{\lvd} (\lwd, \ytrev (\gG)) \iso \svw \te_{\sv} \fl(\ytrev (\gG)) \]
and by diagram (\ref{1resVW}) in Subsection \ref{1subsp}
there is an isomorphism of complexes of $\svw$-modules
\[ \svw \te_{\sv} \fl (\ytrev (\gG)) \iso 
\flvw \circ \res^{\lvd}_{\lvwd}(\ytrev (\gG)). \]
By the previous Corollary \ref{6LEres} there is a homotopy equivalence of 
complexes of $\lvwd$-modules
\[ \res^{\lvd}_{\lvwd} \circ \ytrev (\gG) \iso \ytrevw ( \lder \gG_{|\pvw}). \]
Now for any $N$ in $D(\svw)$ there is a quasi-isomorphism of complexes
of $\svw$-modules
\[ \flvw \circ  \gsvw (N) \pil N. \]

Then by letting $N$ be $\rgstar (\pvw,\lder \gG_{|\pvw})$ 
we get a quasi-isomorphism
\[ \flvw \circ \ytrevw (\lder \gG_{|\pvw}) \pil 
\rgstar (\pvw, \lder \gG_{| \pvw}) \]
and this proves the proposition.
\end{proof}

The meaning of the previous proposition is most transparent in the
following case. 

\begin{corollary} Let $\gF$ be a coherent sheaf and suppose 
$\lfvw $ is isomorphic to $\gF_{|\pvw}$. Then
\[ H^p \gF_{|\pvw} (q) \iso  H^{p+q} \Hom_{\lvd}(\lwd, \ytrev (\gF))_{-q} . \]
\end{corollary}

\remark In conjunction with the theory in Subsection \ref{1sekWFG} we see that 
a necessary condition for a coherent sheaf $\gF$ on $\pvw$ to lift to $\pv$ 
is that there is a structure of $\lwd$-module on $\oplus_{p,q \in \hele} H^p 
\gF(q)$ such that the actions of $\svw$ and $\lwd$ commute.
\remarkfin

\subsection{How properties of a coherent sheaf translate to the
Tate resolution.} \label{6lokalSek}
Now we turn to the following question. Given a coherent sheaf $\gF$. How does 
one determine its rank, the dimensions of its degeneracy loci, and its local 
projective dimension, from $\ytrev (\gF)$? The key is the following.

\begin{theorem} \label{6key}
Let $U \sus W$ be of codimension $1$ so $P = \pvw \inpil \pv$ is a point. Let 
$\gF$ be a coherent sheaf  on $\pv$ and $\gF_P$ its localization in $P$. 
Then for any $p$ we have a functorial isomorphism of vector spaces
\[ \Tor_{-p}^{\gO_{\pv,P}}(k(P), \gF_P) \iso 
H^{p+q} \Hom_{\lvd} (\lwd, \ytrev (\gF))_{-q} \]
for all $q$. In particular for a fixed $a = p+q$ one has for every $p$
\[ \Tor_{-p}^{\gO_{\pv,P}}(k(P), \gF_P) \iso  
 H^a \Hom_{\lvd}(\lwd, \ytrev (\gF))_{p-a}. \]
\end{theorem}

\begin{proof} First note that $\gF_{k(P)}$ and $\gF_{|P}$ both 
denote the fiber of $\gF$ at the point $P$. 
Let
\[ 0 \pil \gE^{-r} \pil \gE^{-r+1} \pil \cdots \pil \gE^0 \pil \gF  \]
be a locally free resolution. Then $\lder \gF_{|P}$ may be identified with the
complex $\gE_{|P} = \gE_{k(P)}$. Since this is sheaves over a point $P$,
$H^p (\rgstar (P,\lder \gF_{|P}))_q$ is isomorphic to 
$H^p (\rgstar (P,\lder \gF_{|P}))_0$ for any $q$ and the latter may be 
identified with 
\[ H^p(\gE_{k(P)}) = \Tor_{-p}^{\gO_{\pv,P}}(k(P), \gF_P). \] 

Since
\[ H^p (\rgstar (P,\lder \gF_{|P}))_q \iso H^{p+q} 
\Hom_{\lvd}( \lwd, \ytrev(\gF))_{-q} \]
by Proposition \ref{6kohytre}, we get the theorem.
\end{proof}

The following corollary shows how we may determine the rank of a coherent 
sheaf $\gF$ at a point $P$ and the projective 
dimension of the localization $\gF_P$, 
from the exterior complex $\ytrev (\gF)$. In fact we are able to 
determine all this by a very local consideration on the exterior complex. 
Given 
any integer $a$ we determine these data just from looking at the terms
\[ \ytrev (\gF)^{a-1} \pil \ytrev (\gF)^a \pil \ytrev (\gF)^{a+1}. \]

\begin{corollary} \label{6lokalP} Fix any integer a. Then the folowing holds.

a. $H^a \Hom_{\lvd} (\lwd, \ytrev (\gF))_{p-a} = 0$ for all $p \geq 1$.

b. The rank of the fiber $\gF_{k(P)}$ is the vector space dimension of
\linebreak
$H^a \Hom_{\lvd} (\lwd, \ytrev (\gF))_{-a} .$

c. The projective dimension of the localization $\gF_P$ is the largest
integer $l$ such that $H^a \Hom_{\lvd} (\lwd, \ytrev 
(\gF))_{-l-a}$ is non-zero.
\end{corollary} 

\begin{proof} 
The projective dimension of $\gF_P$ is the largest integer $l$ such that
\[ \Tor_l^{\gO_{\pv,P}}(k(p), \gF_P) \]
is non-zero, by \cite[Cor. 19.5]{Ei}. Hence c. follows, 
and a. and b. are clear. 
\end{proof}

\begin{corollary} \label{6bunt}  Given an integer $a$. Then
$\gF$ is a vector bundle of rank $r$ if and only if
\[ \dim_k H^a \Hom_{\lvd} (\lwd, \ytrev (\gF))_{-a} = r \]
for all $U \sus W$ of codimension one.
\end{corollary}

\begin{proof} Immediate. \end{proof}

The following shows how to determine the dimension of the degeneracy
locus of the coherent sheaf $\gF$. 

\begin{theorem}
Let $U_0 \sus W$ be a subspace of codimension $d$, so $\pvwo \sus \pv$ is a 
linear subspace of dimension $d-1$. If for all $U_0 \sus U \sus W$ 
where $U$ is of codimension one in $W$
we have 
\begin{equation} \label{6dimWF}
 \dim_k H^a \Hom_{\lvd}(\lwd, \ytrev (\gF))_{p-a} =    
\begin{cases}  r,  & \text{ for } p = 0 \\
              0,  & \text{ for } p \neq 0
\end{cases} \end{equation}
then the degeneracy locus of $\gF$ is of codimension greater or equal to $d$ 
and $\rank(\gF) = r$. Conversely if its degeneracy locus is of codimension 
greater or equal to $d$ and $\rank(\gF) = r$ then (\ref{6dimWF}) holds if 
$W_0$ is a general subspace of $V$ of codimension $d$.
\end{theorem}

\begin{proof}
Let $S$ be the degeneracy locus of $\gF$. 
By Proposition \ref{6FtilS},
$P $ is not in $S$ for any point $P$ in $\pvwo$. Hence $S \cap 
\pvwo$ is empty, and $\codim S \geq d$ and $\rank(\gF) = r$. The converse 
part is also clear.
\end{proof}

\subsection{How properties of a Tate resolution translate to a
coherent sheaf.} \label{6EtilFSek}
Now often one would probably not start from a coherent sheaf $\gF$ on $\pv$. 
Rather one would start with a homomorphism $d : E^\prime \pil E$ of cofree 
left $\lvd$-modules of finite corank. By Remark \ref{3dtilG} 
this gives  rise to an 
object $\gG$ in $\dbco$.  So how can one determine if $\gG$ is actually 
(isomorphic to) a coherent sheaf ? And if so, which properties does $\gG$ 
have ? We consider this in Theorem \ref{6dtilF}, 
but first we do some preparatory 
work.

\begin{lemma} \label{6lpdl} 
Let $\gF$ be a coherent sheaf on $\pv$ and let $l$ be an integer.
Then there is an open subset $\op$ with $\codim (\pv \backslash \op) \geq l+1$
such that the projective dimensions of the localizations $\gF_P$ 
is $\leq l$ for all $P $ in $\op$.
\end{lemma}

\begin{proof}
This follows by Theorem 20.9 and Proposition 18.2 in \cite{Ei}.
\end{proof}

\begin{lemma} \label{6lemkohvan}
Let
\[ \gE^{-l} \mto{d^{-l}} \gE^{-l+1} \pil \cdots \pil \gE^0 \]
be a complex of vector bundles on $\pv$. Suppose $H^p(\gE)$ has support in 
codimension $\geq l$ when $p < 0$. Then $H^p(\gE) = 0$ for $p < 0$.
\end{lemma}

\begin{proof}
Note that if $\gF$ is a coherent sheaf on $\pv$, then 
the local projective dimension
\[ \lpd(\gF) = \max \{ i 
\, | \, 
\gext ^i (\gF, \gR) \neq 0 \text{ for some } \gR \text{ in } \coh \}. \]

Now if $H^p (\gE) \neq 0$ for some $p < 0$, we may assume that the 
codimension of the support of $H^p(\gE)$ is equal to $l$ for some $p$. (Else 
we may replace $l$ with a larger integer.) 


By Lemma \ref{6lpdl} there is an open subset $\op$ of $\pv$ 
with the codimension of $\pv \backslash \op $ greater or equal to $l + 1$ 
such that the 
projective dimension of all localizations $H^p (\gE_P)$ is less or equal to
$l$ for all integers $p$ and points $P$ in $\op$.

Let $\gK^i = \ker d^i$ and $\gB^i = \im d^{i-1}$.
Also let $p < 0$ be minimal such that $H^p(\gE_{|\op}) \neq 0$. 
(Note that such a $p$ exists by our assumptions on $\op$ and $H^p(\gE)$ for 
$p < 0$.)
We have an exact sequence 
\begin{equation}  0 \pil \gB^p_{|\op} \pil \gK^p_{|\op} \pil H^p(\gE_{|\op}) \pil  
0. \label{6bkh}
\end{equation}
Since $H^i(\gE_{|\op}) = 0$ for $i < p$ we get $\lpd(\gB^p_{|\op}) \leq l-2$. 
By the Auslander-Buchsbaum Theorem \cite[Thm. 19.9]{Ei}, 
the local projective dimension
of $H^p(\gE_{|\op})$ is $l$. 
Hence $\lpd(\gK^p_{|\op})$ is also $l$. 

By the exact sequence
\[ 0 \pil \gK^{p}_{|\op} \pil \gE^{p}_{|\op} \pil \gB^{p+1}_{|\op} \pil 0 \]
we get $\lpd(\gB^{p+1}_{|\op}) = l+1$. 
By the exact sequence (\ref{6bkh}) for $p+1$ 
we get that either $\lpd(H^{p+1}(\gE_{|\op})) \geq l+2$ or $\lpd(\gK^{p+1}_{|\op}) 
\geq l+1$. In the latter case we may continue and eventually get 
$\lpd(H^{p+r}(\gE_{|\op})) \geq l+r+1$ for some $r \geq 1$.  
This is not possible 
however since $\lpd(H^{p+r}(\gE_{|\op})) \leq l$.

Thus we have proved the lemma.
\end{proof}

\begin{lemma} \label{6lemdeg}
Let $\gG$ be in $\dbco$ and let $a$ be an integer. Let $U_0 \sus 
W$ have codimension $d$. Suppose that for all $ U_0 \sus U \sus W$ where $U$ 
has codimension $1$, the following holds for some integer $p_0$
\begin{equation} \dim_k H^a \Hom_{\lvd} (\lwd, \ytrev (\gG))_{p-a} = 
\begin{cases} r, &\text{ if } p = p_0. \\
 0, & \text{ otherwise. }
\end{cases} \label{6r} 
\end{equation}
Then $H^p(\gG)$ is a sheaf of rank $r$ for $p = p_0$ and rank $0$ for $p \neq 
p_0$. 
Furthermore if $H^p(\gG)$ degenerates in a locus $S^p$, then $S^p \cap \pvwo$ 
is empty, and so $\codim S^p \geq d$.
\end{lemma}

\begin{proof} The condition (\ref{6r}) says, as in Theorem \ref{6key} b., that 
when $p \neq p_0$ then
$H^p(\lder 
\gG_{|P}) = 0$ for all points $P = \pvw \inpil \pvwo$. 
Thus when $p \neq p_0$, clearly $S^p \cap \pvwo$ is empty and so 
$\codim S^p  \geq d$.  For $P = \pvw$ a point in 
$\op = \pv \backslash \cup_{p \neq p_0} S_p$ the rank 
of $H^p (\lder \gG_{|\pvw})$ must then be constant for all $p$. 
Since $\pvwo \sus \op$ 
we get that $H^{p_0} (\gG)$ is a sheaf of rank $r$.
\end{proof}

\pagebreak

\begin{theorem} \label{6dtilF}
 Let 
\[ E ^{\cdot} \, : \, E^\prime \pil E \pil E^{\prime \prime} \]
be a complex of cofree $\lvd$-modules of finite corank which is exact in the 
middle. In the following let $U \sus W$ denote a subspace of codimension 
one.  Suppose the following holds for some integers $a$ and $l$.
\begin{enumerate}[a.]
\item $ H \Hom_{\lvd}(\lwd, E^\cdot)_q$ is  $ \begin{cases} = 0 &\text{ 
for } q \text{ not in } [-a-l,-a] \text{ and all } U \sus W. \\
\neq 0 & \text{ for some } U \sus W \text{ when } q = -a
\end{cases}$
\item Let $U_0 \sus W$ be some subspace of codimension $l$. Suppose the
vector space dimension of 

$ H \Hom_{\lvd} (\lwd, E^\cdot)_q$ is $  
\begin{cases}   r & \text{for all } U_0 
\sus U \sus W \text{ when } q = -a. \\
0 & \text{for all } U_0 \sus U \sus W \text{ when } q \neq -a. 
\end{cases} $
\end{enumerate}
By taking a free resolution of the kernel of $E^\prime \pil E$ and a 
cofree resolution of the cokernel of  $E \pil E^{\prime \prime}$ 
we may complete 
$E^\cdot$ to an acyclic complex $I$ of cofree modules of finite corank. Index 
this complex so $E = I^a$ is in component $a$. Then this complex $I$ 
corresponds to a coherent sheaf $\gF$ on $\pv$ of rank $r$ with 
$\lpd(\gF) \leq l$. Furthermore $\gF$ degenerates in codimension $\geq l$. 
\end{theorem}

\begin{proof} Let $I$ correspond to $\gE$ in $\dbvb$. 
The condition {\it a.} gives that the ranks of the maps
\[ \gE^{p-1}_{|\pvw} \pil \gE^p_{|\pvw} \]
are constant for all $p \leq -l$ and all $p \geq 1$. Thus the kernels and 
cokernels of these maps are bundles. We may therefore assume we have a complex
\[ \gE^{-l} \pil \gE^{-l+1} \pil \cdots \pil \gE^0. \]

Now we use condition {\it b.} We may apply Lemmata \ref{6lemkohvan} and 
\ref{6lemdeg} and get that 
$H^p(\gE) = 0$ for $p \neq 0$. Thus letting $\gF = H^0(\gE)$ 
we have the theorem.
\end{proof}

\begin{corollary} If also the cohomology 
$H \Hom_{\lvd}(\lwd, E^\cdot)_q$ is zero for $q \leq 
-a-l$ and all $U \sus W$, then $\gF$ is a torsion free sheaf of local 
projective dimension less than $l$.
\end{corollary}

\begin{proof} We may assume that $\gE$ in the proof above is
\[ \gE^{-l+1} \pil \cdots \pil \gE^0 . \]
Thus $\lpd(\gF) \leq l-1$. Furthermore there is an exact sequence
\begin{equation} 0 \pil \gT \pil \gF \pil \overline{\gF} \pil 0 \label{6tf} 
\end{equation} 
where $\overline{\gF}$ is torsion free and $\gT$ is the torsion part of $\gF$.
Let $k$ be the codimension of the support $\text{Supp}\, \gT$. 
Since $\gF$ degenerates in codimension $\geq l$ by Theorem
\ref{6dtilF} we must have $k \geq l$. 
Also the projective
dimension of $\gT_P$ is  $\geq k$ for
$P $ in $\text{Supp}\, \gT$ by the Auslander-Buchsbaum theorem 
\cite[Thm. 19.9]{Ei}.
Then the projective dimension of $\overline{\gF}_P \geq k +1 $ for $P$ in
$\text{Supp}\, \gT$. But this is impossible since the locus of $P $ in 
$\pv$ such that the projective dimension of $\overline{\gF}_P$ is 
 $\geq k +1$ has codimension $ \geq k + 1$ according to Lemma \ref{6lpdl}.
 Hence $\text{Supp}\, \gT = \emptyset$ and $\gF = \overline{\gF}$.

\end{proof}

\section{ Projections}

Let $ U \sus W$ be a linear subspace and let $\op$ be the open
subset $\pv - \pvw$. Se we get a projection morphism $\op \lpil \pw$. 
If $\gF$ is a coherent sheaf on $\pv$ such that $\supp \gF$ is 
disjoint from $\pvw$, then $p_* \gF$ is a coherent sheaf on $\pv$ and
we show that its cohomology is related to that of $\gF$ by
\[ H^i(\pv, \gF) = H^i (\pw, p_* \gF), \]
or, if we consider them as $\sv$- and $\sw$-modules, then
\[ \res^{\sv}_{\sw} H^i_* (\pv, \gF) = H^i_* (\pw, p_* \gF). \]
It is then quite immediate that, if $\ytrev(\gF)$ is the Tate resolution of
$\gF$, then the Tate resolution of $p_* \gF$ is 
\[ \ytrev(p_* \gF) = \Hom_{\lWd}(\lwd, \ytrev(\gF)). \]
We shall give the full derived category version of this.

\subsection{Projections and derived categories.}

Let 
\[ \dbcoM,\, \dbqcM, \,  \text{ and } \dbqcoM \]
be the full subcategories of
\[ \dbco,\, \dbqc,\,  \text{ and } \dbqco \]
consisting of $\gG$ such that $H^i(\gG)$ is a coherent sheaf on $\pv$
and $\supp H^i(\gG)$ is disjoint from $\pvw$. Also let 
$\kolvdcfMw$ be the full subcategory of $\kolvdcf$ consisting of $I$
such that $\Hom_{\lvd}(\lwd, I)$ is acyclic.
\medskip
From Proposition \ref{6kohytre} it is clear that via the equivalence
\[ \dbco \mto{T} \kolvdcf, \]
the subcategory $\dbcoM$ maps to $\kolvdcfMw$ giving an equivalence
of subcategories.
\medskip
Now there are functors
\begin{eqnarray*} && \dbcoM  \mto{i} \dbqcM \\
& \mto{|_\op} & \dbqcoM  \mto{Rp_*} \dbqcw.
\end{eqnarray*}

\begin{proposition} If $\gF$ is a coherent sheaf in $\pv$ with $\supp \gF$ 
disjoint from $\pvw$, then 
\[ Rp_* \circ |_\op \circ i (\gF) \iso p_* \gF \]
in $\dbqcw$. As a consequence,
the image of $Rp_* \circ |_\op \circ i$ is in the 
full subcategory $\dbcqw$. \label{7cohproj} \end{proposition}

\begin{proof} The latter follows from the former since 
$\dbcoM$ is generated by such coherent sheaves $\gF$ and for these
$p_* \gF$ is coherent. 

Consider the functors
\[ \qc \overset{|_\op}{\underset{i_*}{\rightleftarrows}} \qco. \]
By \cite[III.Ex.3.6]{Ha} if $\gI$ is injective in $\qc$, then $\gI_{|_\op}$
is injective. Also since $i_*$ is right adjoint to $|_\op$ and the latter
is exact, if $\gI$ is injective in $\qco$ then $i_* \gI$ is injective.

Let $\gF \pil \gI$ be an injective resolution in $\qc$. 
Then $\gF \pil \gI_{|\op}$ is an injective
resolution in $\qco$ and $\gF \pil i_*(\gI_{|\op})$ is an injective
resolution in $\qc$.  

We now claim that $p_* \gF \pil p_* (\gI_{|\op})$
is a resolution in $\qcw$.
To see this, note that
$\Gamma(\pv, \gF(n)) \pil \Gamma(\pv, i_*(\gI_{|\op})(n))$ is a resolution for 
$n \gg 0$. Hence it follows that $p_* \gF \pil p_* \gI_{|\op}$ must be a 
resolution (see \cite[II.Ex.5.15 e)]{Ha}). Since $p_* \gF$ is coherent 
on $\pw$, we are done.
\end{proof}

Now let $\dbcqw \mto{j} \dbcow$ be a functor quasi-inverse to the natural
inclusion.

\begin{theorem} The diagram of functors
\[ \begin{CD}  \dbcoM  @>{\ytrev}>>  \kolvdcfMw \\
     @V{j \circ Rp_* \circ {|_\op} \circ i}VV  @VV{\Hom_{\lvd}(\lwd, -)}V \\
\dbcow  @>{\ytrevw}>>  \kolwdcf
\end{CD} \] 
gives a natural isomorphism of functors 
\[ \Hom_{\lvd}(\lwd, -) \circ \ytrev  \iso \ytrevw \circ 
(j \circ Rp_* \circ {|_\op} \circ i). \]
\end{theorem}

\begin{proof}

   Consider the diagram

\[ \begin{CD}\dbqcM @>{R\Gamma_*}>>  D(\sv) \\
 @VV{|_\op}V   \\
\dbqcoM  @.  @VV{\text{res}}V \\ 
@VV{Rp_*}V \\
\dbcqw  @>>{R\Gamma^\prime_*}>  D(\sw) 
\end{CD} \]

We first claim that there is a natural isomorphism of functors
\begin{equation} \text{res} \circ  R\Gamma_* 
\pil R\Gamma^\prime_* \circ Rp_* \circ |_{\op}. \label{7natiso}
\end{equation}

Let $\gG \pil \gI$ be an injective resolution in $\qc$.
By the first part in the proof of Proposition \ref{7cohproj},
$\gG_{|\op} \pil \gI_{|\op}$ is an injective resolution in
$\qco$ and $i_* (\gG_{|\op}) \pil i_* (\gI_{|\op})$ is an injective
resolution in $\qcw$. Since the natural map $\gG \pil i_* (\gG_{|\op})$
is a quasi-isomorphism, $\gI \pil i_* (\gI_{|\op})$ is a homotopy 
equivalence and so $\phi : \Gamma_*(\gI) \pil \Gamma_*(\gI_{|\op})$ is
an isomorphism in $D(\sv)$.

Now $Rp_* (\gG_{|\op})$ is $p_* (\gI_{|\op})$ and since the latter is 
a complex of flasque sheaves, there is a quasi-isomorphism
\[ \Gamma^\prime_* \circ p_* (\gI_{|\op}) 
\mto{\psi} R\Gamma^\prime_* \circ Rp_* (\gG_{|\op}). \] 
Composing $\text{res}\, (\phi)$ with $\psi$ shows the claim (\ref{7natiso}).

The theorem now follows from the diagram
\[ \begin{CD} D(\sv) @<{\til}<< K(\sv)  @>{\gs}>> \klvdcf \\
@V{\text{res}}VV @V{\text{res}}VV  @VV{\Hom_{\lvd}(\lwd,-)}V \\
D(\sw) @<{\til}<< K(\sw) @>{\gsw}>> \kolwdcf
\end{CD} \]
where the second square gives a natural isomorphism of functors
\[ \gsw \circ \text{res} \iso \Hom_{\lvd}(\lwd,-) \circ \gs. \]
\end{proof}

\section{ $G$-equivariant quasi-coherent sheaves}

Let $G$ be a linear algebraic group acting on $W$, so $\pv$ comes with a 
$G$-action. 
This section contains generalities on $G$-equivariant quasi-coherent
sheaves on $\pv$ and the correspondence with $\sv,G$-modules. 
The things here are
standard but since it is hard to find a precise reference for the
relation between $G$-equivariant quasi-coherent sheaves on $\pv$
and $\sv,G$-modules, we present the arguments in some detail.
In particular, during the arguments we shall extensively use the language
of Hopf-algebras and comodules (see \cite{Ka}) 
when considering $\sv,G$-modules.

Consider the diagram
\[  \begin{array}{ccc}
G \times \pv & \mto{\eta}  & G \times \pv \\
        & q \searrow & \downarrow p \\
        &  & \pv 
\end{array} \]
where $p$ is the projection, $q$ is the action map $q(g,p) = g.p$
and $\eta(g,p) = (g, g.p)$. Following \cite{MuFo}, a $G$-equivariant
quasi-coherent sheaf $\gF$ on $\pv$ is a quasi-coherent sheaf $\gF$ on 
$\pv$ together with an isomorphism
\[ \phi : q^*(\gF) \pil p^*(\gF). \]
This isomorphism must satisfy the following identity of morphisms on
$G \times G \times \pv$. 
\begin{equation}
 (p_{23}^* \phi) \circ (1 \times q)^* \phi = (m \times 1)^* \phi 
\label{7kosykel}
\end{equation}
where $m : G \times G \pil G$ is the multiplication map and $p_{23}$
is the projection on the second and third factor. We let $\gqc$ be the
{\it category of $G$-equivariant quasi-coherent sheaves on $\pv$}.

\medskip

\subsection{ Notation and terminology.}

Let $\alpha : A \pil B$ be a ring homomorphism and $M$ an $A$-module.
To make explicit that in $B \te_A M$ the $A$-module structure on $B$
comes from $\alpha$ we denote this by $B {}^{\alpha} \te_A M$. 


Now the coordinate ring $k[G]$ is a Hopf algebra. 
Denote by 
\[ \Delta_{k[G]} \, : \, k[G] \lpil k[G] \te k[G] \]
the coalgebra map coming from the multiplication 
$G \times G \pil G$ and by 
\[ m_{k[G]} \, : \, k[G] \te k[G] \lpil k[G] \]
the algebra map (coming from the diagonal map).

If $M$ is an $\svgmod$, 
recall from Subsection \ref{1sekGr} that this is equivalent to $M$ being a 
{\it i.} $k[G]$-comodule {\it ii.} $\sv$-module and {\it iii.} the module map
$\sv \te M \pil M$ is a $k[G]$-comodule map. 

We denote by 
\[ \Delta_M : M \lpil k[G] \te M \]
the $k[G]$-comodule map. Given {\it i.} and {\it ii.}, then 
{\it iii.} is easily seen to 
be equivalent to the assumption that 
$\Delta_M$ is an $\sv$-module map where the $\sv$-module structure
on $k[G] \te M$ is given via the $k[G]$-comodule map 
\[ \Delta_{\sv} : \sv \pil k[G] \te \sv. \]

For the map
\[ \Delta_{k[G]} : k[G] \pil k[G] \te k[G] \]
we use Sweedler's sigma notation (see \cite{Ka}.) 
\[ \Delta_{k[G]}(\mu) = \sum_{(\mu)} \mu^\prime \te \mu^{\prime\prime} \]
and if $m$ is in $M$ then
\[ \Delta_M(m) = \sum_{(m)} m_{k[G]} \te m_M. \]

Also denote the image of $m$ by the commutative diagram
\[  \begin{CD}
M @>{\Delta_M}>>  k[G] \te M \\
@V{\Delta_M}VV @VV{k[G] \te \Delta_M}V  \\
k[G] \te M  @>{\Delta_{k[G]} \te \id_M}>> k[G] \te k[G] \te M
\end{CD} \]
as 
\[ 
\sum_{(m)} m^\prime_{k[G]} \te m^{\prime \prime}_{k[G]} \te m_{M} \]
which is equal to the two expressions
\[ \sum_{(m)} ( \sum_{(m_{k[G]})} m^\prime_{k[G]} \te 
 m^{\prime \prime}_{k[G]}) \te m_{M}, \quad
 \sum_{(m)} m_{k[G]} \te (\sum_{(m_M)} (m_M)_{k[G]} \te (m_M)_M) . \]

 \subsection{Correspondence between equivariant sheaves and modules}
 
The following lemma is standard.

\begin{lemma} 
a. If $\gF$ is a quasi-coherent sheaf on $\pv$, then
$\Gamma_*(G \times \pv, p^*(\gF)) = k[G] \te \Gamma_*(\pv, \gF)$.

b. If $M$ is an $\sv$-module, then $p^*(\tilov{M}) =
(k[G] \te M)^\til$.
\end{lemma}

\begin{proposition} \label{7gekv-svg}
a. Let $\gF$ be a $G$-equivariant quasi-coherent sheaf on $\pv$. Then
$\Gamma_*(\pv,\gF)$ is an $\svgmod$.

b. Let $M$ be an $\svgmod$. Then $\tilov{M}$ is a $G$-equivariant 
quasi-coherent sheaf on $\pv$. 
\end{proposition}

\begin{proof}
For a coherent sheaf $\gF$ on $\pv$ denote $\gstar(\pv, \gF)$ for 
short as $\gstar(\gF)$.

Giving an isomorphism
$ \phi : q^*(\gF) \pil p^*(\gF)$ is equivalent to giving an isomorphism
(where $i_{\sv}$ is the homomorphism given by $i_{\sv}(s) = 1 \te s$)
\begin{equation} 
( k[G] \te \sv){}^{\delsv} \te_{\sv} \Gamma_*(\gF) 
\mto{\Gamma_*(\phi)} (k[G] \te \sv){}^{i_{\sv}} \te_{\sv} \Gamma_*(\gF) \label{7gf}
\end{equation}
of $k[G] \te \sv$-modules (with action on the first factors).

Now the map $\delsv$ induces a map of $\sv$-modules (with $\sv$-action on
the right module via $\delsv$)
\[ \sv \te_{\sv} \Gamma_*(\gF) \pil (k[G] \te \sv)^{\delsv} \te_{\sv} 
\Gamma_*(\gF). \]
Composing with $\Gamma_*(\phi)$ this gives a map 
\[ \gaf \mto{\Delta_{\gaf}} k[G] \te \gaf \]
which is an $\sv$-module map where the $\sv$-module structure on
$k[G] \te \gaf$ is given via $\delsv$.

We next need to show that this map $\Delta_{\gaf}$ makes $\gaf$ into a
$k[G]$-comodule. I.e. that the following diagram commutes.
 \[ \begin{CD} 
 \gaf @>{\Delta_{\gaf}}>> k[G] \te\gaf \\
 @V{\Delta_{\gaf}}VV  @VV{\Delta_{k[G]} \te \id_{\gaf}}V \\
 k[G] \te \gaf @>>{ \id_{k[G]} \te \Delta_{\gaf}}> k[G] \te k[G] \te \gaf.
 \end{CD} \]
This may be verified by translating the cocycle condition
(\ref{7kosykel}) to statements for $k[G] \te k[G] \te \sv$-modules.


%
 
 \medskip
 
 b. The 
automorphism $\eta$ of $G \times \pv$ is induced from the automorphism
$\strdelsv$ which is the composition
\begin{equation} k[G] \te \sv  \mto{\id \te \Delta_{\sv}}
k[G] \te k[G] \te \sv \mto{m_{k[G]} \te \id_{\sv}} k[G] \te \sv. 
\label{7strdelsv} \end{equation}

Suppose now $M$ is an $\svgmod$. Then the isomorphism
 \[ u_M \, : \, k[G] \te M \mto{ \id \te \Delta_M} k[G] \te k[G] \te M
 \mto{m_{k[G]} \te \id_M} k[G] \te M \]
 is a map of $k[G] \te \sv$-modules, where the $k[G] \te \sv$-module
 structure on the right is given via the map $\strdelsv$ in 
(\ref{7strdelsv}). The composition, which we denote by $\theta$,
 \begin{eqnarray} \label{7theta} && (k[G] \te \sv){}^{\Delta_{\sv}}
 \te_{\sv} M \\
& \mto{ (\strdelsv \te \id_M)^{-1} }& 
 (k[G] \te \sv){}^{i_{\sv}} \te_{\sv} M \notag \\
&\ispil & k[G] \te M
 \mto{u_M} k[G] \te M  \notag
 \end{eqnarray}
 now gives an isomorphism of $k[G] \te \sv$-modules (the action on the first
 module is given by acting naturally on the first factor). 
We note that this map
 is determined by
 \[ 1 \te 1 \te m  \mapsto \Delta_M(m) = \sum_{(m)} m_{k[G]} \te m_M . \]
 Sheafifying (\ref{7theta}) we get an isomorphism 
 \[ q^*(\tilov{M}) \mto{\tilov{\theta}} p^*(\tilov{M}) . \]
 We next need to show that this morphism fulfills the cocycle condition
 (\ref{7kosykel}). For this it will be enough to show that the
 following map $\alpha$ of $k[G] \te k[G] \te \sv$-modules
 \begin{eqnarray*}
 & (k[G] \te k[G]){}^{\Delta_{k[G]}} \te_{k[G]} (k[G] \te \sv){}^{\delsv}
 \te_{\sv} M  \\
 \mto{(k[G] \te k[G]){}^{\delkg} \te_{k[G]} \theta} &
 (k[G] \te k[G]){}^{\delkg} \te_{k[G]} (k[G] \te M) = k[G] \te k[G] \te M 
 \end{eqnarray*}
 coincides with the map $\beta $ of $k[G] \te k[G] \te \sv$-modules which
 is the composition
 \begin{eqnarray*}
  &&(k[G] \te \sv)^{\delsv} \te_{\sv} (k[G] \te \sv)^{\delsv} \te_{\sv}
 M  \\
  &\mto{(k[G] \te \sv){}^{\delsv} \te_{\sv} \theta } &
 (k[G] \te \sv)^{\delsv} \te_{\sv} (k[G] \te \sv)^{i_{\sv}} \te_{\sv} M \\
 & \isopil & k[G] \te (k[G] \te \sv)^{\delsv} \te_{\sv} M \\
 & \mto{k[G] \te \theta} & k[G] \te k[G] \te M. 
 \end{eqnarray*}
 
 But the map $\alpha$ sends
 \begin{equation} 1_{k[G]} \te 1_{k[G]} \te 1_{k[G]} \te 1_{\sv} \te m 
 \mapsto (\delkg \te \id_M) \circ \Delta_M (m) \label{7alsend}
 \end{equation}
 and the map $\beta$ sends
 \begin{eqnarray} 1_{k[G]} \te 1_{k[G]} \te 1_{k[G]} \te 1_{\sv} \te m
 & \mapsto & \sum_{(m)} m_{k[G]} \te 1_{k[G]} \te 1_{\sv} \te m_{\sv} 
\notag \\
 & \mapsto & \sum_{(m)} m_{k[G]} \te \Delta_M(m_{\sv}).  \label{7betasend}
 \end{eqnarray}
 Since $M$ is a $k[G]$-comodule the expressions (\ref{7alsend}) and 
 (\ref{7betasend}) are equal.
 \end{proof}
 
 As a consequence of Proposition \ref{7gekv-svg} we get that there are functors
 \[  \vmod{\sv,G}  \bihom{\til}{\gags(\pv,-)} \gqc . \]
For short when $\gF$ is in $\gqc$,
we denote $\gags(\pv, \gF)$ as $\gags(\gF)$. 

 \begin{lemma} The functor $\til$ is left adjoint to the functor
 $\gags$, i.e. for $M$ in $\modv{\sv, G}$ and $\gF$ in $\gqc$ there
 is a natural isomorphism
 \[ \Hom_{\gqc}(\tilov{M}, \gF) \iso
\Hom_{\modv{\sv, G}}(M, \gags(\gF)). \]
 \end{lemma}
 
 \begin{proof} 
 Given $M \pil \gags(\gF), $ an $\svgmod$ map, we get by the
functoriality of  (\ref{7theta})
 a commutative diagram of $k[G] \te \sv$-modules
 \[ \begin{CD}
  (k[G] \te \sv){}^{\delsv} \te_{\sv} M  @>>> 
(k[G] \te \sv){}^{\delsv} \te_{\sv}
  \gags(\gF) \\
  @VVV @VVV \\
  k[G] \te M @>>> k[G] \te \gags(\gF)
  \end{CD}  \]
 Sheafifying this and using that $\tilov{\Gamma}_{G,*}(\gF) \iso \gF$ 
we get a morphism in $\gqc$.
 Conversely, given a morphism $\phi \, : \, \tilov{M} \pil \gF$ in 
 $\gqc$. This gives a commutative diagram of quasi-coherent sheaves on 
 $G \times \pv$. 
 \[ \begin{CD}
 q^*(\tilov{M}) @>{q^*(\phi)}>>  q^*(\gF) \\
 @VVV @VVV \\
 p^*(\tilov{M}) @>{p^*(\phi)}>> p^*(\gF).
 \end{CD} \]
 Observe that $p^*(\tilov{M})$ is the sheafification of $k[G] \te M$ and
 $\Gamma_*(G \times \pv, p^*(\gF)) = k[G] \te \gstar(\gF)$. Also 
 $q^*(\tilov{M})$ is the sheafification of $(k[G] \te \sv){}^{\delsv} 
 \te_{\sv} M$ and 
 \[ \gstar (G \times \pv, q^*(\gF)) = (k[G] \te \sv){}^{\delsv}
 \te_{\sv} \gaf. \]

 Now we are going to apply the adjunction (\ref{2adj}) noting that it is valid
 on $\text{Spec} A \times_{\text{Spec} k} \pv$ for any $k$-algebra $A$.
 Thus we get a diagram
 \[ \begin{CD}
 (k[G] \te \sv){}^{\delsv} \te_{\sv} M @>>> (k[G] \te \sv){}^{\delsv} \gaf \\
 @VVV @VVV \\
 k[G] \te M @>>> k[G] \te \gaf. 
 \end{CD} \]
 
 Combined with the commutative diagram of $\sv$-modules (where the $\sv$-module
 structure on the lower row is given via $\delsv$)
 \[ \begin{CD}
 M @>>> \gaf \\
 @VVV @VVV  \\
 (k[G] \te \sv){}^{\delsv} \te_{\sv} M @>>> (k[G] \te \sv){}^{\delsv} \te_{\sv}
 \gaf 
 \end{CD} \]
 this gives an $\svgmod$ map $M \pil \gags(\gF)$.
 \end{proof}
 
 \subsection{ Injectives in $\gqc$}
 
 In order to develop a $G$-equivariant version of our theory we need to 
 establish that the category $\gqc$ has enough injectives, which are 
 {\it also} acyclic considered as objects in $\qc$ for the functor
 \[ \gstar \, : \, \qc \pil \modv{\sv}. \]
 Then these injectives can be used to compute the derived functors of
 both $\gags$ and $\gstar$. We will construct these injectives in $\gqc$ as
 $p_* q^*(\gI)$ where $\gI$ is an injective in $\qc$. First we show the 
 following.
 
 \begin{proposition} Let $\gF$ be in $\qc$. Then $p_*   q^* (\gF)$ is in 
 $\gqc$.
 \end{proposition}
 
 \begin{proof}
 
 First $p_*   q^*(\gF)$ is quasi-coherent by \cite[II.5.8]{Ha}. Let $M = 
 \gaf$. Since $q^*(\gF)$ is the sheafification of 
the $k[G] \te \sv$-module
\[ (k[G] \te \sv){}^{\delsv} 
 \te_{\sv} M,\] 
then $p_* q^*(\gF)$ is the sheafification
 of the $\sv$-module 
\[ (k[G] \te \sv){}^{\delsv} \te_{\sv} M \] 
 where $\sv$ acts by multiplication on the second factor.
 
 By Proposition \ref{7gekv-svg} b. we thus need to show that this module 
 is naturally a 
 $\svgmod$. Let $\sigma : k[G] \pil k[G]$ be the antipode, i.e. the
 algebra homomorphism corresponding to $G \pil G$ given by
 $g \mapsto g^{-1}$. We claim that there is a map 
 \begin{equation} \label{deltaM}
(k[G] \te \sv){}^{\delsv} \te_{\sv} M \mto{\Delta} k[G] \te 
 (k[G] \te \sv){}^{\delsv} \te_{\sv} M \end{equation}
 making the first module into a $k[G]$-comodule. The map $\Delta$ is given by
 \[ \mu \te s \te m \mapsto \sum_{(\mu) (s)} \sigma(\mu^{\prime \prime})
 s_{k[G]} \te \mu^{\prime} \te s_{\sv} \te m. \]
 
 To check that the map is well defined, one verifies
(using that $k[G]$ is a Hopf algebra)
 \[ \Delta(1 \te 1 \te sm) = \Delta( \sum_{(s)} s_{k[G]} \te s_{\sv} \te m). \]
 It is also easily seen that $\Delta$ is an $\sv$-module map 
 (with $\sv$ acting via
 $\Delta_{\sv}$ on the second module and then on the first and third
factor). Finally it is easily seen that 
$\Delta$ is a $k[G]$-comodule map. 
 Hence by Proposition \ref{7gekv-svg}, $p_*q^*(\gF)$, 
 which is the sheafification of 
 $(k[G] \te \sv){}^{\delsv} \te_{\sv} M$ is a $G$-equivariant 
 quasi-coherent sheaf.
 \end{proof}
 
 \begin{proposition} The forgetful functor
 \[ o \, : \, \gqc \lpil \qc  \]
 is left adjoint to the functor
 \[ p_* q^* \, : \, \qc \lpil \gqc. \]
 \end{proposition}
 
 \begin{proof}
 Let $\gA$ be in $\gqc$ and let $\gF$ be in $\qc$.
 Given a morphism $o(\gA) \pil \gF$, we get a $G$-equivariant morphism
 \begin{equation}
 p_*   q^*(o(\gA)) \pil p_* q^*(\gF). \label{7oaf}  \end{equation}
 Now given any morphism of schemes $f : X \pil Y$, the functor 
 $f^*$ is left adjoint to the functor
 $f_*$. Hence the isomorphism $p^*(\gA) \pil q^*(\gA)$ 
corresponds to a morphism
 $\gA \pil p_*  q^* (o(\gA))$ which is $G$-equivariant (as can be checked). 
 Composing with (\ref{7oaf}) gives us the $G$-equivariant morphism 
 $\gA \pil p_* q^*(\gF)$. 
 
 Conversely given a $G$-equivariant morphism $\gA \pil p_* q^*(\gF)$
 this corresponds to a morphism $p^*(o(\gA)) \pil q^*(\gF)$. 
 Taking the fiber at $1$ in $G$ gives a morphism
 $o(\gA) \pil \gF$. 
 
 Now starting out from $\alpha : o(\gA) \pil \gF$, we get a map 
 $\gA \pil p_* q^*(\gF)$ and again a map $o(\gA) \pil \gF$ which is 
 easily seen to be $\alpha$.
 
 Conversely, starting with a $G$-equivariant map 
$\beta : \gA \pil p_* q^*(\gF)$
 we get $o(\gA) \pil \gF$ and then again by our construction a map
 $\gamma : \gA \pil p_* q^*(\gF)$. We must show that $\gamma = \beta$. 
 So let $\delta = \beta - \gamma$. Then $\delta : \gA \pil p_* q^*(\gF)$ is a
 $G$-equivariant map such that the associated $o(\gA) \pil \gF$ is zero.
 We must show that then $\delta = 0$. 
 
 So let $M = \Gamma_*(\gF)$ and $A = \gags(\gA)$. Then 
 \[ \gstar(p_* q^*(\gF)) = (k[G] \te \sv){}^{\delsv} \te_{\sv} M \]
 with $\sv$ acting
 by multiplication on the second factor. The map $\delta$ corresponds to a
 map of $\svgmod$s
 \[ \gstar(\delta) : A \lpil (k[G] \te \sv){}^{\delsv} \te_{\sv} M \]
 such that composing with the counit $\eta : k[G] \pil k$ gives the 
 zero map $A \pil M$. Since $\gstar(\delta)$ is a $k[G]$-comodule map there is
 a commutative diagram
 \begin{equation} \begin{CD}
 A @>{\gstar(\delta)}>> (k[G] \te \sv){}^{\delsv} \te_{\sv} M   \\
 @VVV  @V{\Delta}VV \\
 k[G] \te A @>{k[G] \te \gstar(\delta)}>> 
 k[G] \te (k[G] \te \sv){}^{\delsv} \te_{\sv} M \\
 @. @V{\id \te \eta \te \id \te \id}VV \\
 @.  (k[G] \te \sv){}^{i_{\sv}} \te_{\sv} M
 \end{CD} \label{7gdel} \end{equation}
 So let $\sum \mu \te s \te m$ be in the image of $\gstar(\delta)$. Then
 \begin{equation} \Delta(\sum \mu \te s \te m) = \sum \sum_{(\mu) (s)} 
 \sigma(\mu^{\prime \prime})s_{k[G]} \te \mu^{\prime} \te s_{\sv} \te m.
 \label{7gdelelm} \end{equation}
 The composition of the maps $k[G] \te \gstar(\delta)$ and
 $\id \te \eta \te \id \te \id$ is zero since this composition is just
$k[G]$ tensor the map $A \pil M$ which is zero.
  Hence (\ref{7gdelelm}) maps to zero
 by $\id \te \eta \te \id \te \id$
\[ 0 = \sum \sum_{(\mu) (s)} \eta(\mu^{\prime}) \sigma(\mu^{\prime \prime})
s_{k[G]} \te s_{\sv} \te m. \]
Since
\[ \sum_{(\mu)} \eta(\mu^\prime) \mu^{\prime\prime} = \mu, \] 
we get that this is
\[ 0 = \sum \sum_{(s)} \sigma(\mu) s_{k[G]} \te s_{\sv} \te m \in
(k[G] \te \sv)^{i_{\sv}} \te_{\sv} M. \]
Composing with the map
\[ (k[G] \te \sv)^{i_{\sv}} \te_{\sv} M
\mto{\sigma \te \id \te \id} 
(k[G] \te \sv)^{i_{\sv}} \te_{\sv} M \]
we get that
\begin{equation} \sum \sum_{(s)} \mu \sigma(s_{k[G]}) \te s_{\sv} \te m = 0.
\label{7sum0}  \end{equation}
Tensoring the map $\strdelsv$ with $M$ we get a map
\[ (k[G] \te \sv){}^{i_{\sv}} \te_{\sv} M 
\lpil (k[G] \te \sv){}^{\delsv} \te_{\sv} M \]
which maps (\ref{7sum0}) to
\[ \sum \sum_{(s)} \mu \sigma(s^{\prime}_{k[G]})s^{\prime\prime}_{k[G]}
\te s_{\sv} \te m \label{7sum00} \]
where we have used that the map
\[ (k[G] \te \Delta_{\sv}) \circ \Delta_{\sv} \, : \, 
\sv \pil k[G] \te k[G] \te \sv \]
maps $s$ to $\sum_{(s)} s^{\prime}_{k[G]} \te s^{\prime\prime}_{k[G]}
\te s_{\sv}$.
But 
\[ \sum_{(s_{k[G]})} \sigma(s^\prime_{k[G]}) s^{\prime \prime}_{k[G]}
= \eta(s_{k[G]}) \]
and 
\[ \sum_{(s)} \eta(s_{k[G]}) s_{\sv} = s \]
and so (\ref{7sum00}) is
\[ 0 = \sum \mu \te s \te m \in (k[G] \te \sv)^{\delsv} \te_{\sv} M. \]
This gives $\gstar(\delta) = 0$. Hence the adjunction is proven.
\end{proof}

\begin{corollary} \label{7nokinj}
If $\gI$ is injective in $\qc$, then $p_* q^*(\gF)$
is injective in $\gqc$.
\end{corollary}

\begin{proof}
This is just the general fact that a right adjoint to an exact functor
between abelian categories takes injectives to injectives.
\end{proof}

\begin{proposition} \label{7ogasyk}
Let $\gI$ be an injective object of $\qc$. Then $p_* q^*(\gI)$, considered as
an object in $\qc$, is acyclic for the functor
$\gstar : \qc \pil \modv{\sv}$.
\end{proposition}

\begin{proof}
Since $p_* q^*(\gI(n)) = (p_* q^*(\gI))(n)$ by the projection formula,
it will be enough to show that the right derived functors
of $\Gamma(\pv, -)$ vanish.
The functor $\Gamma(G \times \pv, -)$ is the composition of the functors
$p_*$ and $\Gamma(\pv, -)$. 
Note first that $p_*$ takes injectives to objects acyclic for 
$\Gamma(\pv,-)$. This is because if $\gJ$ is an injective
quasi-coherent sheaf on $G \times \pv$, then $\gJ$ is flasque 
\cite[III.2]{Ha} and hence $p_*(\gJ)$ is also flasque and hence
acyclic for $\Gamma(\pv,-)$. Also $p_*(\gJ)$ is quasi-coherent
since $G \times \pv$ is Noetherian. 

The composition of $p_*$ and $\Gamma(\pv,-)$ therefore gives a
Grothendieck spectral sequence for the sheaf $q^*(\gI)$. 
\[ H^r (\pv, {\bf R}^s p_*(q^*(\gI))) \Rightarrow 
H^{r+s} (G \times \pv, q^*(\gI)). \]
Now since $p$ is an affine morphism, ${\bf R}^s p_* (q^*(\gI)) = 0$ for
$s > 0$. Thus
\[ H^r(\pv, p_*(q^*(\gI))) = H^r(G \times \pv, q^*(\gI)). \]
We show that the latter is zero for $r>0$ and this will finish the
proof.
But $q^*(\gI) = \eta^* p^*(\gI)$ and $\eta$ is an isomorphism. 
Thus it is enough to show that
\[ H^r( G \times \pv, p^*(\gI)) = 0 \]
for $r > 0$.
But the same argument as above also gives an isomorphism
\begin{equation} H^r(\pv, p_*(p^*(\gI))) = H^r( G \times \pv, p^*(\gI)).
\label{7gro} \end{equation}
By the projection formula $p_*(p^*(\gI)) = \gI \te k[G]$. Since $\pv$ is
a Noetherian scheme, a direct sum of injectives is injective, and so 
(\ref{7gro}) vanish for $r > 0$.
\end{proof}

\section{$G$-equivariant versions}

Suppose the vector space $W$ comes equipped with the action of an
algebraic group $G$. Remark \ref{3dtilG} suggests the following natural way of
constructing objects in $\dbco$. Let $A$ be a representation of $G$ and 
let $B$ be a quotient representation of $W \te A$. This gives
a morphism $\lv(-1) \te A \mto{d} \lv \te B$ and by
Remark \ref{3dtilG} we thus get an object $\gG$ in $\dbco$. Now the map
$d$ is a map of $\lvd,G$-modules. Therefore $\gG$ should be a 
$G$-equivariant complex of coherent sheaves.

In this section we develop the theory for such conclusions. The 
main results are generalizations of Theorems \ref{3ekv} and 
\ref{3ekspl} to the
$G$-equivariant case. In order to have good generalizations we shall
assume that {\it the category of $G$-modules is semi-simple}.
This holds for instance if char $k = 0$ and $G$ is a finite or a 
semi-simple group.

\subsection{G-equivariant versions. }
The category $\gqc$ has full subcategories
\[ \gvb \sus \gco \sus \gqc \]
consisting of $G$-equivariant locally free sheaves and $G$-equivariant coherent
sheaves respectively. We get derived categories
\[ D^b(\gvb),\; D^b(\gco),\; D^b_{\gco}(\gqc), \] 
and
\[ D_{b, \gco}(\gqc). \]
Proposition \ref{2Pderekv}
easily generalizes so all these categories are equivalent.

In the previous section we showed, Corollary \ref{7nokinj}, that the category
$\gqc$ has enough injectives, so we get a right derived functor
\[ \rgags \, : \, D_{b, \gco}(\gqc) \lpil D_G(\sv). \]
We can compose this with the Koszul functor
\[ G_{\sv} \, : \, D_G(\sv) \lpil D^R_G(\lvd). \]

\begin{proposition} Let $\gQ$ be in $D_{b, \gco}(\gqc)$. Then
$G_{\sv} \circ \rgags (\gQ)$ is acyclic.
\end{proposition}

\begin{proof}
Let $\gI$ be an injective resolution of $\gQ$ in $\gqc$
such that $\gI^i$ is acyclic for the functor 
$\gstar : \qc \pil \vmod{\sv}$. 
Such a resolution exists by Corollary \ref{7nokinj} and 
Proposition \ref{7ogasyk}.

Then $\gags(\gI)$ considered as an object in $D(\sv)$ is isomorphic
to  $\rgstar(\gQ)$. By Proposition \ref{2acext} we then get that
$G_{\sv} \circ \gags(\gI)$ is acycylic, which proves the proposition.
\end{proof}

Let $K_G^\circ(\lvdcf)$ be the full subcategory of
$K_G(\lvdcf)$ consisting of acyclic complexes whose components have finite
corank. We now get the main theorem of this section which generalizes
Theorem \ref{3ekv}.

\begin{theorem}
Let $G$ be a linear algebraic group such that the category of 
$G$-modules is semi-simple. Let $V$ be a finite dimensional $G$-module.

There is a functor
\[ \gsmin \circ  \trgags \, : \, D_{b, \gco}(\gqc) \lpil 
K^\circ_{G}(\lvdcf) \]
and a functor
\[ \til \circ \fl \, : \, K^\circ_G(\lvdcf) \lpil D_{b, \gco}(\gqc). \]
These functors give an adjoint equivalence of triangulated categories,
with $\til \circ \fl$ left adjoint.
\end{theorem}

\begin{proof}
This is analogous to the proof of Theorem \ref{3ekv}.
\end{proof}

Let 
\[ o_1 \, : \, \modv{\sv, G} \lpil \modv{\sv} \]
and 
\[ o_2 \, : \, \gco \lpil \coh \]
be the forgetful functors. If $\gF$ is in $\gco$ we see by 
Corollary \ref{7nokinj} and Proposition \ref{7ogasyk} that
\[ o_1(H^p( \rgags (\gF))) = H^p_* o_2(\gF). \]
Hence when $\gF$ is in $\gqc$ then $H^p_* o_2(\gF)$ comes with the
structure of an $\svgmod$
which we still denote as $H^p_* \gF$. Then the analog of Theorem \ref{3ekspl} 
holds.

\begin{proposition}
Let $\gF$ be in $\gco$. Then $\gsmin \circ  \trgags (\gF)$ is a minimal
complex with
\[ (\gsmin \circ  \trgags (\gF))^p = \oplus_{r = 0}^{\pdim}
\lv(p-r) \te H^r \gF(p-r). \]
\end{proposition}

\section{Examples}

In this section we shall illustrate the theory presented in this paper
in some examples. The main example will be the construction of the
Horrocks-Mumford bundle. We will show how the construction of this
bundle becomes very natural. One does almost not have to use any
cleverness in constructing it.

\subsection{ Constructing the Horrocks-Mumford bundle}

Let $\dim_k V = 5$. The Horrocks-Mumford bundle on $\pv$ is a rank $2$ 
bundle $\gE$
with Chern classes $c_1(\gE) = -1 $ and $c_2(\gE) = 4$. From this one 
easily calculates that the Hilbert polynomial of $\gE$ is 

\[ \chi \gE(n) = 2 \binom{n+4}{4} - \binom{n+3}{3} - 4 \binom{n+2}{2} 
                - 2 \binom{n+1}{1}. \]
From this again we calculate the following values
\vskip 3mm
\hskip 1.5cm \begin{tabular}{c | c | c | c | c | c | c | c}
$n$  & -5 & -4 & -3 & -2 & -1 & 0 & 1 \\ \hline
$\chi \gE(n)$ & -10 & -5 & 0 & 2 & 0 & -5 & -10 
\end{tabular} 
\vskip 3mm
Now the Hilbert polynomial is
\[ \chi \gE(n) = h^0 \gE(n) - h^1 \gE(n) + h^2 \gE(n) - h^3 \gE(n) 
+ h^4 \gE(n). \]
If the cohomology $H^i \gE(n)$ behaves as nicely as possible, one would
expect for a given $n$ that at most one $H^i \gE(n) \neq 0$ for
$0 \leq i \leq 4$. If things are as nice as possible one would therefore 
expect the following cohomology table (a ``$\cdot$'' indicates a zero)
\vskip 3mm
\hskip 1.5cm \begin{tabular}{c | c | c | c | c | c | c | c}
 $n$  & -5 & -4 & -3 & -2 & -1 & 0 & 1 \\ \hline
$h^0 \gE(n)$ &  $\cdot$ & $\cdot$ & $\cdot$ & $\cdot$ & $\cdot$ & $\cdot$ & $\cdot$ \\ \hline
$h^1 \gE(n)$ &  $\cdot$ & $\cdot$ & $\cdot$ & $\cdot$ & 0  & 5 & 10 \\ \hline
$h^2 \gE(n)$ &  $\cdot$ & $\cdot$ & $\cdot$ & 2 & $\cdot$ & $\cdot$ & $\cdot$ \\ \hline
$h^3 \gE(n)$ &  10 & 5 & 0 & $\cdot$ & $\cdot$ & $\cdot$ & $\cdot$ \\ \hline
$h^4 \gE(n)$ &  $\cdot$ & $\cdot$ & $\cdot$ & $\cdot$ & $\cdot$ & $\cdot$ & $\cdot$ \\
\end{tabular}
\vskip 3mm
If this holds the exterior complex $\ytrev(\gE)$ would have components in
degrees $-1, 0$, and $1$ as follows
\[ \lv(-4)^5 \mto{d^{-1}} \lv(-2)^2 \mto{d^0} \lv^5. \]
Now $d^0$ is given by a $5 \times 2$ matrix of quadratic exterior forms
on $V$. A tempting guess is that the columns of this matrix are cyclic
permutations of exterior forms $e_i \wedge e_j$ where $e_0, e_1, e_2, 
e_3, e_4$ is a basis for $V$. So let (check that $d^0 \circ d^{-1} = 0$)

\[ d^0 = \left( \begin{array}{cc} 
   e_0 \wedge e_1 & e_2 \wedge e_4 \\
   e_1 \wedge e_2 & e_3 \wedge e_0 \\
   e_2 \wedge e_3 & e_4 \wedge e_1 \\
   e_3 \wedge e_4 & e_0 \wedge e_2 \\
   e_4 \wedge e_0 & e_1 \wedge e_3
   \end{array} \right ),
  \quad d^{-1} = \left( \begin{array}{cc}
   e_2 \wedge e_4 & e_1 \wedge e_0 \\
 e_3 \wedge e_0 & e_2 \wedge e_1 \\
 e_4 \wedge e_1 & e_3 \wedge e_2 \\
   e_0 \wedge e_2 & e_4 \wedge e_3 \\
  e_1 \wedge e_3 & e_0 \wedge e_4 
\end{array} \right )^T. \]
  
\begin{lemma} \label{9E} \hfill

a. The complex 
\[ E \, : \, \lv(-4)^5 \mto{d^{-1}} \lv(-2)^2 \mto{d^0} \lv^5 \]
is exact.

b. We have the following table
\vskip 2mm {\rm
\hskip 1.5cm \begin{tabular}{c | c | c | c | c | c | c | c}
$n$  & $\leq$ -2  & -1 & 0 & 1 & 2 & $\geq$ 3 \\ \hline
$\dim_k (\ker d^0)_n$ & 0 & 5 & 15 & 10 & 2 & 0
\end{tabular}. }
\vskip 2mm

\end{lemma}

\begin{proof} a. 
It is clear that the complex is exact in degrees $\leq -3$. It is also
easily seen that it is exact in degree $-2$, since there are no relations
of $d^0$ of linear forms in $V$.

So consider degree $-1$. It is straight forward to check that 
the image of $(d^0)_{-1}$ has dimension $15$ and so the sequence is exact in degree
$-1$. 

Now taking the graded dual of the complex and using the identification
$(\lv)^{\gd} \iso \lvd \iso \lv(-5)$ we see that the complex becomes
isomorphic to the original complex twisted with $-1$. That the complex
is exact in degrees $\geq 0$ follows then from this fact. 

b. This is clear from the considerations above.

\end{proof}

By taking a free resolution of the kernel of $d^{-1}$ and a cofree
resolution of the cokernel of $d^0$ we then get a complex $I$ in
$K^{\circ}(\lvdcf)$, which corresponds to an object $\gF$ in 
$D^b(\coh)$. By Lemma \ref{9E} b. and Theorem \ref{5ker} 
the Hilbert polynomial of $\gF$ is
\begin{eqnarray*} \chi \gF(n) & = & 2 \binom{n+6}{4} - 10 \binom{n+5}{4} 
+ 15 \binom{n+4} {4} - 5 \binom {n+3}{4} \\
  & = &  2 \binom{n+4}{4} - \binom{n+3}{3} - 4 \binom{n+2}{2} 
                - 2 \binom{n+1}{1} 
\end{eqnarray*}
which coincides with the Hilbert polynomial of the Horrocks-Mumford
bundle. We next need to show that $\gF$ is isomorphic to a rank $2$ vector
bundle. For this we invoke the theory of Section 6.

\begin{lemma} Let $U \sus W$ have codimension $1$. Consider the
complex $\Hom_{\lvd}(\lwd, E)$
\begin{equation} 
  \lw(-4)^5 \mto{\Hom_{\lvd}(\lwd, d^{-1})} \lw(-2)^2
            \mto{\Hom_{\lvd}(\lwd, d^0)} \lw^5. \label{9lw}
\end{equation}
Then the dimension of
$H^0 \Hom_{\lvd}(\lwd, E)_q$ is $2$ for $ q = 0$, and 
zero for $q \neq 0$ regardless of $U$. \label{9WV}
\end{lemma}

\begin{proof} $U$ corresponds to a line $(u)$ in $V$. The matrix of 
$\Hom_{\lvd}(\lwd, d^0)$ is obtained from the matrix of $d^0$ by letting
the entries map to quadratic forms in $V/(u)$. Let $\overline{e_i}$
be the image of $e_i$ in $V/(u)$. Clearly in degrees $\leq -2$ the 
cohomology of (\ref{9lw}) vanishes. 

\medskip
Consider the complex in degree $-1$. We must prove that
the differential 
$\Hom_{\lvd}(\lwd, d^0)$ has no linear relations in $V/(u)$. 
By a suitable permutation of $\ove{0}, \ldots, \ove{4}$ we may assume that a 
relation between the $\ove{0}, \ldots, \ove{4}$ involves only
$\ove{i}$ for $i = 0, \ldots, r$. We may find such a permutation so that
the new matrix $d^{0\prime}$ obtained from $d^0$ also may be obtained
from $d^0$ by permuting its rows and columns.
Now for each $r = 0, \ldots, 4$ it is then easily checked that
there is no linear relation of 
$\Hom_{\lvd}(\lwd, d^0)$.
%
Hence
the complex is exact in degree $-1$.

  So consider the complex in degree $0$. We shall show that \linebreak
$\Hom_{\lvd}(\lwd, d^0)$ is surjective in degree $0$. 
First a piece of notation. If $i_0 < i_1 < i_2 < i_3$ and 
$i$ is the element in  
$\{ 0,1,2,3,4 \}$ not in $\{ i_0, i_1, i_2, i_3 \}$,
denote $\ove{i_0} \wedge \ove{i_1} \wedge \ove{i_2} \wedge \ove{i_3}$
by $\cove{i}$.

Assume that $\{ \ove{0}, \ove{1}, \ove{2}, \ove{3} \}$ is independent.
The images of 
\[ \begin{pmatrix} \ove{0} \wedge \ove{3} \\ 0 \end{pmatrix},
   \begin{pmatrix} \ove{2} \wedge \ove{3} \\ 0 \end{pmatrix},
   \begin{pmatrix} \ove{0} \wedge \ove{1} \\ 0 \end{pmatrix},
   \begin{pmatrix}  0 \\ - \ove{1} \wedge \ove{3} \end{pmatrix},
   \begin{pmatrix}  0 \\ - \ove{0} \wedge \ove{2} \end{pmatrix}
\]
by $\Hom_{\lvd}(\lwd, d^0)$ are the following :
\[ \begin{pmatrix} 0 \\ \cove{4} \\ 0 \\ 0 \\ 0  \end{pmatrix},
   \begin{pmatrix} \cove{4} \\0 \\0 \\0 \\ -\cove{1} \end{pmatrix},
   \begin{pmatrix} 0 \\ 0 \\ \cove{4} \\ \cove{2} \\ 0 
   \end{pmatrix},
   \begin{pmatrix}  \cove{0} \\ 0 \\ 0 \\ \cove{4} \\ 0 
   \end{pmatrix},
   \begin{pmatrix}  0 \\ 0 \\ -\cove{3} \\ 0 \\ \cove{4} 
   \end{pmatrix}.  
\]   
Now if $\{ \ove{1}, \ove{2}, \ove{3}, \ove{4} \}$ is dependent, then
$\cove{0} = 0$, and we easily see that the image of 
$\Hom_{\lvd}(\lwd, d^0)$ is the whole of $(\lv^5)_0$. 

If $\{ \ove{1}, \ove{2}, \ove{3}, \ove{4} \}$ is independent, then the image
of 
\[ \begin{pmatrix} \ove{1} \wedge \ove{4} \\ 0 \end{pmatrix} \]
is the transpose of $\begin{pmatrix} 0 & 0 & \cove{0} & 0 & 0 
\end{pmatrix}$, which is nonzero. This also implies that the image of 
$\Hom_{\lvd}(\lwd, d^0)$ is the whole of $(\lv^5)_0$.

Now by taking the graded dual of the complex $\Hom_{\lvd}(\lwd, E)$ and
using the isomorphism $\lw^\gd \iso \lwd \iso \lw(-4)$, we get a complex
isomorphic to 
$\Hom_{\lvd}(\lwd, E)$. This gives that 
$\Hom_{\lvd}(\lwd, d^{-1})$ is injective in degree $0$ and 
$\Hom_{\lvd}(\lwd, E)$ is exact in degrees $\geq 1$.

In degree $0$ we see that the cohomology has dimension $2$.
\end{proof}

We then obtain the following.

\begin{theorem}
Consider the complex $E$ in Lemma \ref{9E}. By taking a free resolution
of $\ker d^{-1}$ and a cofree resolution of ${\rm \coker} d^0$ so we get
an object in $K^{\circ}( \lvdcf)$, this object corresponds to an 
object in $D^b(\coh)$ isomorphic to a rank two vector bundle
$\gE$ on $\pv = \pfour$ with Chern classes $c_1(\gE) = -1$ and 
$c_2(\gE) = 4$.
\end{theorem}

\begin{proof}
That $\gE$ is a rank $2$ bundle follows from Lemma \ref{9WV} in 
conjunction with Theorem \ref{6dtilF} and Corollary \ref{6bunt}. 
   By Lemma \ref{9E} b. and Theorem \ref{5ker} we find that the 
Hilbert polynomial
\[ \chi \gE (n) = 2 \binom{n+4}{4} - \binom{n+3}{3} - 4 \binom{n+2}{2} 
                - 2 \binom{n+1}{1} \]
and by standard computations this gives $c_1(\gE) = -1$ and $c_2(\gE) = 4$.
\end{proof}

\rem It is known \cite{HM} that the Horrocks-Mumford bundle is acted upon by
the Heisenberg group of order $125$, and even (assume $k = {\bf C}$)
by the the normalizer
$N$ of $H$, when $H$ is considered as a subgroup of $SL_5(k)$. 
The order of $N$ is $15 000$. The complex $E$ therefore is an 
$N$-equivariant complex
\[ \lv(-4) \te V_1 \mto{d^{-1}} \lv(-2) \te U \mto{d^0} \lv \te V_2 \]
where $V_1, V_2$, and $V$ are representations of $N$ of degree $5$ and
$U$ is a representation of $N$ of degree $2$.

\subsection{$GL(W)$-equivariant sheaves}

Let $\Omega_{\pv}$ be the sheaf of differentials on $\pv$. For
a partition 
\[ \lambda_1 \geq \lambda_2 \geq \cdots \geq \lambda_v \]
we get a Schur bundle $S_\lambda (\Omega_{\pv}(1))$ which is 
a $GL(W)$ equivariant bundle on $\pv$. Hence its cohomology groups
$H^p S_\lambda (\Omega_{\pv}(1))(q)$ are representations of $GL(W)$.

The irreducible representations of $GL(W)$ are parametrized
by partitions \cite[15.5]{FuHa} 
\[ \mu_0 \geq \mu_1 \geq \cdots \geq \mu_{v} \]
and we denote the corresponding representation as $S_\mu W$.

In \cite{Fl3}
all the cohomology groups of $S_\lambda (\Omega_{\pv}(1))$ are computed
and can be described as follows.

Set $\lambda_0 = +\infty$ and $\lambda_{v+1} = -\infty$.
Let $f_\lambda : \hele \pil \{0, \ldots, v\}$ be given by
$f_\lambda (a) = r$ if $\lambda_r > a \geq \lambda_{r+1}$ and
in this case let the partition $\lambda_a$ be given by
\[ \lambda_1 - 1 \geq \cdots \lambda_r - 1 \geq a \geq
\lambda_{r+1} \geq \cdots \geq \lambda_v. \]
Then the following holds.

\begin{theorem}
\[ H^r S_\lambda (\Omega_{\pv}(1))(a-r) = 
\begin{cases}
S_{\lambda_a} W, & r = f_{\lambda}(a) \\
0, & r \neq f_\lambda(a).
\end{cases} \]
\end{theorem}
In particular we see that all components of the exterior complex of 
$S_\lambda (\Omega_{\pv}(1))$ are of the simple form
\[ \lv(-r) \te S_\mu W. \]

\end{document}